\documentclass[a4paper,12pt,twoside]{article}
\usepackage[pdftex,colorlinks,bookmarksnumbered,bookmarksopen,citecolor=red,urlcolor=red]{hyperref}
\usepackage{graphicx}\usepackage{a4wide}
\usepackage{amssymb,amsmath,amsthm}
\usepackage{subfigure,color,colordvi,graphicx,xcolor}
\usepackage[only,llbracket,rrbracket,interleave]{stmaryrd}
\usepackage[superscript,sort,compress]{cite}
\usepackage{subfigure,color,colordvi,graphicx}
\newcommand{\bm}[1]{\text{\boldmath $#1$\unboldmath}}
\newcommand{\bddot}{\operatorname{\bm{:}}}
\newcommand{\vect}[1]{\mathbf{#1}}
\newcommand{\mat}[1]{\mathbf{#1}}
\newcommand{\Div}{{\bm{\nabla}{\cdot}}}
\newcommand{\grad}{\bm{\nabla}}
\newcommand{\sobo}[1][1]{\ensuremath{\mathcal{H}^{#1}}}
\newcommand{\eltwo}{\ensuremath{\mathcal{L}_2}}
\newcommand{\nen}  {\ensuremath{\texttt{n}_{\texttt{en}}}}
\newcommand{\nfn}  {\ensuremath{\texttt{n}_{\texttt{fn}}}}

\newcommand{\nsd}  {\ensuremath{\texttt{n}_{\texttt{sd}}}}
\newcommand{\numel}{\ensuremath{\texttt{n}_{\texttt{el}}}}

\newcommand{\bu}{\bm{u}}
\newcommand{\bn}{\bm{n}}
\newcommand{\bx}{\bm{x}}
\newcommand{\bL}{\bm{L}}
\newcommand{\bt}{\bm{t}}
\newcommand{\bs}{\bm{s}}
\newcommand{\hu}{\bm{\hat{u}}}
\newcommand{\bI}{\bm{I}}
\newcommand{\bxi}{\bm{\xi}}
\newcommand{\bet}{\bm{\eta}}
\newcommand{\blam}{\bm{\lambda}}

\newcommand{\nodalu}{\vect{u}}
\newcommand{\nodaluh}{\hat{\vect{u}}}
\newcommand{\nodalL}{\mat{L}}
\newcommand{\nodalp}{\text{p}}

\newcommand{\nodaluV}{\vect{u}}
\newcommand{\nodaluhV}{\hat{\vect{u}}}
\newcommand{\nodalLV}{\vect{L}}
\newcommand{\nodalpV}{\vect{p}}

\newcommand{\tu}{\bm{v}}
\newcommand{\tL}{\bm{G}}
\newcommand{\that}{\bm{\mu}}
\newcommand{\tdiv}{w}

\newcommand{\spP} {\mathcal{P}}

\newcommand{\scae}[3]{\bigl({#2},{#3}\bigr)_{#1}}
\newcommand{\scab}[3]{\bigl\langle{#2},{#3}\bigr\rangle_{#1}}
\newcommand{\bv}{\bm{v}}
\newcommand{\bU}{\bm{U}}
\newcommand{\bV}{\bm{V}}

\newcommand{\jump}[1]{\llbracket #1\rrbracket}

\newenvironment{keywords}{\begin{quote}\emph{\textbf{Keywords:}}}{\end{quote}}

\theoremstyle{definition}

\newtheorem{remark}{Remark}

\newcommand{\RS}[1]{#1}

\graphicspath{{figsSubmitted/}}

\begin{document}
%==========================================================================
\title{HDG-NEFEM with degree adaptivity for Stokes flows}

\author{Ruben Sevilla\\[-1ex]
             \small Zienkiewicz Centre for Computational Engineering, \\[-1ex]
             \small College of Engineering, Swansea University, Wales, UK \\[1em]
             Antonio Huerta\\[-1ex]
             \small Laboratori de C\`alcul Num\`eric (LaC\`aN), \\[-1ex]
             \small ETS de Ingenieros de Caminos, Canales y Puertos, \\[-1ex]
             \small Universitat Polit\`ecnica de Catalunya, Barcelona, Spain}

\date{December 1, 2017}
%________________________________________________________________________
\maketitle

%==========================================================================
\begin{abstract}
The NURBS-enhanced finite element method (NEFEM) combined with a hybridisable discontinuous Galerkin (HDG) approach is presented for the first time. The proposed technique completely eliminates the uncertainty induced by a polynomial approximation of curved boundaries that is common within an isoparametric approach and, compared to other DG methods, provides a significant reduction in number of degrees of freedom. In addition, by exploiting the ability of HDG to compute a postprocessed solution and by using a local a priori error estimate valid for elliptic problems, an inexpensive, reliable and computable error estimator is devised. The proposed methodology is used to solve Stokes flow problems using automatic degree adaptation. Particular attention is paid to the importance of an accurate boundary representation when changing the degree of approximation in curved elements. Several strategies are compared and the superiority and reliability of HDG-NEFEM with degree adaptation is illustrated.
\end{abstract}

%
%________________________________________________________________________
 
%________________________________________________________________________
%\subclass{65N08 \and 65N12 \and 65N22 \and 65N30}

%==========================================================================
\section{Introduction} \label{sc:intro}
%==========================================================================

Early work on mesh and degree adaptivity schemes for the finite element method~\cite{zienkiewicz1989effective,oden1989toward,johnson1992adaptive} already showed the advantages of adaptive schemes to achieve a required accuracy in an economic manner. The use of mesh adaptive methods is substantially more extended due to the popularity of low-order methods in the computational mechanics community. This is largely due, as discussed later, to the fact that mesh adaptation is easier to implement, compared to degree adaptivity, in standard finite element codes. But, with recent needs on high fidelity simulations for fluids and wave propagation phenomena~\cite{wang2013high,cohen2003higher,hesthaven2003high}, the interest in degree adaptive (or the combination of mesh and degree adaptivity) processes has increased~\cite{burbeau2005dynamic,giorgiani2013hybridizable,karban2013advanced,giorgiani2014hybridizable}. 

One of the main reasons for the increasing popularity of degree adaptive schemes in the last years is the rise of discontinuous Galerkin (DG) methods as a viable alternative for convection dominated flow and wave propagation problems~\cite{cockburn2001runge,nguyen2007rans,shu2016discontinuous,hesthaven2002nodal,dumbser2006arbitrary,de2008interior}. In a standard continuous Galerkin framework, the implementation of variable degree of approximation is cumbersome, whereas its application in a DG context is straightforward due to the weak imposition of the continuity of the solution by means of numerical fluxes. Despite traditional DG methods have not been able to consistently prove its superiority against low-order techniques traditionally employed in industry (e.g. finite volume methods), the recently proposed hybridisable DG (HDG)~\cite{Jay-CGL:09} has shown its superiority compared to traditional DG methods~\cite{cockburn2009hybridizable,Cockburn-KSC:11,AA-HARP:13}. The ability to substantially reduce the number of degrees of freedom combined with the possibility to obtain a post-processed solution that converges at a faster rate to the exact solution are the two main properties of HDG methods behind its superiority compared to other DG methods~\cite{nguyen2010hybridizable,cockburn2014devising,cockburn2009derivation,sevilla2016tutorial}. Moreover, this is achieved while preserving the well-known advantages of DG for stabilizing convection and circumventing the so-called Ladyzhenskaya-Babu{\v s}ka-Brezzi  (LBB) condition in the incompressible limit.

A key aspect in any adaptive scheme is the ability to devise cheap and reliable error measures for a given numerical solution in order to decide the regions where a more accurate solution is required~\cite{bangerth2013adaptive}. Error indicators and error estimators are typically employed to asses the error of a simulation with an adaptive framework~\cite{huerta1999adaptive}. Error indicators are computationally inexpensive but they are problem dependent whereas error estimators are considerably more expensive but more general~\cite{pares2006subdomain,diez2007equilibrated,pares2009exact}. A cheap, general and reliable error estimator was proposed in~\cite{giorgiani2013hybridizable,giorgiani2014hybridizable} by exploiting the ability of the HDG method to construct a post-processed solution, more accurate than the HDG solution.

One of the aspects that is normally ignored when devising degree adaptive schemes is the geometric representation of domains with curved boundaries. Despite it is now well known that a poor representation of the geometry can have an important effect on the results of a finite element simulation~\cite{Bassi-BR:97,cirak2000subdivision,sevilla2011nurbs,soghrati2016nurbs}, the most extended practice consists on maintaining the shape of the elements during the degree adaptive process~\cite{giorgiani2013hybridizable,karban2013advanced,giorgiani2014hybridizable}. In the majority of cases, a polynomial representation of the boundary is selected whereas the polynomial degree of the functional approximation changes at each iteration of the degree adaptive scheme. 

This work analyses and discusses three approaches to perform a degree adaptive process in domains with curved boundaries. The first one corresponds to the approach typically employed in practice, consisting of fixing the shape of the curved elements and changing the degree of the functional approximation as dictated by the degree adaptivity procedure. The second approach proposed in this work is to employ the so-called NURBS-enhanced finite element method (NEFEM) that enables to exactly represent the geometry of the computational domain, given by a CAD model, irrespectively of the degree of the polynomials used to approximate the solution. The third approach, despite not considered useful from a practical point of view, consists of changing the geometry representation of the computational domain to represent with the same degree of polynomials both the geometry and the solution at each iteration of the degree adaptive process. This approach is not considered of interest from a practical point of view because it requires communication with the CAD model at each iteration and re-generation of nodal distributions for curved elements. 

The second approach proposed here considers, for the first time, the combination of the so-called NURBS-enhanced finite element method (NEFEM) and the HDG rationale. The resulting method combines all the advantages of both methods, that is the efficiency of HDG and the ability of NEFEM to decouple the functional approximation from the geometric representation, usually tied in traditional isoparametric implementations. 

A number of numerical examples is considered in order to compare the different degree adaptivity approaches. Furthermore, this work presents a simple idea to verify computational methods that are able to use different degrees of approximation for the solution in different elements. The idea is based on an existing local a priori error estimator developed in~\cite{diez1999unified} for elliptic problems. 

The remainder of the paper is organised as follows. Sections~\ref{sc:HDGstokes} briefly presents the model problem considered (i.e. Stokes flows)  and the HDG formulation. The spatial discretisation of the HDG weak formulation is presented in Section~\ref{sc:spatialDiscretisation} for both isoparametric and NEFEM, with particular emphasis on the differences between both formulations. The details about the proposed error estimator and degree adaptivity process proposed are presented in Section~\ref{sc:adaptivity}, including a discussion of the three approaches considered to perform a degree adaptive process. In Section~\ref{sc:validation} a simple technique to verify the implementation of a solver with variable degree of approximation is presented and used to test the implementation of the HDG code for Stokes flows with isoparametric and NEFEM. Section~\ref{sc:comparison} presents a comparison of the different degree adaptive approaches and a number of numerical examples are used in Section to show the potential of the proposed approach. Finally, Section~\ref{sc:conclusions} summarises the main conclusions of the work that has been presented.

%==========================================================================
\section{Hybridisable discontinuous Galerkin for Stokes flow} \label{sc:HDGstokes}
%==========================================================================

%--------------------------------------------------------------------------
\subsection{Problem statement}  \label{sc:stokes}
%--------------------------------------------------------------------------

Let us consider an open bounded domain $\Omega\in\mathbb{R}^{\mathtt{n_{sd}}}$ with boundary $\partial\Omega$, where $\mathtt{n_{sd}}$ the number of spatial dimensions. The strong form of the stationary Stokes problem is obtained by neglecting the transient and convective effects in the full incompressible Navier-Stokes equations~\cite{Donea-Huerta}. The so-called velocity-pressure formulation is obtained by invoking the Stoke's law and results in 
\begin{equation}\label{eq:StokesStrong}
\left\{\!\begin{aligned}
-\Div (\nu\grad \bu - p\bI)      &= \bs    &&  \text{in $\Omega$},  \\
\grad\cdot \bu                   &= 0      &&  \text{in $\Omega$},  \\
\bu                              &= \bu_D  &&  \text{on $\Gamma_D$}, \\
\bn \cdot (\nu\grad \bu - p\bI)  &= \bt    &&  \text{on $\Gamma_N$}, \\
\end{aligned}\right.\\
\end{equation}
where $\bu$ is the velocity vector, $\nu$ is the kinematic viscosity, $p$ denotes the dynamic pressure, $\bs$ is a body force, $\bu_D$ is the imposed velocity on the Dirichlet boundary $\Gamma_D$, $\bn$ is the outward unit normal vector to $\partial\Omega$ and $\bt$ is the \textit{pseudo-traction} vector imposed on the Neumann boundary $\Gamma_N$.
The disjoint boundaries $\Gamma_D$ and $\Gamma_N$ satisfy $\partial\Omega = \overline{\Gamma}_D \cup \overline{\Gamma}_N$.

In what follow, $\scae{D}{\cdot}{\cdot}$ denotes the $\eltwo$ scalar product in a generic subdomain $D \subset \Omega$, that is
\begin{equation*}
\scae{D}{u}{v}     = \int_{D}uv d\Omega \qquad 
\scae{D}{\bu}{\bv} = \int_{D}\bu \cdot \bv d\Omega \quad \text{and} \quad
\scae{D}{\bU}{\bV} = \int_{D}\bU \bddot \bV d\Omega,
\end{equation*}
for scalars, vectors and second order tensors respectively. Analogously, $\scab{S}{\cdot}{\cdot}$ denotes the $\eltwo$ scalar product in any domain $S\subset\Gamma\cup\partial\Omega$.

The free divergence condition in Equation~\eqref{eq:StokesStrong} induces the compatibility condition
\begin{equation} \label{eq:compatibilityCondition}
\scab{\Gamma_D}{\bu_D \cdot \bn}{1} + \scab{\Gamma_N}{\bu \cdot \bn}{1} = 0.
\end{equation}

It is worth noting that, if only Dirichlet boundary conditions are considered (i.e. $\Gamma_N = \emptyset$), an additional constraint on the pressure must be imposed to avoid its indeterminacy. It is common~\cite{nguyen2010hybridizable,cockburn2014devising,cockburn2009derivation} to impose the mean pressure on the element boundary, namely
\begin{equation} \label{eq:constraintDirichlet}
\scab{\partial \Omega}{p}{1} = 0.
\end{equation}
%

%--------------------------------------------------------------------------
\subsection{HDG weak formulation}  \label{sc:HDGweak}
%--------------------------------------------------------------------------

The domain $\Omega$ is assumed partitioned in $\numel$ disjoint subdomains $\Omega_e$
%
%\begin{equation*}
%\overline{\Omega}=\bigcup_{i=1}^{\numel}\overline{\Omega}_e,\quad\Omega_e\cap\Omega_l=\emptyset\; \text{for}\;e\neq l,
%\end{equation*}
%
with boundaries $\partial\Omega_e$, which define an internal interface $\Gamma$
\begin{equation} \label{eq:meshSkeleton}
\Gamma : =\left[\bigcup_{e=1}^{\mathtt{n_{el}}}\partial\Omega_e\right]\setminus \partial\Omega.
\end{equation}

The corresponding strong form of the Stokes system given in Equation~\eqref{eq:StokesStrong} can be written in mixed form and in the broken computational domain as
\begin{align}
	\centering
	\begin{cases}
		\label{eq:StokesMixedBrokenDomain}
		\begin{aligned}
			\bL +\grad \bu               &= \bm{0} &&  \text{in $\Omega_e$}, \\
			\grad\cdot(\nu\bL + p\bm{I}) &= \bs    &&  \text{in $\Omega_e$}, \\			
			\grad\cdot \bu               &= 0      &&  \text{in $\Omega_e$}, \\
			\bu                          &= \bu_D  &&  \text{on $\partial\Omega_e\cap\Gamma_D$}, \\
			\bu                          &= \hu    &&  \text{on $\partial\Omega_e\setminus \Gamma_D$} , \\			
			\bn\cdot (\nu\bL+p\bm{I})    &= -\bt   &&  \text{on $\partial\Omega_e\cap\Gamma_N$}, \\
\jump{\bu\otimes\bn}         &= \bm{0} &&  \text{on $\Gamma$} , \\
			\jump{\bn\cdot(\nu\bL+p\bI)} &= \bm{0} &&  \text{on $\Gamma$}, 
		\end{aligned}
	\end{cases}
\end{align}
for $e=1,\dotsc,\numel$, where $\bI$ is the identity tensor of dimension $\nsd$, $\bL = -\grad\bu$ is a new variable (the second order velocity gradient tensor) which is introduced after splitting the second order momentum conservation equation in two first order equations and $\hu$ is an independent variable representing the trace of the solution in $\partial\Omega_e\setminus \Gamma_D$.

The free divergence condition in Equation~\eqref{eq:StokesMixedBrokenDomain} induces the compatibility condition
\begin{equation} \label{eq:compatibilityConditionElem}
\scab{\partial \Omega_e \cap \Gamma_D}{\bu_D \cdot \bn}{1} + \scab{\partial \Omega_e \setminus \Gamma_D}{\hu \cdot \bn}{1} = 0.
\end{equation}

The last two equations in \eqref{eq:StokesMixedBrokenDomain} impose the continuity of velocity and continuity of the normal component of the pseudo-stress across the interior faces respectively, where the \textit{jump} $\jump{\cdot}$ operator has been introduced following the definition in~\cite{AdM-MFH:08}, such that, along each portion of the interface $\Gamma$ it sums the values from the element on the left and right of say, $\Omega_e$ and $\Omega_l$ , namely
\begin{equation*}
\jump{\odot} = \odot_e + \odot_l.
\end{equation*}

The HDG method solves problem \eqref{eq:StokesMixedBrokenDomain} in two stages, see for instance~\cite{Jay-CGL:09,Cockburn-CDG:08,Nguyen-NPC:09,Nguyen-NPC:09b,Nguyen-NPC:11}. First, an element-by-element problem is defined with $(\bL,\bu,p)$ as unknowns. This is the so-called local problem and is given by
\begin{align}
\centering
\begin{cases}
\label{eq:stokesLocalProblem}
\begin{aligned}
\bL +\grad \bu               &= \bm{0} &&  \text{in $\Omega_e$}, \\
\grad\cdot(\nu\bL + p\bm{I}) &= \bs    &&  \text{in $\Omega_e$}, \\
\grad\cdot \bu               &= 0      &&  \text{in $\Omega_e$}, \\
\bu                          &= \bu_D  &&  \text{on $\partial\Omega_e\cap\Gamma_D$}, \\
\bu                          &= \hu    &&  \text{on $\partial\Omega_e\setminus \Gamma_D$} , \\
\scab{\partial \Omega_e}{p}{1} & = \rho_e,  &\quad &
\end{aligned}
\end{cases}
\end{align}
for $e=1,\dotsc,\numel$, where the last equation in~\eqref{eq:stokesLocalProblem} has been introduced to remove the indeterminacy of the pressure and $\rho_e$ denotes the mean pressure on the boundary of element $\Omega_e$. The local problem is used to obtain the solution in each element, $\bL$, $\bu$ and $p$, for $e=1,\dotsc,\numel$, in terms of $\hu$ and $\rho$ along the interface $\Gamma \cup \Gamma_N$. 

Second, a global problem is defined to determine the traces of the velocity and the mean pressure, denoted by $\hu$ and $\rho$, on the element boundaries. This is given by
\begin{align}
\centering
\begin{cases}
\label{eq:stokesGlobalProblemWithCont}
\begin{aligned}
\jump{\bu\otimes\bn}         &= \bm{0} &&  \text{on $\Gamma$}, \\
\jump{\bn\cdot(\nu\bL+p\bI)} &= \bm{0} &&  \text{on $\Gamma$}, \\
\bn\cdot (\nu\bL+p\bm{I})    &= -\bt   &&  \text{on $\Gamma_N$}
\end{aligned}
\end{cases}
\end{align}
where the first equation is automatically satisfied due to the unique definition of the hybrid variable $\hu$ on each face of the mesh skeleton and the condition $\bu=\hu$  on $\Gamma$, as imposed in the local problem. %Hence, the global problem reads
%
%\begin{align}
%\centering
%\begin{cases}
%\label{eq:stokesGlobalProblem}
%\begin{aligned}
%\jump{\bn\cdot(\nu\bL+p\bI)} &= \bm{0} &&  \text{on $\Gamma$}, \\
%\bn\cdot (\nu\bL+p\bm{I})    &= -\bt   &&  \text{on $\Gamma_N$}
%\end{aligned}
%\end{cases}
%\end{align}

The weak formulation for each element equivalent to~\eqref{eq:stokesLocalProblem} is as follows: for $e=1,\dotsc,\numel$, given $\bu_D$ on $\Gamma_D$ and $\hu$ on $\Gamma \cup \Gamma_N$, find $(\bL, \bu, p) \in \left[\sobo(\Omega_e)\right]^{\nsd\times\nsd} \times \left[\sobo(\Omega_e)\right]^{\nsd} \times \eltwo(\Omega_e)$ that satisfies
\begin{subequations}\label{eq:stokesWeakLocalWithNumFlux}
\begin{align} 
- \scae{\Omega_e}{\tL}{\bL} & + \scae{\Omega_e}{\Div \tL}{\bu}  = \scab{\partial \Omega_e \cap \Gamma_D}{\bn \cdot \tL}{\bu_D} + \scab{\partial \Omega_e \setminus \Gamma_D }{\bn \cdot \tL}{\hu} \label{eq:stokesWeakLocalWithNumFlux1}
\\
 -\scae{\Omega_e}{\grad \tu}{\nu \bL} & -\scae{\Omega_e}{\Div \tu}{p}  + \scab{\partial\Omega_e}{\tu}{\bn \cdot \widehat{(\nu\bL + p \bI)} } = \scae{\Omega_e}{\tu}{\bs} \label{eq:stokesWeakLocalWithNumFlux2}
\\
\scae{\Omega_e}{\grad \tdiv}{\bu} & = \scab{\partial \Omega_e \cap \Gamma_D}{\tdiv}{\bu_D \cdot \bn} + \scab{\partial \Omega_e \setminus \Gamma_D }{\tdiv}{\hu \cdot \bn} \label{eq:stokesWeakLocalWithNumFlux3}
\\
\scab{\partial \Omega_e}{p}{1} &  = \rho_e
\label{eq:stokesWeakLocalWithNumFlux4}
\end{align}
\end{subequations}
for all $(\tL, \tu, \tdiv) \in \left[\sobo(\Omega_e)\right]^{\nsd\times\nsd} \times \left[\sobo(\Omega_e)\right]^{\nsd} \times \eltwo(\Omega_e)$, where the numerical trace of the normal flux is defined as
\begin{equation} \label{eq:stokesNumericalFlux}
\bn \cdot \widehat{(\nu\bL + p \bI)} :=
\begin{cases}
\bn \cdot (\nu\bL + p \bI) +\tau (\bu - \bu_D) & \text{on $\partial\Omega_e\cap\Gamma_D$}\\
\bn \cdot (\nu\bL + p \bI) +\tau (\bu - \hu  ) & \text{elsewhere},
\end{cases}
\end{equation}
with $\tau$ being a stabilization parameter, whose selection has an important effect on the stability, accuracy and convergence properties of the resulting HDG method. The influence of the stabilization parameter has been studied extensively by Cockburn and co-workers, see for instance~\cite{Jay-CGL:09,Cockburn-CDG:08,Nguyen-NPC:09,Nguyen-NPC:09b,Nguyen-NPC:10,Nguyen-NPC:11}.

Introducing the definition of the numerical trace of the normal flux in Equation~\eqref{eq:stokesWeakLocalWithNumFlux} leads to the weak form of the local problem: for $e=1,\dotsc,\numel$, find $(\bL, \bu, p) \in \left[\sobo(\Omega_e)\right]^{\nsd\times\nsd} \times \left[\sobo(\Omega_e)\right]^{\nsd} \times \eltwo(\Omega_e)$ such that
\begin{subequations}\label{eq:stokesWeakLocal}
	\begin{align} 
	- \scae{\Omega_e}{\tL}{\bL}   & + \scae{\Omega_e}{\Div \tL}{\bu}  = \scab{\partial \Omega_e \cap \Gamma_D}{\bn \cdot \tL}{\bu_D} + \scab{\partial \Omega_e \setminus \Gamma_D }{\bn \cdot \tL}{\hu}  ,
	\label{eq:stokesWeakLocal1}
	\\
	\begin{split}
	-  \scae{\Omega_e}{\grad \tu}{\nu \bL}  &  + \scab{\partial\Omega_e}{\tu}{ \bn \cdot   \nu\bL}  + \scab{\partial\Omega_e}{\tu}{\tau \bu } -\scae{\Omega_e}{\Div \tu}{p}  + \scab{\partial\Omega_e}{\tu}{p \bn }  
	\\  
	&  = \scae{\Omega_e}{\tu}{\bs} + \scab{\partial\Omega_e \cap \Gamma_D}{\tu}{\tau \bu_D } + \scab{\partial\Omega_e \setminus \Gamma_D }{\tu}{\tau \hu }, 
	\label{eq:stokesWeakLocal2}
	\end{split}
	\\
	\scae{\Omega_e}{\grad \tdiv}{\bu}  &= \scab{\partial \Omega_e \cap \Gamma_D}{\tdiv}{\bu_D \cdot \bn} + \scab{\partial \Omega_e \setminus \Gamma_D }{\tdiv}{\hu \cdot \bn} 
	\label{eq:stokesWeakLocal3}
	\\
	\scab{\partial \Omega_e}{p}{1}  & = \rho
	\label{eq:stokesWeakLocal4}
	\end{align}
\end{subequations}
for all $(\tL, \tu, \tdiv) \in \left[\sobo(\Omega_e)\right]^{\nsd\times\nsd} \times \left[\sobo(\Omega_e)\right]^{\nsd} \times \eltwo(\Omega_e)$.

For the global problem, the weak formulation equivalent to~\eqref{eq:stokesGlobalProblemWithCont} is: find $\hu \in \left[\eltwo(\Gamma \cup \Gamma_N)\right]^{\nsd}$ and $\rho \in \mathbb{R}^{\numel}$ that satisfies 
\begin{subequations} \label{eq:stokesWeakGlobalWithNumFlux}
\begin{align}
%\begin{split}
\sum_{e=1}^{\numel} \Big\{ \scab{\partial\Omega_e\setminus\partial\Omega}{\that}{\bn\cdot \widehat{(\nu\bL+p\bI)}}  & + \scab{\partial \Omega_e \cap \Gamma_N}{\that}{\bn\cdot \widehat{(\nu\bL+p\bI)} + \bt}  \Big\}= 0
\label{eq:stokesWeakGlobalWithNumFlux1}
\\
\scab{\partial \Omega_e \setminus \Gamma_D }{\hu \cdot \bn}{1} & = - \scab{\partial \Omega_e \cap \Gamma_D}{\bu_D  \cdot \bn}{1}  \qquad \text{for $e=1,\dotsc,\numel$}, \label{eq:stokesWeakGlobalWithNumFlux2}
\end{align}
\end{subequations}
for all $\that \in \left[\eltwo(\Gamma \cup \Gamma_N)\right]^{\nsd}$.

Introducing the definition of the numerical trace of the normal flux in Equation~\eqref{eq:stokesWeakGlobalWithNumFlux} leads to the weak form of the global problem: find $\hu \in \left[\eltwo(\Gamma \cup \Gamma_N)\right]^{\nsd}$ and $\rho \in \mathbb{R}^{\numel}$ such that, for all $\that \in \left[\eltwo(\Gamma \cup \Gamma_N)\right]^{\nsd}$, 
\begin{subequations}\label{eq:stokesWeakGlobal}
	\begin{align} 
	\begin{split}
	\sum_{e=1}^{\numel} \Big\{  \Big.
	 \scab{\partial\Omega_e\setminus \Gamma_D }{\that}{\bn\cdot \nu\bL}   & +
	\scab{\partial\Omega_e\setminus (\Gamma_D \cup \Gamma_S)}{\that}{ \tau \bu } +  \scab{\partial\Omega_e\setminus \Gamma_D ) }{\that}{p \bn} 
	\\
	&   
	 \Big.  - \scab{\partial\Omega_e\setminus \Gamma_D }{\that}{ \tau \hu}   \Big\}  
	 = 	-\sum_{e=1}^{\numel} \left\{ \scab{\partial\Omega_e\cap\Gamma_N}{\that}{\bt}  \right\},
	\label{eq:stokesWeakGlobal1}
	\end{split}
	\\
	 & \scab{\partial \Omega_e \setminus \Gamma_D}{\hu \cdot \bn}{1}   = - \scab{\partial \Omega_e \cap \Gamma_D}{\bu_D  \cdot \bn}{1}. 
	\label{eq:stokesWeakGlobal2}
	\end{align}
\end{subequations}

%==========================================================================
\section{Spatial discretisation} \label{sc:spatialDiscretisation}
%==========================================================================

This section presents the discretisation of the HDG weak forms derived in the previous section. Both the standard isoparametric and the so-called NEFEM formulations are presented. Special attention is paid to the differences between both formulations as this represents the first time NEFEM is considered in an HDG framework. 

%--------------------------------------------------------------------------
\subsection{Isoparametric elements} \label{sc:isoparametric}
%--------------------------------------------------------------------------

Standard isoparametric formulations map each element $\Omega_e$ and face $\Gamma_e$ in the physical domain into a reference element, $\widetilde{\Omega}$, and a reference face, $\widetilde{\Gamma}$, where polynomial functional approximations characterize the discrete finite dimensional spaces. Namely, $\spP^k (\widetilde{\Omega})$ and $\spP^{\hat{k}} (\widetilde{\Gamma})$ are the spaces of polynomial functions of degree at most $k \geq 1$ and $\hat{k} \geq 1$ in the reference element and the reference face respectively. Finally, the approximations for each variable are defined as 
\begin{align} 
\label{eq:veloApproxFEM}
\bu(\bxi) \simeq \bu^h(\bxi) &= \sum_{j=1}^{\nen} \nodalu_j N_j(\bxi)
                                           &&\in\bigl[\{ v\in\eltwo(\Omega); \, v\vert_{\Omega_e} \in \spP^k (\widetilde{\Omega}) \}\bigr]^{\nsd} ,\\
\label{eq:pressApproxFEM}
p(\bxi) \simeq p^h(\bxi)       &= \sum_{j=1}^{\nen} \nodalp_j N_j(\bxi)
                                           &&\in\{ q\in\eltwo(\Omega); \, q\vert_{\Omega_e} \in \spP^k (\widetilde{\Omega}) \} ,\\
\label{eq:LApproxFEM}
\bL(\bxi) \simeq \bL^h(\bxi) &= \sum_{j=1}^{\nen} \nodalL_j N_j(\bxi)
                                           &&\in\bigl[\{ v\in\eltwo(\Omega); \, v\vert_{\Omega_e} \in \spP^k (\widetilde{\Omega}) \}\bigr]^{\nsd\times\nsd} ,\\
\label{eq:traceApproxFEM}
\hu(\bet) \simeq \hu^h(\bet) &= \sum_{j=1}^{\nfn} \nodaluh_j \hat{N}_j(\bet)
                                           &&\in\bigl[\{ v\in\eltwo(\Gamma); \, v\vert_{\Gamma_e} \in\spP^{\hat{k}} (\widetilde{\Gamma}) \}\bigr]^{\nsd} ,
\end{align} 
where $\nodalu_j$ , $\nodalp_j$, $\nodalL_j$ and $\nodaluh_j$ are nodal values, $N_j$ are polynomial shape functions of order $k$ in the reference element, $\nen$ is the number of nodes per element, $\hat{N}_j$ are polynomial shape functions of order $\hat{k}$ in the reference face and $\nfn$ is the corresponding number of nodes per face. Note that equal interpolation is used for all element variables (i.e. velocity, pressure and gradient of velocity). Recall that HDG allows for equal interpolation because of the numerical fluxes and the stabilization parameter $\tau$. They ensure solvability and stability, see \cite{Cockburn-CDG:08}, without the need of an enriched space for the gradient variable, or a reduced space for the trace variable.

An isoparametric mapping is used to link the reference element $\hat{\Omega}$ and the computational element $\Omega_e^h$
\begin{equation}\label{eq:isoparametricMapping}
\begin{aligned}
\bm{\varphi}\; :
\widetilde{\Omega} \subset \mathbb{R}^{\nsd} \ &\longrightarrow \Omega_e^h \subset \mathbb{R}^{\nsd} \\
\bxi &\longmapsto
\bm{\varphi}(\bxi):= \sum_{j=1}^{\nen} \bx_j N_j(\bxi),
\end{aligned}
\end{equation}
where $\bx_j$ are the nodal coordinates of the computational element $\Omega_e^h$. 

It is worth noting that in general, when the physical element $\Omega_e$ is curved, the isoparametric mapping is non-linear and the approximation defined in the reference element do not induce a polynomial interpolation in the physical space. In addition, the computational element $\Omega_e^h$ is just an approximation of $\Omega_e$, see~\cite{sevilla2011comparison} for a detailed discussion.

Similarly, an isoparametric mapping is used to link the reference face $\hat{\Gamma}$ and the computational face $\Gamma_e^h$
\begin{equation}\label{eq:isoparametricMappingFace}
\begin{aligned}
\bm{\psi}\; :
\widetilde{\Gamma} \subset \mathbb{R}^{\nsd-1} &\longrightarrow \Gamma_e^h \subset \mathbb{R}^{\nsd} \\
\bet &\longmapsto
\bm{\psi}(\bet):= \sum_{j=1}^{\nfn} \bx_j \hat{N}_j(\bet),
\end{aligned}
\end{equation}
where $\bx_j$ denote the face nodal coordinates.

Using the mappings in Equations~\eqref{eq:isoparametricMapping} and \eqref{eq:isoparametricMappingFace}, the integrals appearing in the weak form of the local problems are transformed to the reference element and reference face/edge respectively. Then, the nodal interpolations given by Equations~\eqref{eq:veloApproxFEM}  to \eqref{eq:traceApproxFEM} are introduced, leading to a system of equations for each element with the following structure
\begin{equation} \label{eq:localProblemFEM}
\left[\begin{array}{cccc}
\mat{A}_{LL} & \mat{A}_{Lu} & \bm{0} & \bm{0}\\
\mat{A}_{uL} & \mat{A}_{uu} & \mat{A}_{up} & \bm{0}\\
\bm{0} & \mat{A}_{pu} & \bm{0} &\bm{a}_{\rho p}^T\\
\bm{0} & \bm{0} &\bm{a}_{\rho p} & 0
\end{array}\right]
\left\{\begin{array}{c}
\nodalLV \\
\nodaluV \\
\nodalpV \\
\zeta
\end{array}\right\} = 
\left\{\begin{array}{c}
\vect{f}_L\\
\vect{f}_u\\
\vect{f}_p\\
0
\end{array}\right\} +
\left[\begin{array}{c}
\mat{A}_{L\hat{u}}\\
\mat{A}_{u\hat{u}}\\
\mat{A}_{p\hat{u}}\\
\bm{0}
\end{array}\right]
\nodaluhV + 
\left\{\begin{array}{c}
\bm{0}\\
\bm{0}\\
\bm{0}\\
1
\end{array}\right\}\rho,
\end{equation}
where $\zeta$ is the Lagrange multiplier corresponding to the constraint of Equation~\eqref{eq:stokesWeakLocal4}.
% and 
%%
%\begin{equation}
%\mat{A} =
%\left[\begin{array}{cccc}
%\mat{A}_{LL} & \mat{A}_{Lu} & \bm{0} & \bm{0}\\
%\mat{A}_{uL} & \mat{A}_{uu} & \mat{A}_{up} & \bm{0}\\
%\bm{0} & \mat{A}_{pu} & \bm{0} &\bm{a}_{\rho p}^T\\
%\bm{0} & \bm{0} &\bm{a}_{\rho p} & 0
%\end{array}\right].
%\end{equation}

Analogously, using the isoparametric mappings, the nodal interpolations of the corresponding variables and introducing the expression of $\nodalLV$, $\nodaluV$ and $\nodalpV$ from Equation~\eqref{eq:localProblemFEM} in the global problem of Equation~\eqref{eq:stokesWeakGlobal}, a global system of equations is obtained
\begin{equation} \label{eq:GlobalSystemStokes}
\widehat{\mat{K}} \hat{\vect{U}} = \hat{\vect{f}},
\end{equation}
where the vector of unknowns $\hat{\vect{U}}$ contains the nodal values of the trace of the velocity on the elementa faces and the mean pressure within each element.

%The global matrix and right hand side are obtained by assembling the element contributions
%%
%\begin{equation}  \label{eq:globalProblemFEM}
%\widehat{\mat{K}}_e = 
%\left[\begin{array}{cc}
%\mat{B} \mat{Z}^f + \mat{A}_{\hat{u}\hat{u}} & \mat{B} \vect{z}^{\rho}\\
%\bm{\mathrm{a}}_{\rho\hat{u}} & 0
%\end{array}\right], 
%\qquad 
%\hat{\vect{f}}_e = 
%\left[\begin{array}{c}
%\vect{f}_{\hat{u}} - \mat{B} \mat{z}^f \\
%\vect{f}_{\rho}
%\end{array}\right],
%\end{equation}
%%
%where
%%
%\begin{equation} \label{eq:LocalProblemMatrixStokesAux}
%\mat{z}^f {=}
%\mat{A}^{-1}
%\left\{\begin{array}{c}
%\vect{f}_L\\
%\vect{f}_u\\
%\vect{f}_p\\
%0
%\end{array}\right\}, 
%%
%\,
%%
%\mat{Z}^f{=}
%\mat{A}^{-1}
%\left[\begin{array}{c}
%\mat{A}_{L\hat{u}}\\
%\mat{A}_{u\hat{u}}\\
%\mat{A}_{p\hat{u}}\\
%\bm{0}
%\end{array}\right],
%%
%\,
%%
%\vect{z}^{\rho} {=}
%\mat{A}^{-1}
%\left\{\begin{array}{c}
%\bm{0}\\
%\bm{0}\\
%\bm{0}\\
%1
%\end{array}\right\}, 
%%
%\,
%%
%\mat{B} {=}
%\left[\mat{A}_{\hat{u}L} \;\; \mat{A}_{\hat{u}u} \;\; \mat{A}_{\hat{u}p}
%\right].
%\end{equation}

%--------------------------------------------------------------------------
\subsection{NEFEM elements}  \label{sc:nefem}
%--------------------------------------------------------------------------

In NEFEM, the boundary of the computational domain $\partial \Omega$ is exactly represented by NURBS. In what follows, in order to simplify the presentation and without loss of generality, the NURBS are restricted to two dimensional problems, see~\cite{sevilla20113d} for a detailed description of the three dimensional case. An edge is given by  $\Gamma_e := \bm{C}([\lambda^e_a,\lambda^e_b])$, where $\bm{C}$ is the NURBS boundary parametrisation and $\lambda_a$ and $\lambda_b$ are the parametric coordinates (in the parametric space of the NURBS) of the end points of $\Gamma_e$.

The discrete approximations are defined now as:
\begin{align} 
\label{eq:veloApproxNEFEM}
\bu(\bx) \simeq \bu^h(\bx) &= \sum_{j=1}^{\nen} \nodalu_j N_j(\bx)
                                          &&\in\bigl[\{ v\in\eltwo(\Omega); \, v\vert_{\Omega_e} \in \spP^k (\Omega_e) \}\bigr]^{\nsd} ,\\
\label{eq:pressApproxNEFEM}
p(\bx) \simeq p^h(\bx)       &= \sum_{j=1}^{\nen} \nodalp_j N_j(\bx)
                                          &&\in\{ q\in\eltwo(\Omega); \, q\vert_{\Omega_e} \in \spP^k (\Omega_e) \} ,\\
\label{eq:LApproxNEFEM}
\bL(\bx) \simeq \bL^h(\bx) &= \sum_{j=1}^{\nen} \nodalL_j N_j(\bx)
                                          &&\in\bigl[\{ v\in\eltwo(\Omega); \, v\vert_{\Omega_e} \in \spP^k (\Omega_e) \}\bigr]^{\nsd\times\nsd} ,\\
\label{eq:traceApproxNEFEM}
\hu(\blam) \simeq \hu^h(\blam) &= \sum_{j=1}^{\nfn} \nodaluh_j \hat{N}_j(\blam)
                                           &&\in\bigl[\{ v\in\eltwo(\Gamma); \, v\vert_{[\lambda^e_a,\lambda^e_b]} \in\spP^{\hat{k}} ([\lambda^e_a,\lambda^e_b]) \}\bigr]^{\nsd} ,
\end{align}
where $\nodalu_j$, $\nodalp_j$, $\nodalL_j$ and $\nodaluh_j$ are nodal values, $N_j$ are polynomial shape functions of order $k$ in the physical element, $\nen$ is the number of nodes per element, $\hat{N}_j$ are polynomial shape functions of order $\hat{k}$ in $[\lambda^e_a,\lambda^e_b]$ and $\nfn$ is the corresponding number of nodes per face.

The main differences of NEFEM with respect to the isoparametric formulation are:
\begin{itemize}
	\item The exact description of the computational domain is considered by means of its NURBS boundary representation.
	\item The approximation of the elemental variables directly in the physical space, with Cartesian coordinates.	
	\item The approximation of the trace of the velocity is defined in the parametric space of the NURBS. It is worth noting that other options could be considered such as defining the approximation directly in the physical space. The main advantage of defining the approximation in the parametric space of the NURBS is that the number of unknowns remains the same as in the isoparametric formulation. In contrast, if the approximation of this variable is selected in the physical space it would require further degrees of freedom~\cite{sevilla2011comparison,sevilla2011nurbs}.
\end{itemize}

In addition, from the computational point of view, NEFEM uses specifically designed numerical quadratures that provide a more efficient alternative to standard quadratures defined in a reference triangle~\cite{sevilla2011numerical,sevilla20113d}. For instance, in two dimensions, the following mapping is introduced between a reference rectangle and the physical element
\begin{equation} \label{eq:mappingRect2Tri}
\begin{aligned}
\bm{\psi}\; :
R \subset \mathbb{R}^{\nsd} &\longrightarrow \Omega_e \subset \mathbb{R}^{\nsd} \\
\blam &\longmapsto \bm{\psi} (\blam):=
(1-\vartheta) \bm{C}(\lambda_1) + \lambda_2 \bx_I,
\end{aligned}
\end{equation}
where $R = [\lambda^e_a,\lambda^e_b] \times [0,1]$ and $\bx_I$ is the internal vertex of $\Omega_e$.

Using the mapping in Equation~\eqref{eq:mappingRect2Tri} and the NURBS boundary representation given by $\bm{C}$, the integrals appearing in the weak form of the local problems are transformed to the reference rectangle and the parametric space of the NURBS respectively. Then, the nodal interpolations given by Equations~\eqref{eq:veloApproxNEFEM}, \eqref{eq:pressApproxNEFEM}, \eqref{eq:LApproxNEFEM} and \eqref{eq:traceApproxNEFEM} are introduced, leading to a system of equations similar to Equation~\eqref{eq:localProblemFEM}. Analogously, the global problem with NEFEM leads to a global system of equations similar to Equation~\eqref{eq:GlobalSystemStokes}.

%==========================================================================
\section{Error estimation and adaptivity} \label{sc:adaptivity}
%==========================================================================

In HDG, the possibility to obtain a postprocessed solution~\cite{Jay-CGL:09} that converges at a higher rate (i.e. $k+2$) than the HDG solution, not only provides a higher accurate solution to the problem at hand but it can also be used to build an inexpensive, reliable and computable error estimator~\cite{giorgiani2014hybridizable,giorgiani2013hybridizable}. In this section, particular attention is paid to the fact that, when the degree of approximation is changed in a curved element, a choice must be made regarding the geometric definition of the element. 

An element by element measure of the error is defined by employing the HDG solution and the postprocessed solution as proposed in~\cite{giorgiani2014hybridizable}
\begin{equation} \label{eq:errorMeasure}
E_e = \left[ \frac{1}{|\Omega_e|} \int_{\Omega_e} \left( \bu^\star - \bu \right) \cdot \left( \bu^\star - \bu \right) d\Omega \right]^{1/2},
\end{equation}
where the normalisation becomes critical when meshes with different element sizes are considered~\cite{diez1999unified}.

For elliptic problems, and by using that the influence of pollution errors becomes negligible if the mesh is sufficiently refined in the area where the pollution error is generated~\cite{huerta2000error}, the following a priori error estimate was derived in~\cite{diez1999unified}
\begin{equation} \label{eq:aPrioriError}
\varepsilon_e = \| \bu - \bu_h \|_{\Omega_e} \leq C h_e^{k_e + 1 + \nsd/2},
\end{equation}
for $e=1,\ldots,\numel$, where $C$ is a constant, $h_e$ the characteristic element size of $\Omega_e$ and $k_e$ the degree of approximation used in $\Omega_e$. \RS{It is worth noting that the error estimate of Equations~\eqref{eq:aPrioriError} was initially derived for the standard finite element method but its extension to the HDG method is straightforward.}

By applying a standard Richardson extrapolation, it is possible to predict the required change in the degree of approximation in order to ensure that the error in each element is lower than a desired accuracy $\epsilon$, namely
\begin{equation} \label{eq:newP}
\Delta k_e = \left\lceil \frac{\log(\varepsilon/E_e)}{\log(h_e)} \right\rceil
\end{equation}
for $e=1,\ldots,\numel$, where $\lceil \cdot \rceil$ denotes the ceiling function.

The adaptive procedure consist on solving the Stokes problem using the HDG formulation as described in Section~\ref{sc:spatialDiscretisation} and estimating the required degree of approximation in each element according to Equation~\eqref{eq:newP}. The process is repeated until convergence is achieved, meaning that the error in each element $\varepsilon_e$ is lower than the desired error $\epsilon$.

%--------------------------------------------------------------------------
\subsection{Geometry update}  \label{sc:geometricUpdate}
%--------------------------------------------------------------------------

The technique described to drive a degree adaptive process only focuses on the degree of approximation used for the functional approximation, but in the presence of curved boundaries it is known that high-order approximations of both the solution and the geometry are required to exploit the full potential of a high-order method~\cite{Bassi-BR:97,Dey-Shephard:97,Krivodonova-KB:06,HO-Fluids}. This aspect is usually ignored as degree adaptive procedures are applied to problems involving polygonal boundaries, see for instance~\cite{demkowicz2005fully,pardo2010multigoal,vsolin2008arbitrary,dumbser2007arbitrary}. Here three options are discussed and assessed and compared later using numerical examples.

\begin{figure}[!tb]
	\centering
	\includegraphics[width=0.9\textwidth]{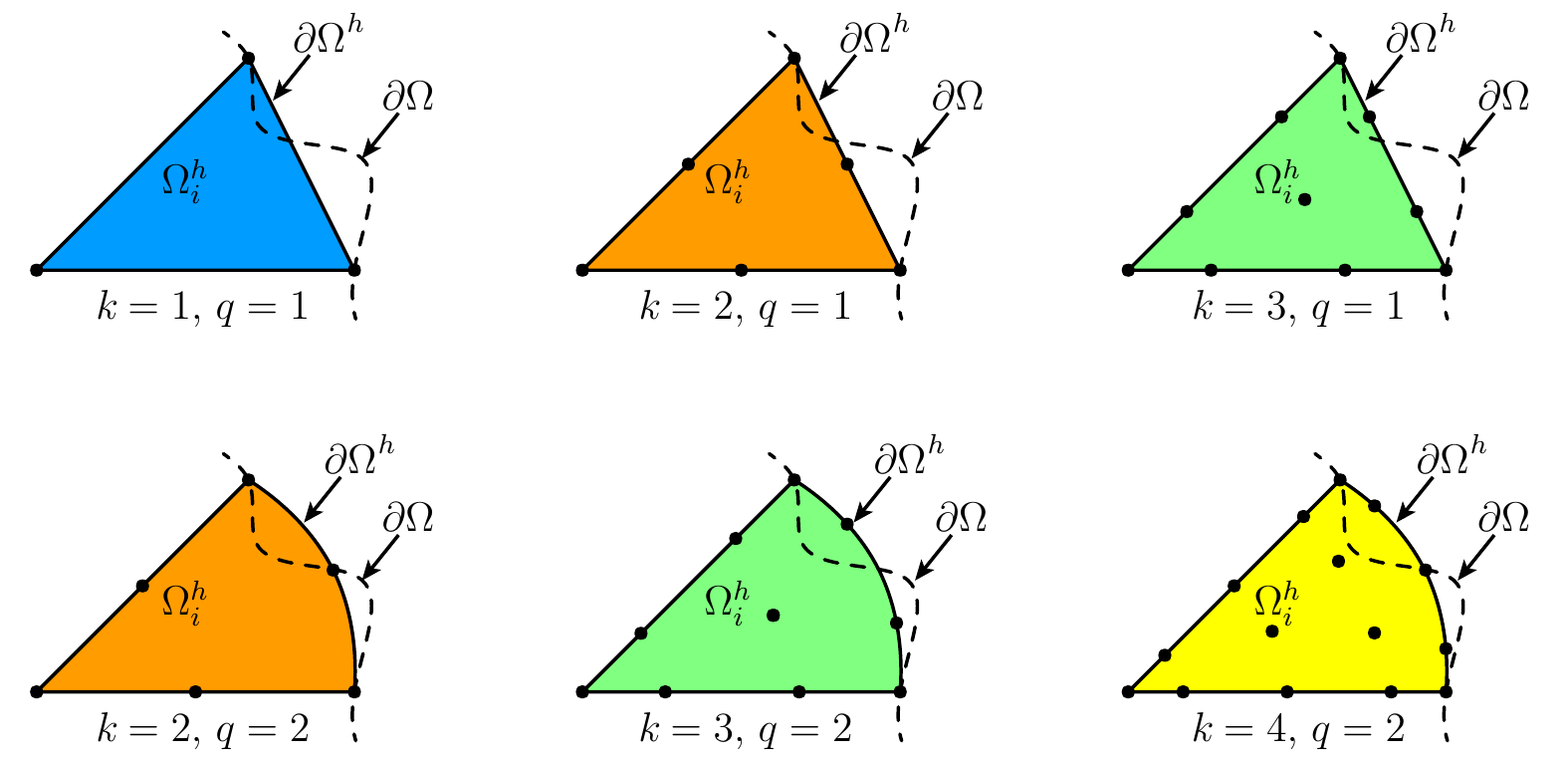}
	\caption{Illustration of a degree adaptation in an element with linear (top) and quadratic (bottom) approximation of the geometry.}
	\label{fig:pAdaptivityFEM0}
\end{figure}

The first, and the one typically considered in a degree adaptive process, technique consists of defining a polynomial representation of the curved boundaries that is maintained during the adaptive process, irrespectively of the degree of approximation used for the solution~\cite{giorgiani2013hybridizable,li2013numerical,giorgiani2014hybridizable}. This option is attractive because when the degree of approximation is changed in an element, there is no need to communicate with a CAD library to re-generate the nodal distributions in curved elements at each iteration of the degree adaptive process. The strategy is illustrated in Figure~\ref{fig:pAdaptivityFEM0}. The first row of plots show triangular elements where the geometric approximation is linear ($q=1$) and the polynomial degree of approximation of the solution increases from $k=1$ to $k=3$. The second row shows a similar situation where the boundary of the computational domain is described using quadratic polynomials ($q=2$) and the degree of approximation of the solution is increased from $k=2$ to $k=4$. The boundary of the computational domain is denoted by $\partial \Omega^h$ whereas the exact boundary is denoted by $\partial \Omega$.

\begin{figure}[!tb]
	\centering
	\includegraphics[width=0.9\textwidth]{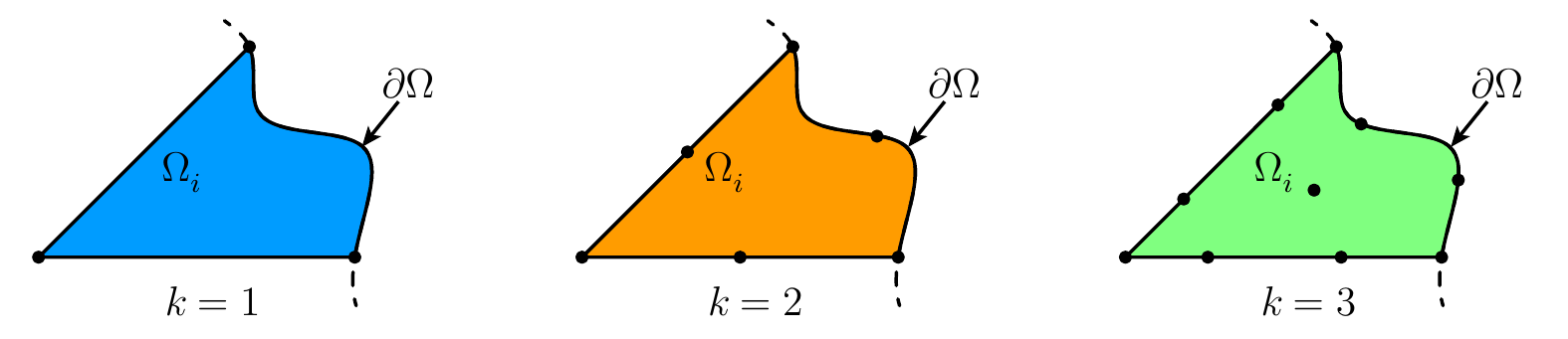}
	\caption{Illustration of a degree adaptation in a NEFEM element.}
	\label{fig:pAdaptivityNEFEM}
\end{figure}

The second alternative, proposed in this work, consists of using NEFEM, where the exact boundary representation of the computational domain is considered irrespectively of the degree of approximation considered for the solution. As NEFEM encapsulates the necessary information to define the approximation and perform the numerical integration in curved elements in contact with a NURBS boundary, communication with a CAD library is avoided. The strategy is illustrated in Figure~\ref{fig:pAdaptivityNEFEM}, showing a NEFEM element where the exact boundary representation is considered and a degree of approximation for the solution being increased from $k=1$ to $k=3$. 

\begin{figure}[!tb]
	\centering
	\includegraphics[width=0.9\textwidth]{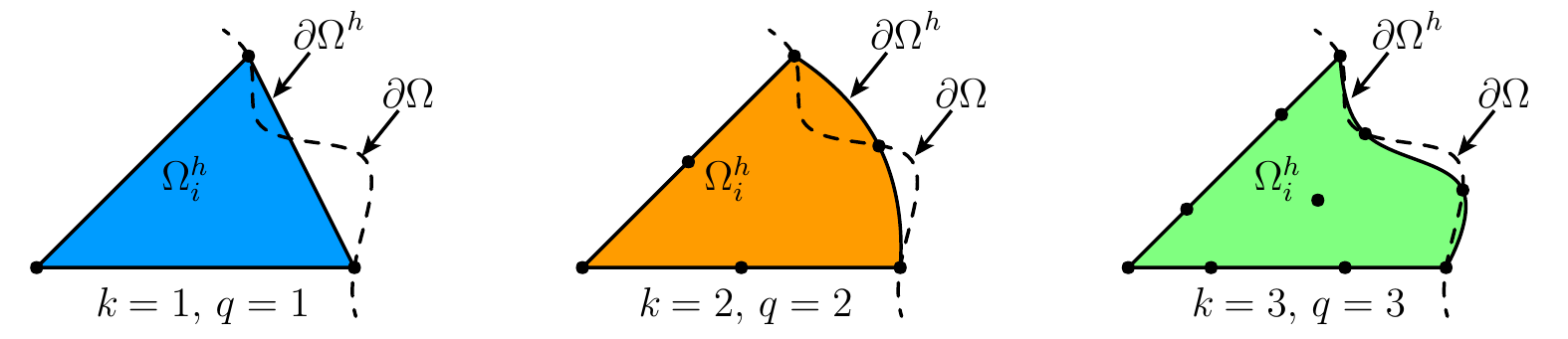}
	\caption{Illustration of a degree adaptation in an element where the same degree of approximation is used for both the solution and the geometry.}
	\label{fig:pAdaptivityFEM1eMeshes}
\end{figure}

A third alternative, not considered in practice, consists of communicating with the CAD model after each iteration of the degree adaptive process in order to re-generate the nodal distribution of curved elements by placing the nodes over the true boundary. The strategy is illustrated in Figure~\ref{fig:pAdaptivityFEM1eMeshes}, showing a triangular element where both the degree of the functional approximation and the degree of the polynomials used to approximate the solution are updated at each iteration.
This strategy has not been considered in practical applications due to the cost associated to communicating with the CAD model at each iteration.

\begin{remark}
It is important to note that the first strategy, where the geometry remains unchanged, does not guarantee the convergence of the numerical solution to the physical solution in domains with curved boundaries because the distance between the computational domain and the physical domain does not converge to zero with as the degree of approximation is increased, see~\cite{raviart1998introduction,ciarlet2002finite} for more details. For the second strategy, proposed here, convergence to the physical solution is guaranteed because no geometrical error is introduced~\cite{sevilla2011nurbs}. Finally, for the third approach, convergence is also guaranteed if the distance between the computational boundary and the physical boundary tends to zero as the order of the approximation is increased and the derivatives of the isoparametric mapping up to order $k+1$ are bounded by $h^s$, for $s=2,\ldots,k+1$~\cite{raviart1998introduction,ciarlet2002finite}, where $h$ denotes the characteristic element size. It is worth noting that a specifically designed nodal distribution for curved elements is required in the third approach to guarantee that the second hypothesis is fulfilled~\cite{ciarlet2002finite}.
\end{remark}

%--------------------------------------------------------------------------
\section{Validation of the HDG formulation with variable degree of approximation}  \label{sc:validation}
%--------------------------------------------------------------------------

The first example provides a novel and simple technique to fully validate a solver that employs variable different degree of approximation in different elements for the solution of elliptic problems. The idea consists of utilising the local a priori error estimate of Equation~\eqref{eq:aPrioriError} that states how the error, measured in an element, decreases when the mesh is refined. 

To illustrate the proposed technique and validate the HDG isoparametric and NEFEM implementations with variable degree of approximation, the Stokes equations are solved in a circle of radius 0.5 centred at $(0.5,0.5)$ with Dirichlet boundary conditions. The viscosity is considered as $\nu=1$ and source and boundary conditions are taken such that the analytical solution is given by
\begin{equation} \label{eq:analyticalStokes}
\bu = \begin{pmatrix} x^2(1-x)^2(2y-6y^2+4y^3)\\ -y^2(1-y)^2(2x-6x^2+4x^3) \end{pmatrix} 
\qquad 
p = x(1-x).
\end{equation}

\begin{figure}[!tb]
	\centering
	\subfigure[Mesh 1]{\includegraphics[width=0.32\textwidth]{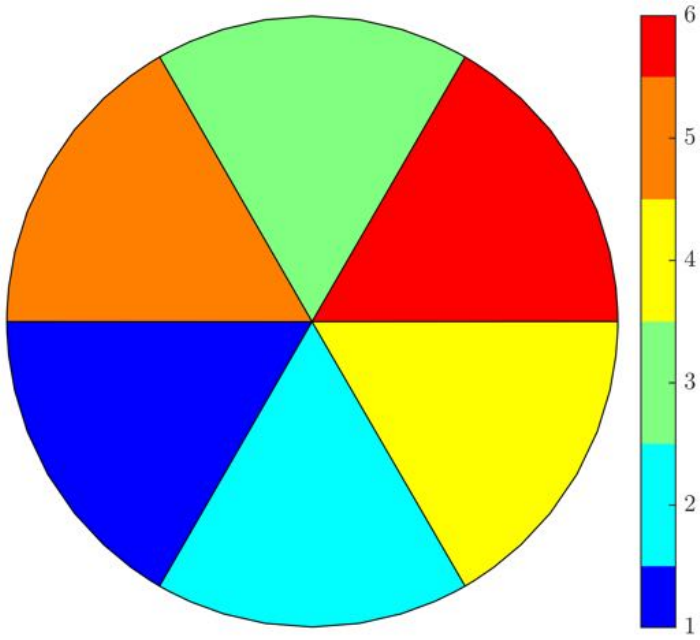}}
	\subfigure[Mesh 2]{\includegraphics[width=0.32\textwidth]{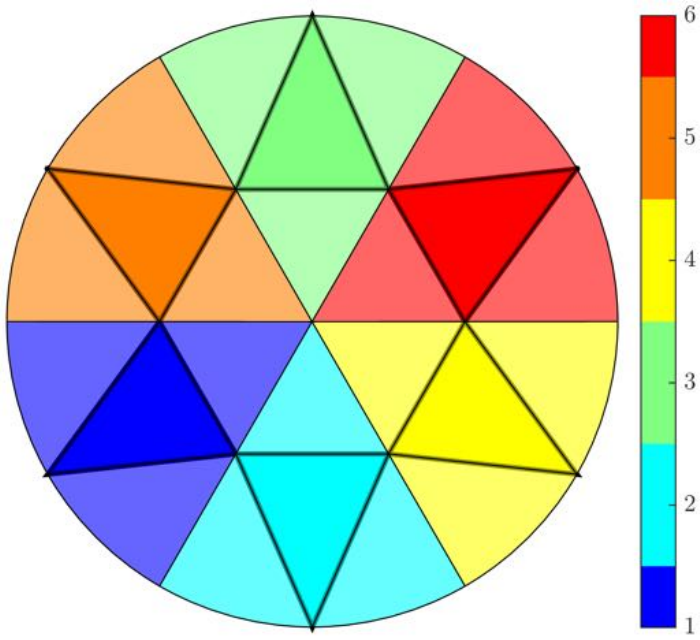}}
	\subfigure[Mesh 3]{\includegraphics[width=0.32\textwidth]{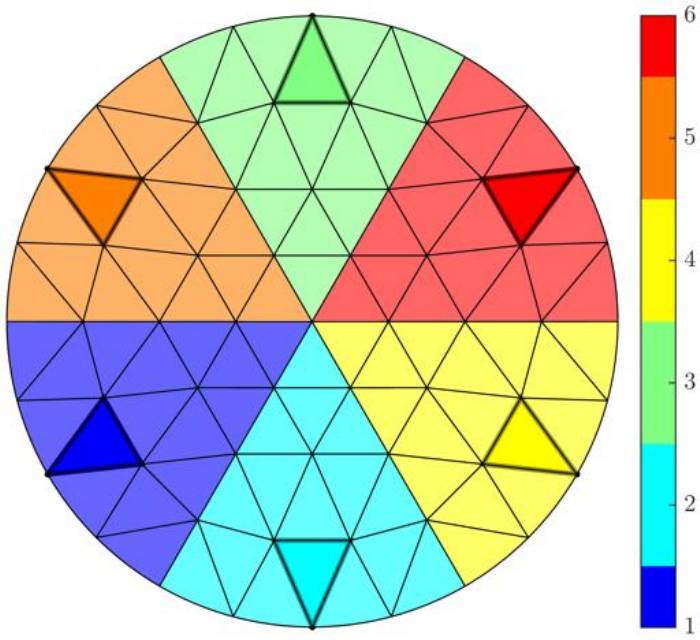}}
	\caption{First three NEFEM meshes where the colour indicates the degree of approximation used in each element and the highlighted element represents the region where the error is measured for each degree of approximation.}
	\label{fig:circleMeshesVarP}
\end{figure}

Six triangular meshes of the domain are generated using nested refinement. The first three meshes are shown in Figure~\ref{fig:circleMeshesVarP}, where the colour of each element represents the degree of approximation used, ranging from $k=1$ to $k=6$. 
%
%The initial mesh contains only six elements and a different degree of approximation is used in each element. In each mesh refinement step, the four elements that are created to split a triangle of size $h$ into four triangles of size $h/2$ inherit the degree of approximation. 
In each mesh, there is one element per degree of approximation highlighted with a thicker line and darker colour, representing the region where the error is measured to test the local a priori error estimate.
It is worth noting that the meshes shown in Figure~\ref{fig:circleMeshesVarP} are NEFEM meshes, as the exact boundary representation is always employed, even for $k=1$, whereas for the computations both NEFEM and isoparametric meshes are employed.

\begin{figure}[!tb]
	\centering
	\subfigure[HDG Isoparametric]{\includegraphics[width=0.49\textwidth]{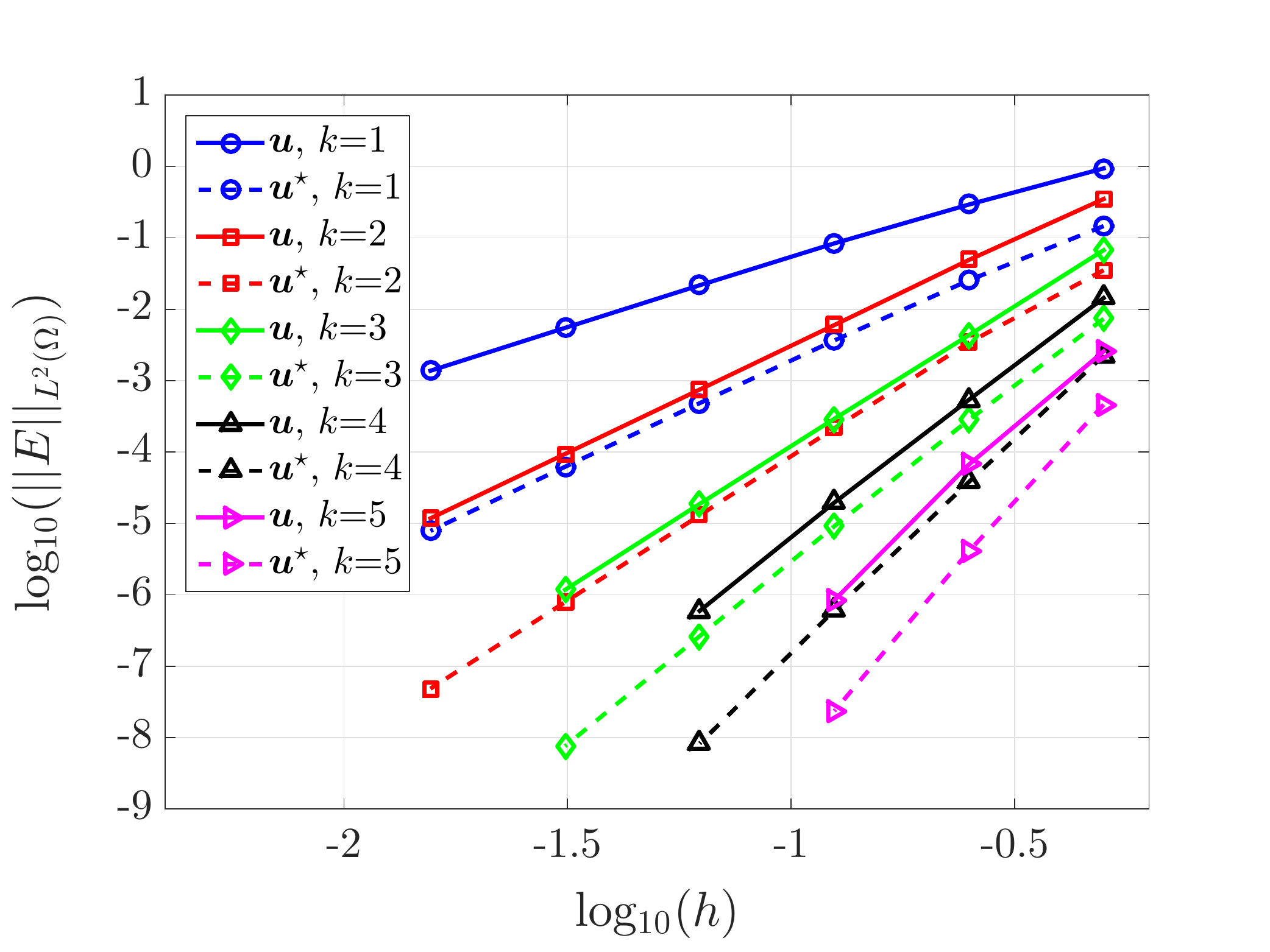}}
	\subfigure[HDG NEFEM]{\includegraphics[width=0.49\textwidth]{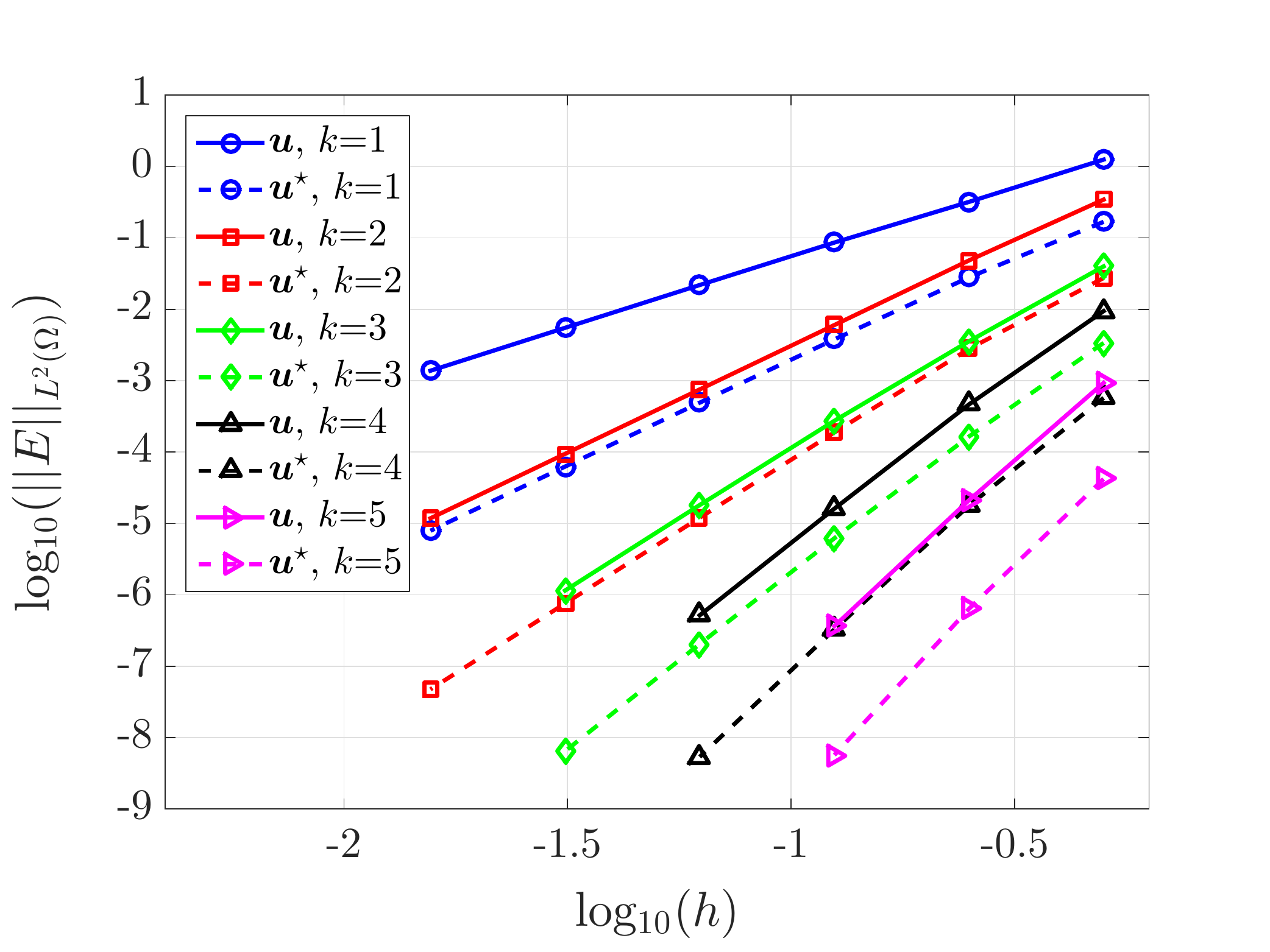}}
	\caption{Error of the solution $\bu$ and the postprocessed solution $\bu^\star$ in the $\eltwo(\Omega_e)$ norm for different degrees of approximation in each element.}
	\label{fig:circleHConvU}
\end{figure}

\begin{figure}[!tb]
	\centering
	\subfigure[HDG Isoparametric]{\includegraphics[width=0.49\textwidth]{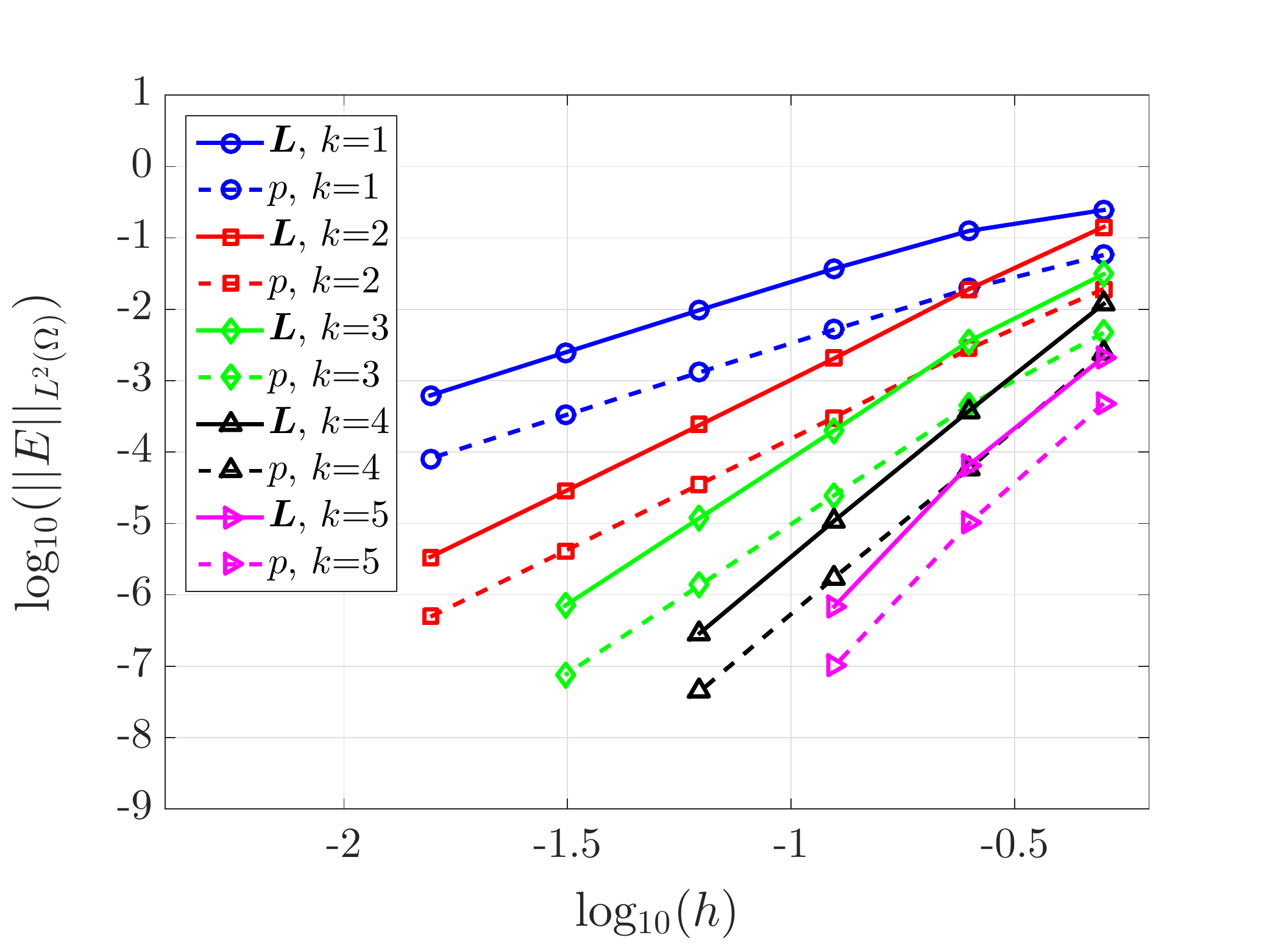}}
	\subfigure[HDG NEFEM]{\includegraphics[width=0.49\textwidth]{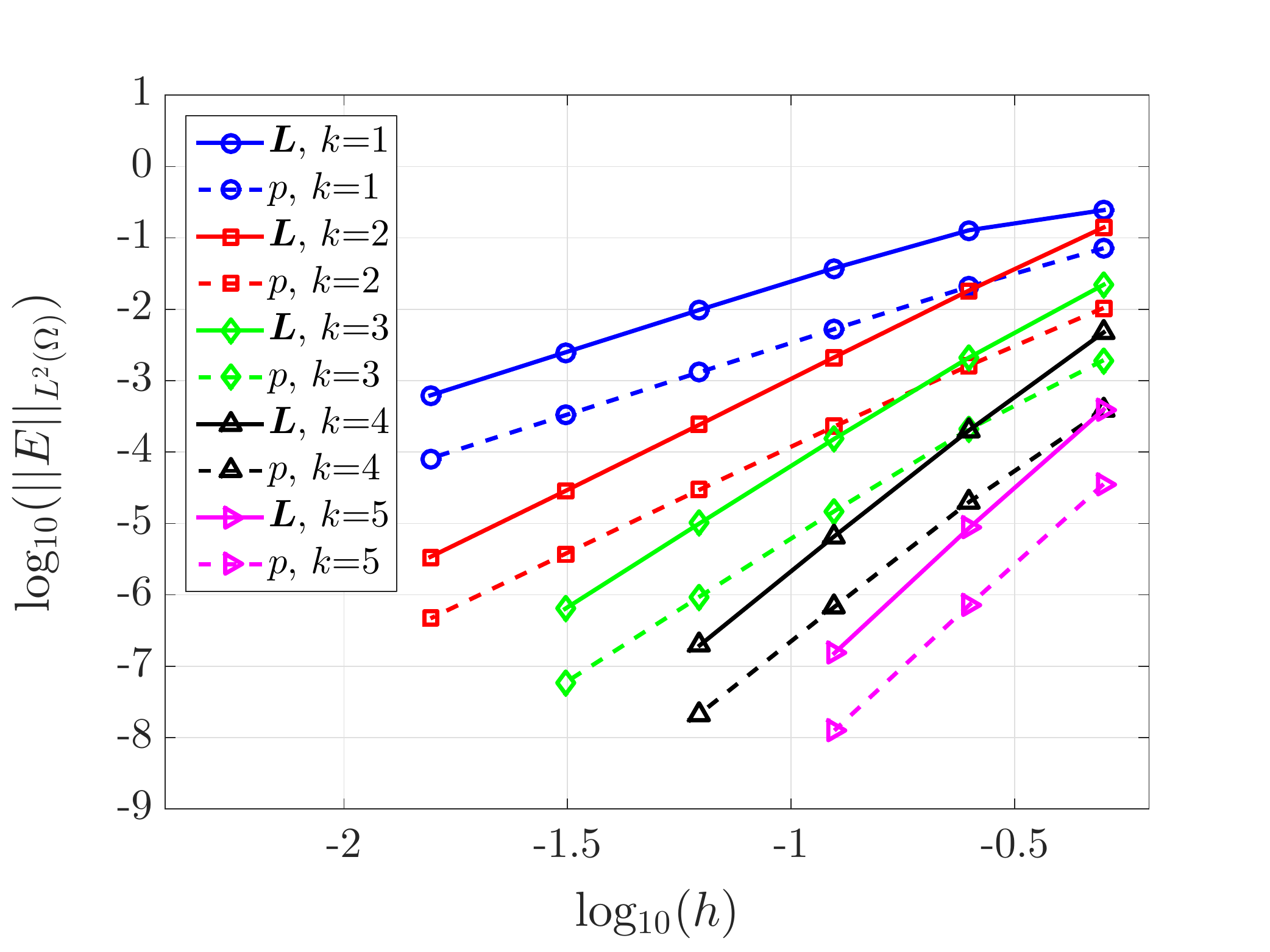}}
	\caption{Error of the dual variable $\bL$ and the pressure $p$ in the $\eltwo(\Omega_e)$ norm for different degrees of approximation in each element.}
	\label{fig:circleHConvP}
\end{figure}

The results of the $h$-convergence study are presented in Figures~\ref{fig:circleHConvU} and \ref{fig:circleHConvP}. Figure~\ref{fig:circleHConvU} shows the error of the solution $\bu$ and the postprocessed solution $\bu^\star$ in the $\eltwo(\Omega_e)$ norm for isoparametric and NEFEM elements for $k=1, \ldots, 5$. 
The optimal rate of convergence is obtained in all cases for both the solution $\bu$ (rate $k+2$) and the postprocessed solution $\bu^\star$ (rate $k+3$).

In Figure~\ref{fig:circleHConvP} a similar analysis is conducted, but now the error is measured for the dual variable $\bL$ and the pressure $p$, also in the $\eltwo(\Omega_e)$ norm and for isoparametric and NEFEM elements. 
Again, the optimal rate of convergence is obtained in all cases for both $\bL$ and $p$ (rate $k+2$).

%--------------------------------------------------------------------------
\section{Comparison of degree adaptivity strategies}  \label{sc:comparison}
%--------------------------------------------------------------------------

\begin{figure}[!b]
	\centering
	\includegraphics[width=0.32\textwidth]{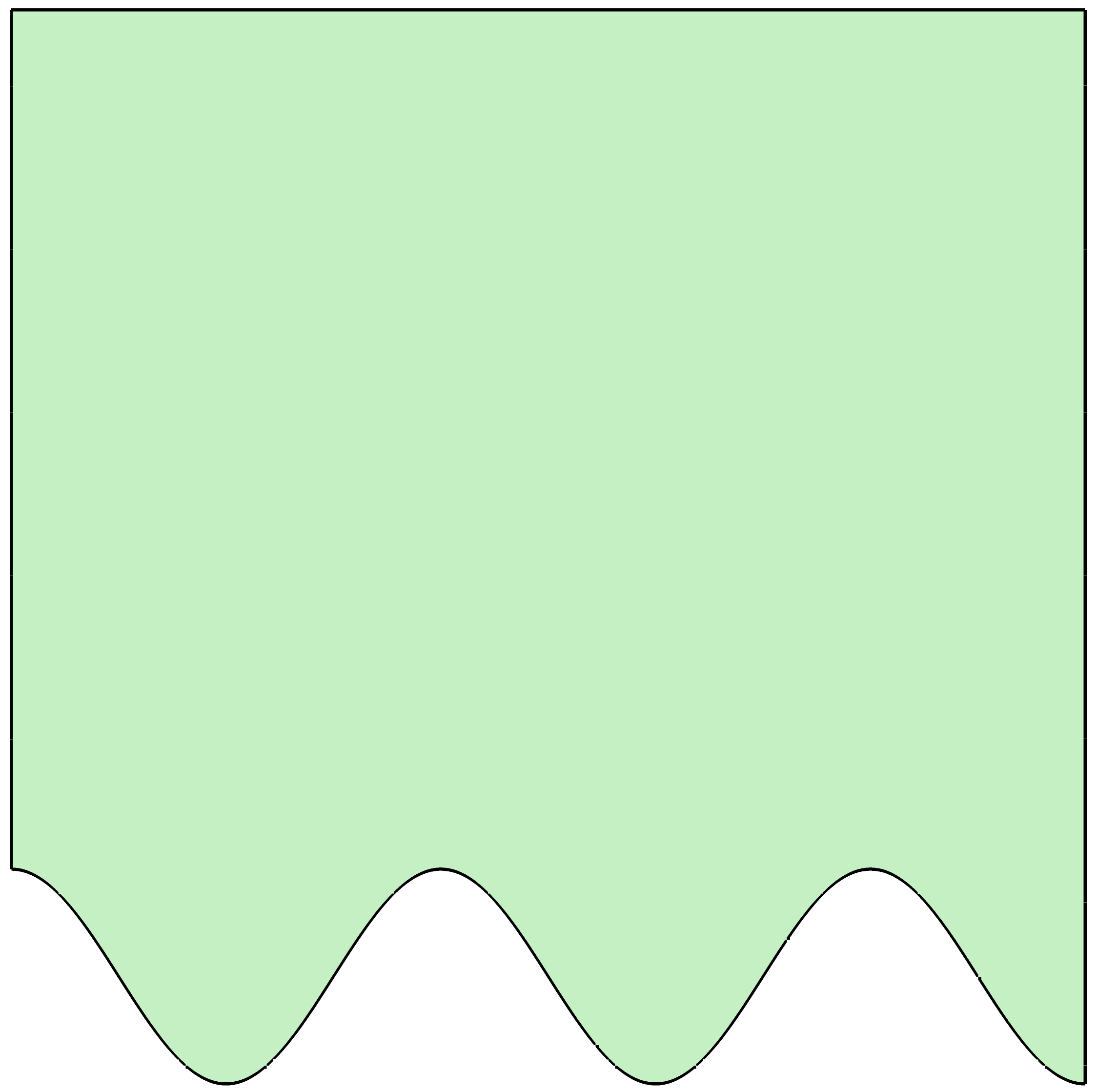}
	\caption{Computational domain for the test problem used to compare the different geometry update options in a degree adaptive process.}
	\label{fig:comparison_Domain}
\end{figure}

The same model problem employed in the previous example is utilised to compare the strategies described in Section~\ref{sc:geometricUpdate} to update the geometry during a degree adaptive process. The computational domain selected, shown in Figure~\ref{fig:comparison_Domain}, features an oscillatory boundary and represents a common problem encountered in biological transport applications, see for instance~\cite{scholle2004creeping}. \RS{More precisely, the curved part of the boundary is given by the curve $f(x) = (1 + \cos(5\pi x))/10$.}

 Dirichlet boundary conditions are imposed on the polygonal part of the boundary whereas a Neumann boundary condition, corresponding to the exact traction derived from the solution in Equation~\eqref{eq:analyticalStokes} is imposed on the oscillatory part of the boundary.

%--------------------------------------------------------------------------
\subsection{No geometric update}  \label{sc:comparison0}
%--------------------------------------------------------------------------

\begin{figure}[!bt]	
	\centering
	\subfigure[Degree]			{\includegraphics[width=0.32\textwidth]{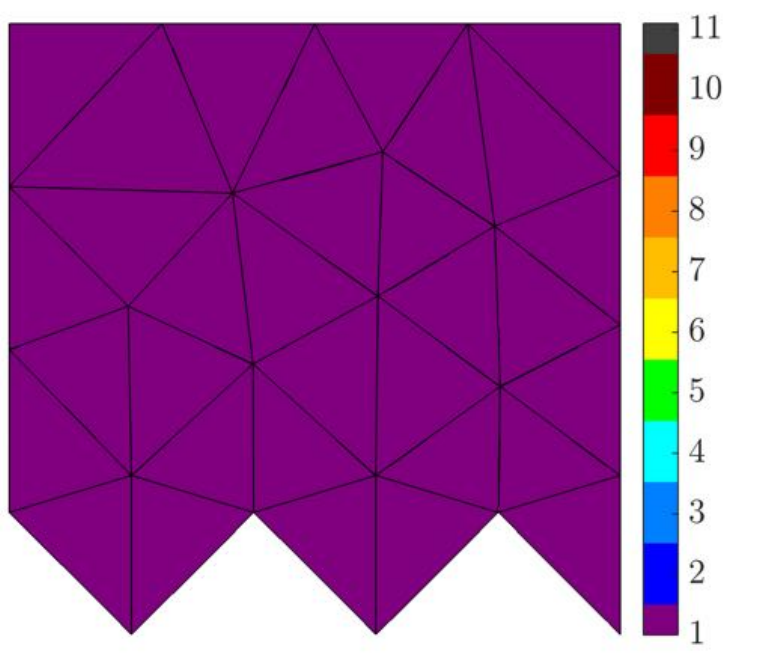}}
	\subfigure[Estimated error]	{\includegraphics[width=0.32\textwidth]{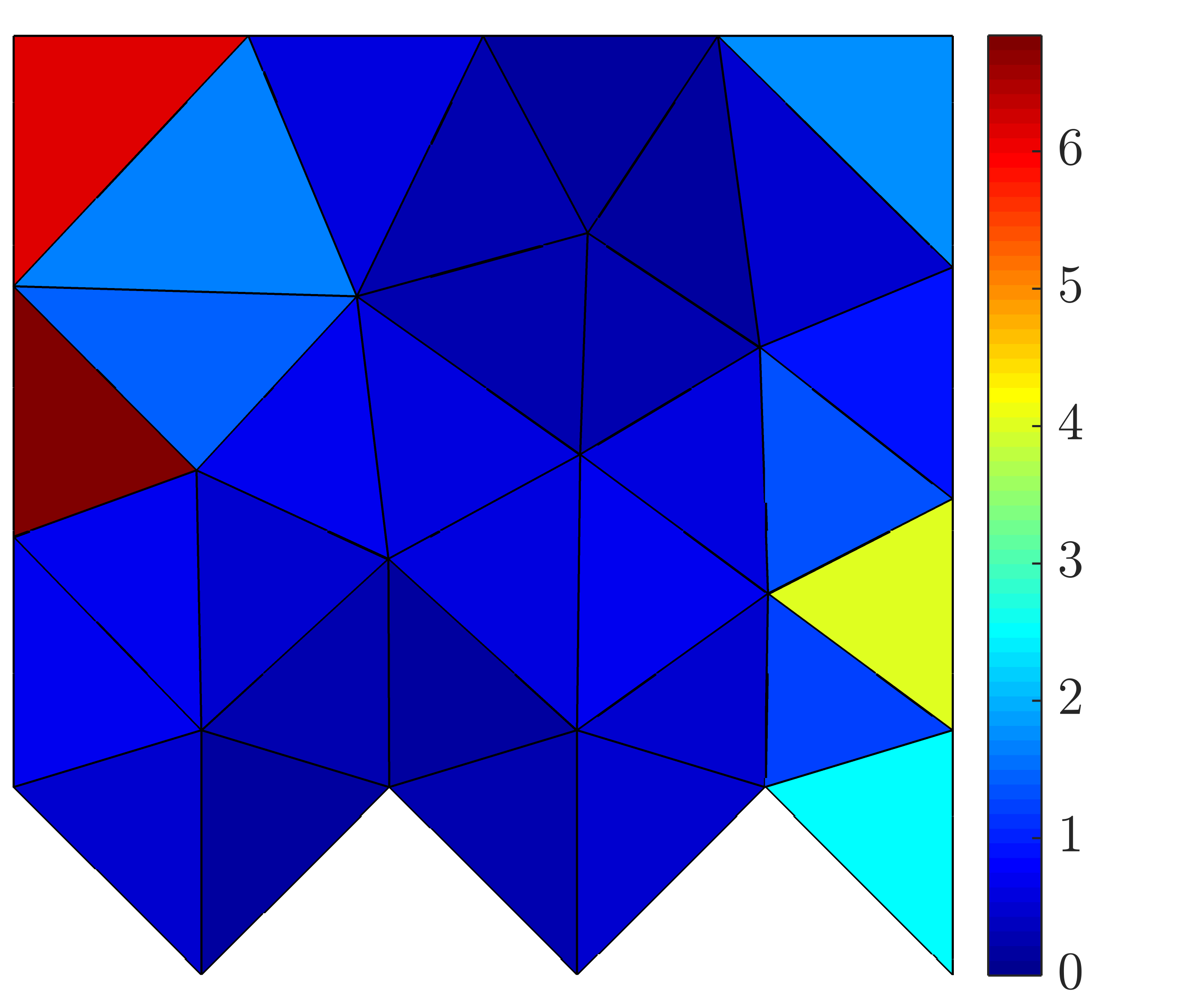}}
	\subfigure[Exact Error]		{\includegraphics[width=0.32\textwidth]{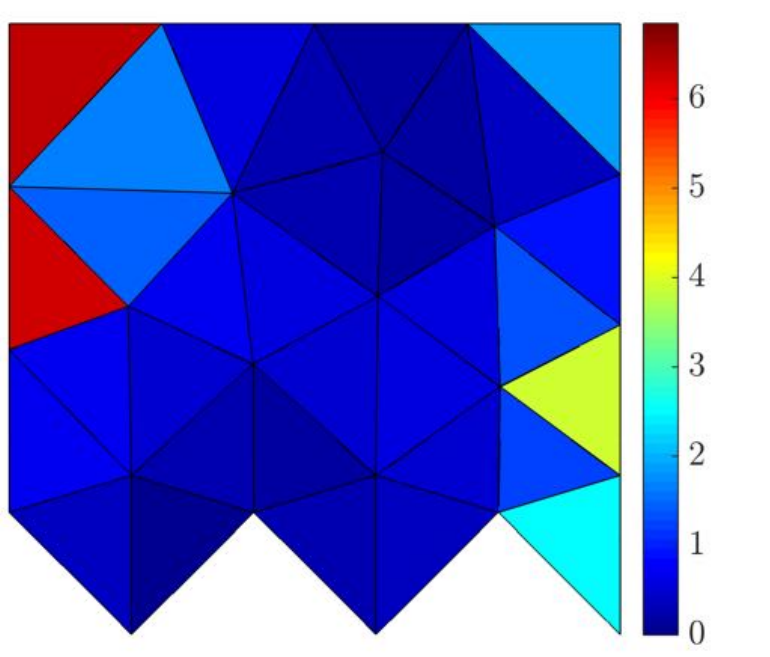}}
	\caption{First iteration of the degree adaptivity procedure with HDG isoparametric elements.}
	\label{fig:comparisonFEM0_Iteration1}
\end{figure}

First, the degree adaptive process with no communication with the CAD model is studied, as illustrated in Figure~\ref{fig:pAdaptivityFEM0}. The process starts with a degree of approximation $k=1$ in all elements. At each iteration the degree of the functional approximation is adapted according to the strategy presented in Section~\ref{sc:adaptivity} whereas a linear approximation of the geometry is kept irrespectively of the degree of the functional approximation. Figure~\ref{fig:comparisonFEM0_Iteration1} shows the original mesh, the estimated error and the exact error, computed by using the known analytical solution. The $\eltwo$ norm of the error is represented as a constant value in each element, showing a good agreement between the estimated and the exact error.

\begin{figure}[!bt]	
	\centering
	\subfigure[Degree]			{\includegraphics[width=0.32\textwidth]{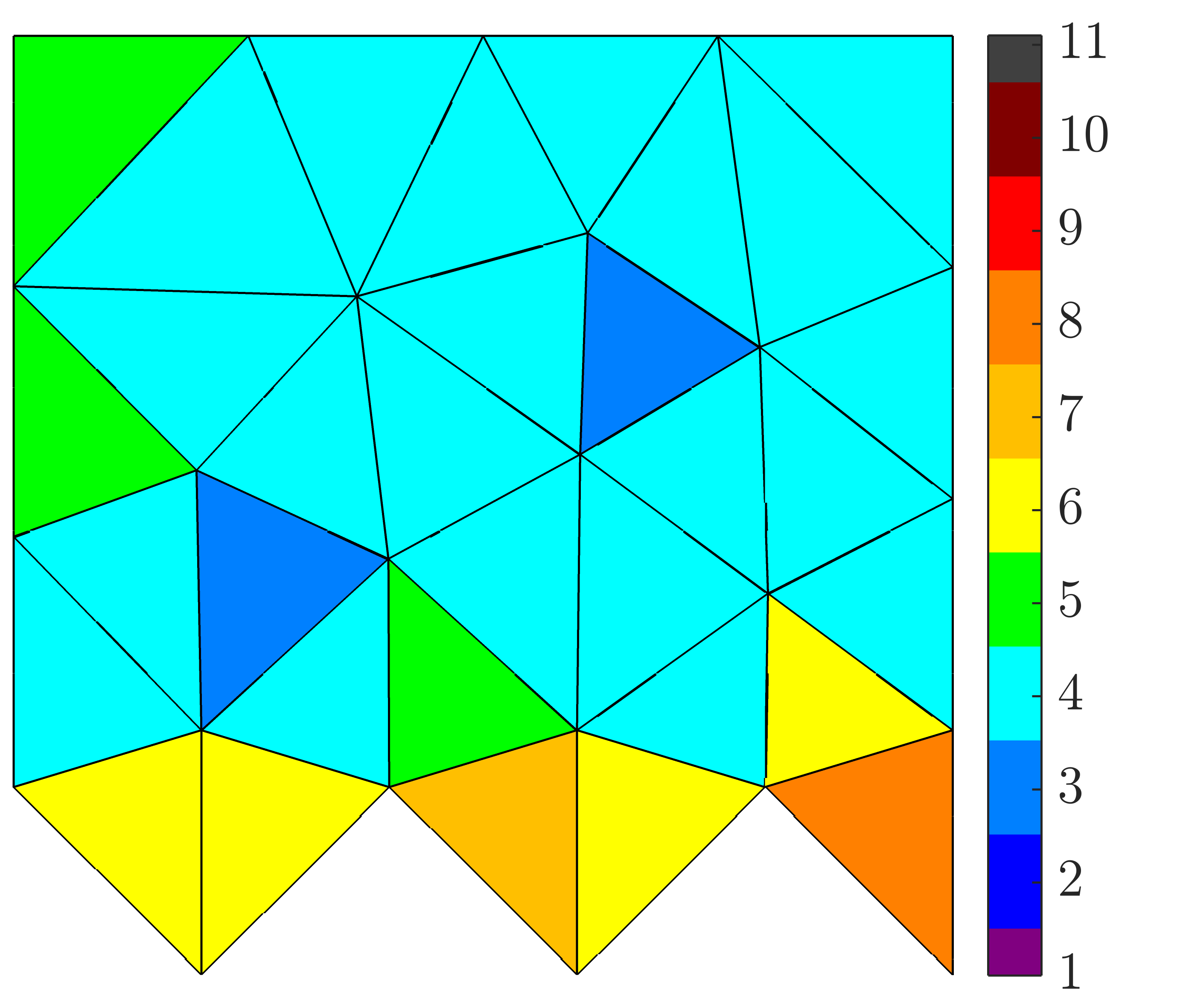}}
	\subfigure[Estimated error]	{\includegraphics[width=0.32\textwidth]{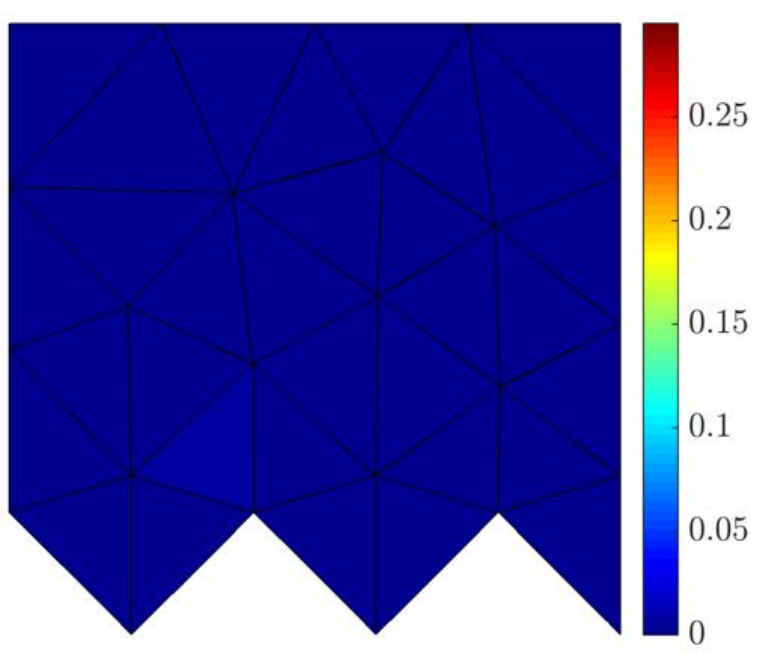}}
	\subfigure[Exact Error]		{\includegraphics[width=0.32\textwidth]{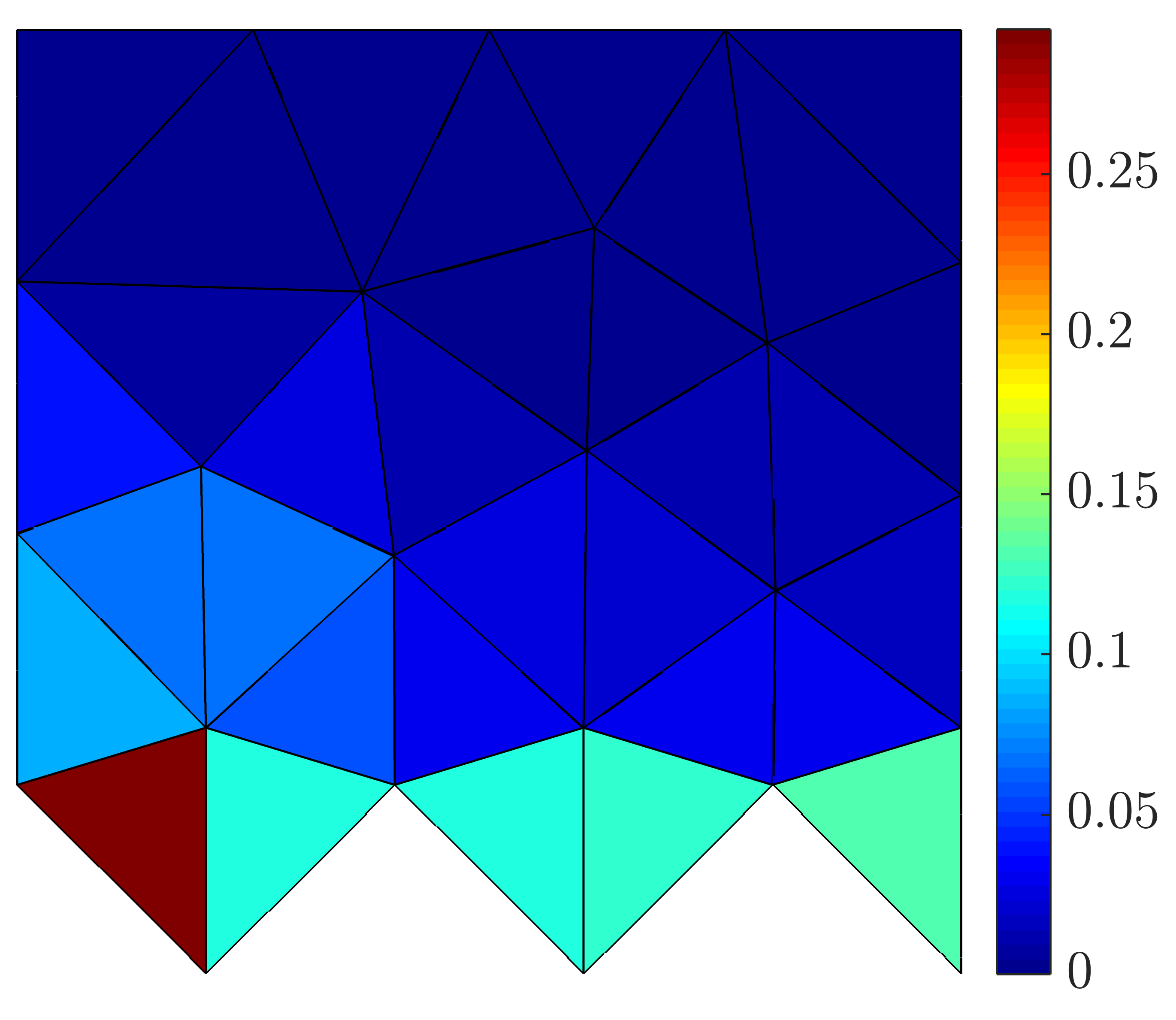}}
	\caption{Sixth iteration of the degree adaptivity procedure with no geometric update using HDG isoparametric elements and $q=1$.}
	\label{fig:comparisonFEM0_Iteration6}
\end{figure}

After six iterations of the adaptivity process, the degree of approximation is adapted in each element as shown in Figure~\ref{fig:comparisonFEM0_Iteration6} (a) but a linear geometric approximation of the curved boundary is still considered. 
The estimated error in each element, shown in Figure~\ref{fig:comparisonFEM0_Iteration6} (b), is below the desired error which is $0.5 \times 10^{-2}$ in this example but the computation of the exact error, shown in Figure~\ref{fig:comparisonFEM0_Iteration6} (c), reveals a significant disparity when compared to the estimated error.

\begin{figure}[!ht]	
	\centering
	\subfigure[$q$=1]{\includegraphics[width=0.49\textwidth]{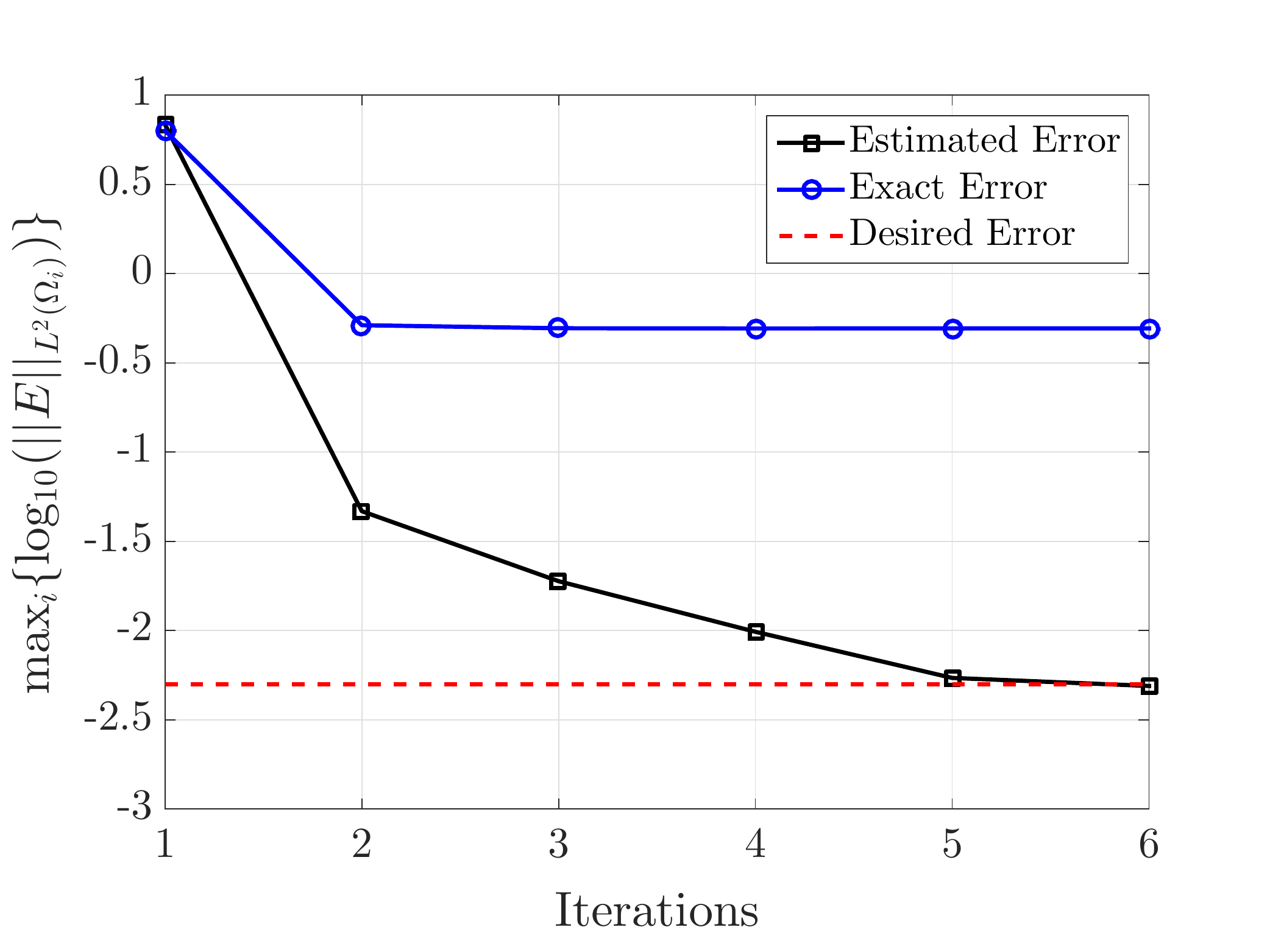}}
	\subfigure[$q$=2]{\includegraphics[width=0.49\textwidth]{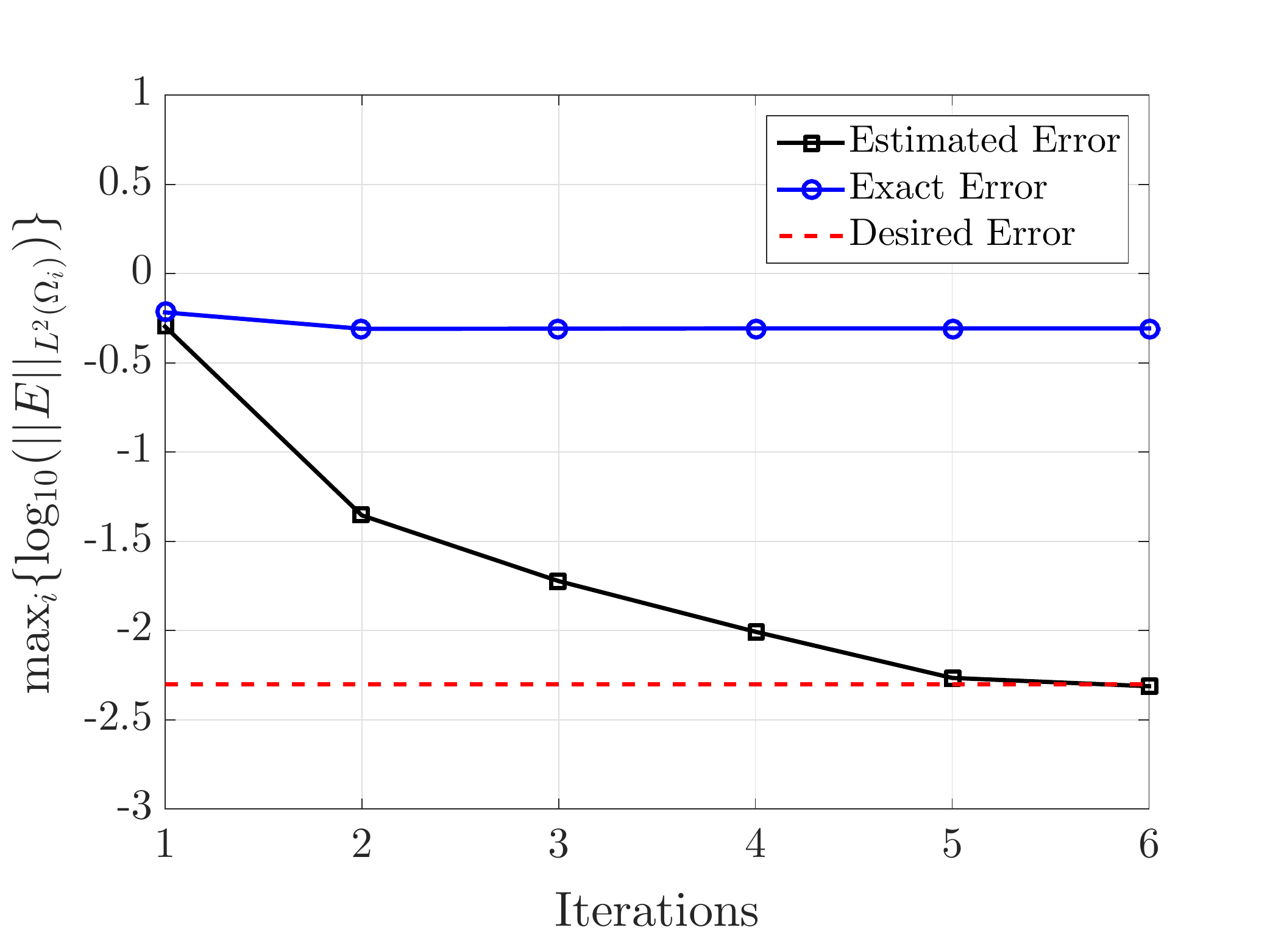}}\\
	\subfigure[$q$=3]{\includegraphics[width=0.49\textwidth]{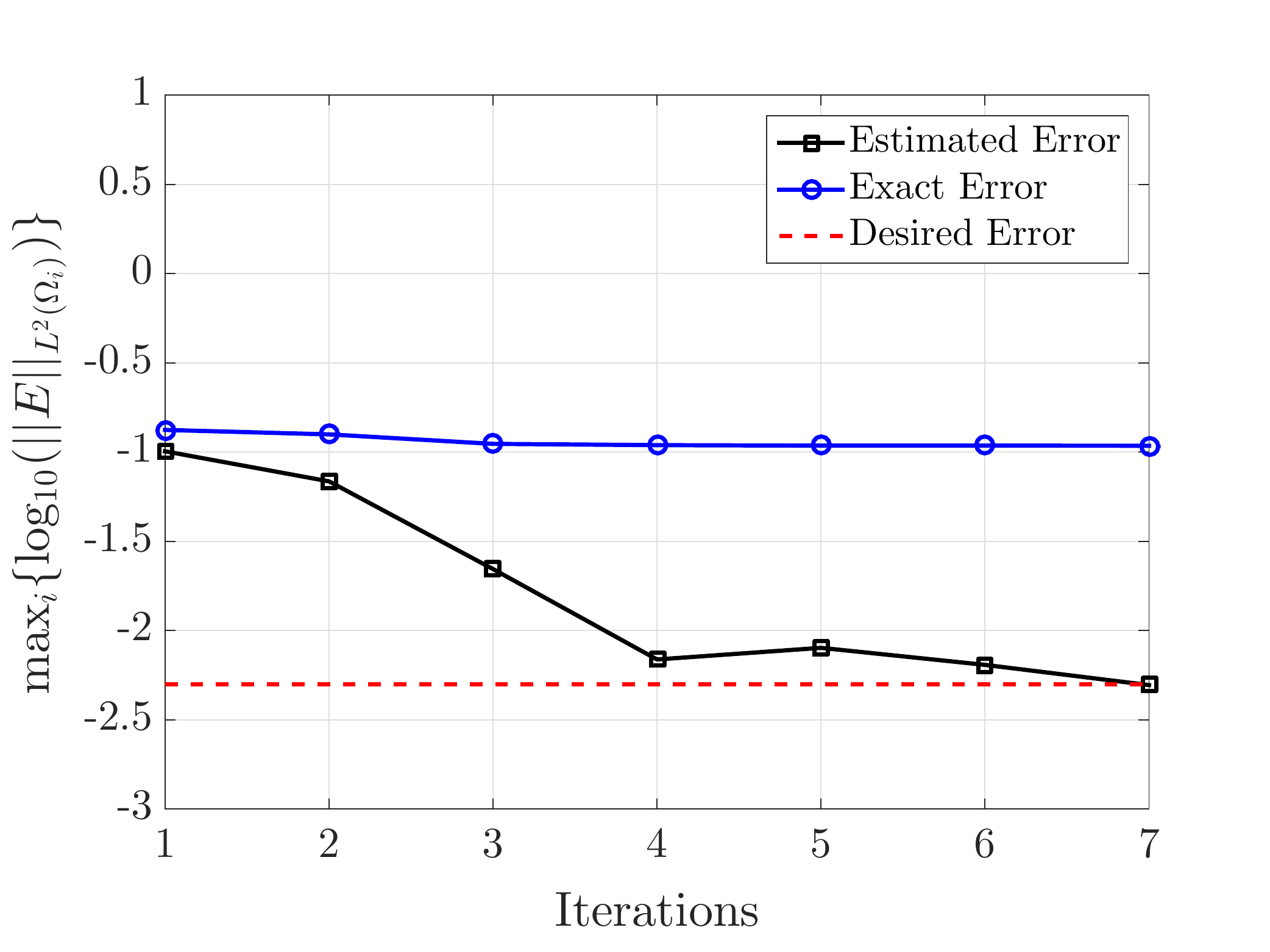}}
	\subfigure[$q$=4]{\includegraphics[width=0.49\textwidth]{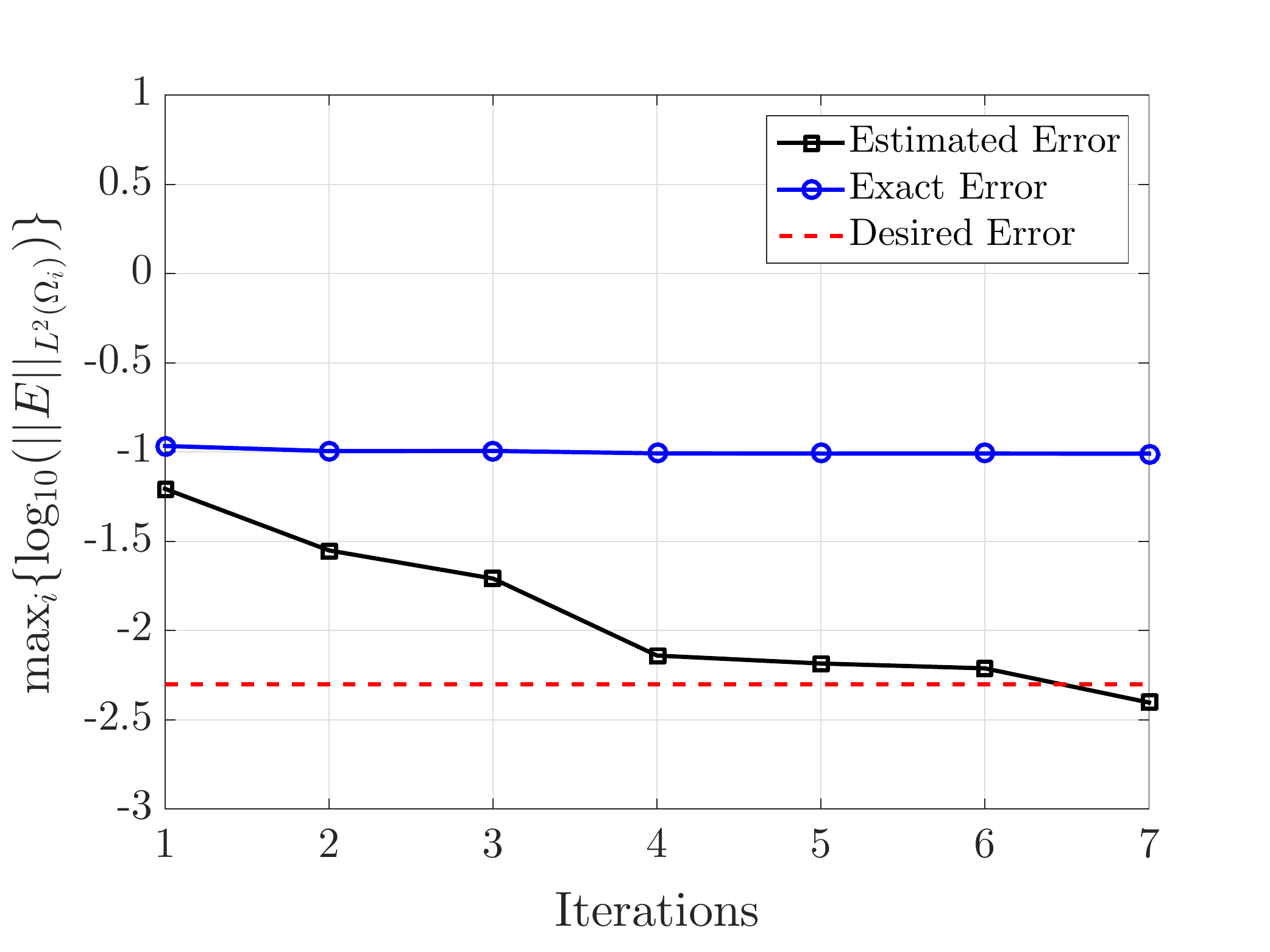}}
	\caption{Evolution of the estimated and exact errors during a degree adaptivity process for different degrees of the polynomials used to approximate the geometry ($q$).}
	\label{fig:adaptivityErrorFEM0}
\end{figure}

To better analyse the results, Figure~\ref{fig:adaptivityErrorFEM0}  shows the evolution of the maximum estimated error in each element and the maximum exact in each element for different geometric approximations of the curved boundary. Figure~\ref{fig:adaptivityErrorFEM0} (a) corresponds to the case illustrated in Figures~\ref{fig:comparisonFEM0_Iteration1} and \ref{fig:comparisonFEM0_Iteration6}, where a linear approximation of the geometry is considered ($q=1$). It can be clearly observed that, as the degree adaptive process evolves, the difference between estimated and exact error becomes more and more sizeable. In the sixth iteration, the adaptive process shows convergence, with an estimated error less than the desired error $0.5 \times 10^{-2}$ but this is two orders of magnitude lower than the exact error, clearly indicating that the error estimator is not reliable because the estimator assumes that the geometry is exactly represented.

As a linear approximation of the geometry is well known to be not suitable when high order functional approximations are considered~\cite{Bassi-BR:97,cirak2000subdivision,sevilla2011nurbs,soghrati2016nurbs}, the same experiment is repeated by using a more accurate boundary representation. The plots in Figure~\ref{fig:adaptivityErrorFEM0} (b), (c) and (d) show the evolution of the maximum estimated error in each element and the maximum exact in each element for quadratic, cubic and quartic approximation of the geometry. In all cases it is clearly observed that the error estimator is not reliable because the adaptive process converges but the exact error is more than one order of magnitude higher than the desired error.

\begin{figure}[!tb]	
	\centering
	\subfigure[Degree]			{\includegraphics[width=0.32\textwidth]{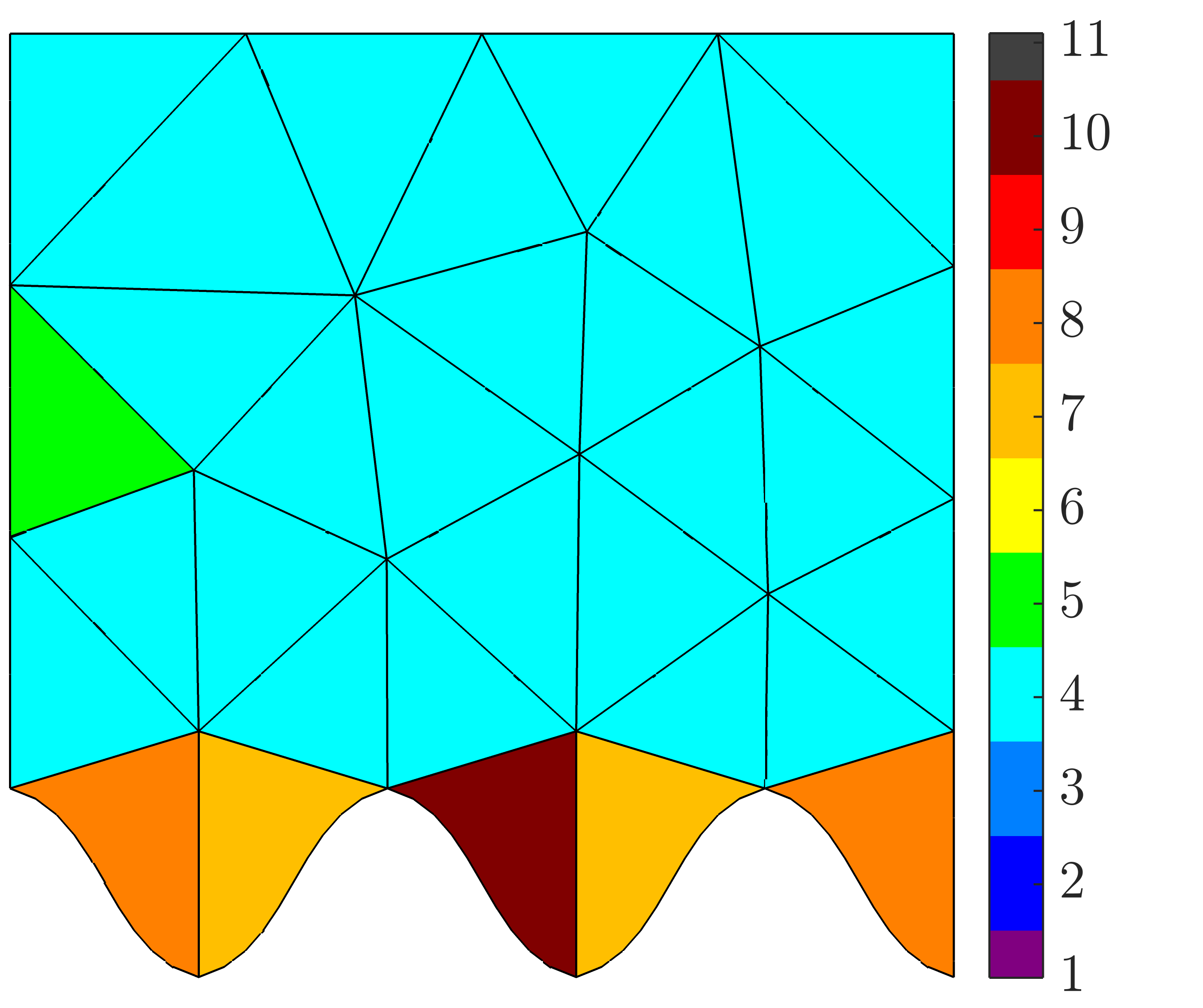}}
	\subfigure[Estimated error]	{\includegraphics[width=0.32\textwidth]{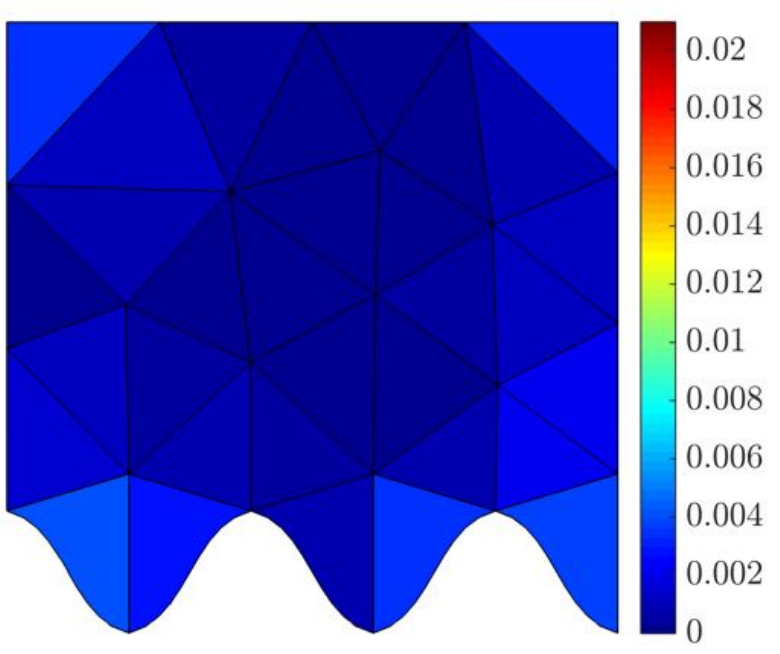}}
	\subfigure[Exact Error]		{\includegraphics[width=0.32\textwidth]{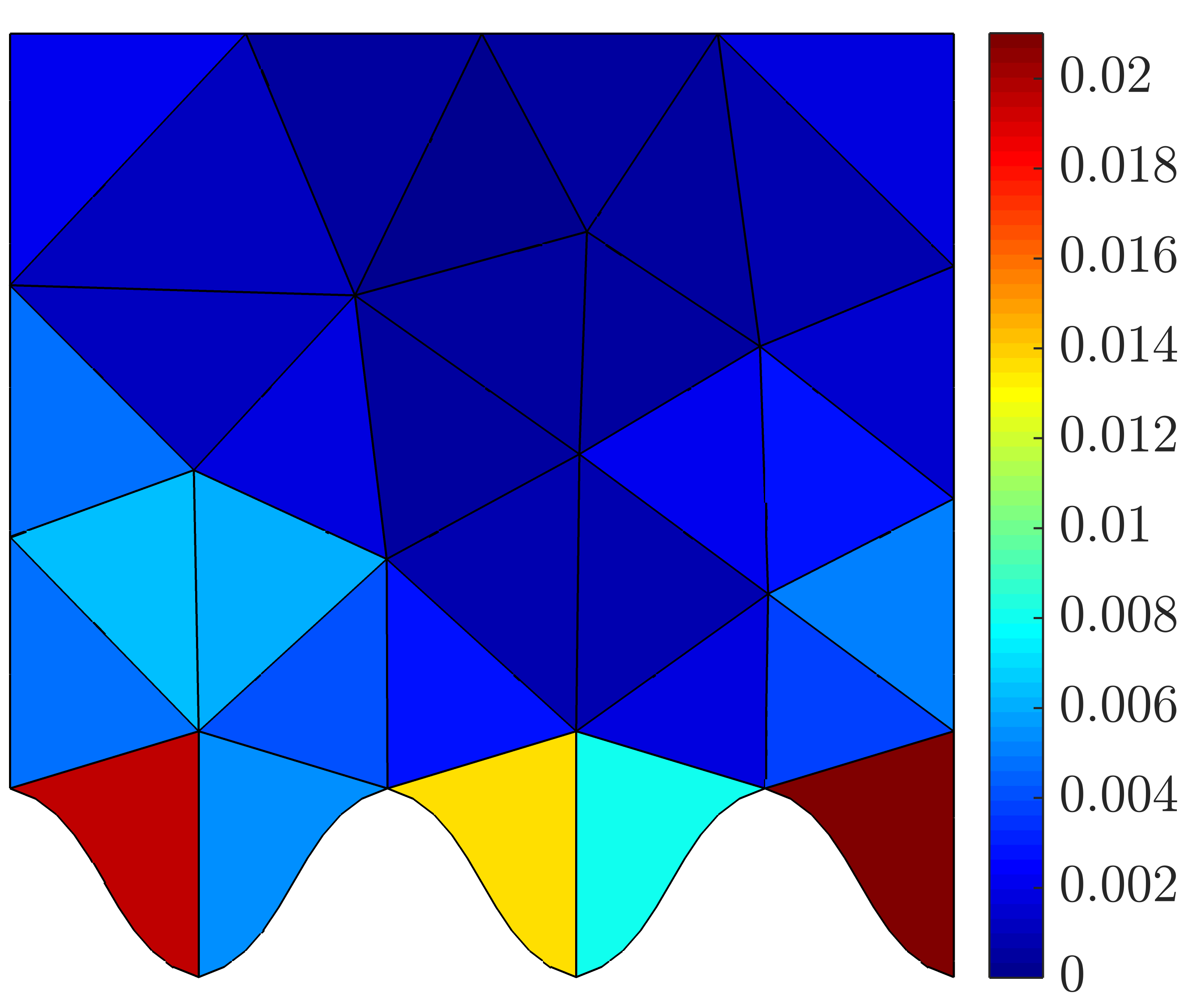}}
	\caption{Sixth iteration of the degree adaptivity procedure with no geometric update using HDG isoparametric elements and $q=4$.}
	\label{fig:comparisonFEM0P4_Iteration6}
\end{figure}

The degree of approximation, estimated error and exact error obtained in the last iteration of the adaptive process for $q=4$ is represented in Figure~\ref{fig:comparisonFEM0P4_Iteration6}.
The results show that even with a more accurate geometric approximation, the exact error in the elements close to the curved boundary is much higher than the estimated error. There is clear evidence that, if no communication with a CAD model is undertaken during the degree adaptive process, the original mesh must be pre-adapted manually in order to ensure that the geometric error is small enough in order to ensure that the error estimator is reliable, clearly compromising the robustness of the whole adaptivity process.

%--------------------------------------------------------------------------
\subsection{NEFEM HDG}  \label{sc:comparison2}
%--------------------------------------------------------------------------

The strategy proposed in this work consists of utilising NEFEM, where the geometry is always given by its CAD boundary representation, irrespective of the degree of the functional approximation. In the context of a degree adaptive process, this means that no communication the CAD model is required as the exact boundary representation is already used by the NEFEM solver. 

The process starts with a degree of approximation $k=1$ in all elements. At each iteration the degree of the functional approximation is adapted according to the strategy presented in Section~\ref{sc:adaptivity} and a new nodal distribution is generated for each curved element. 

\begin{figure}[!b]	
	\centering
	\subfigure[Degree]			{\includegraphics[width=0.32\textwidth]{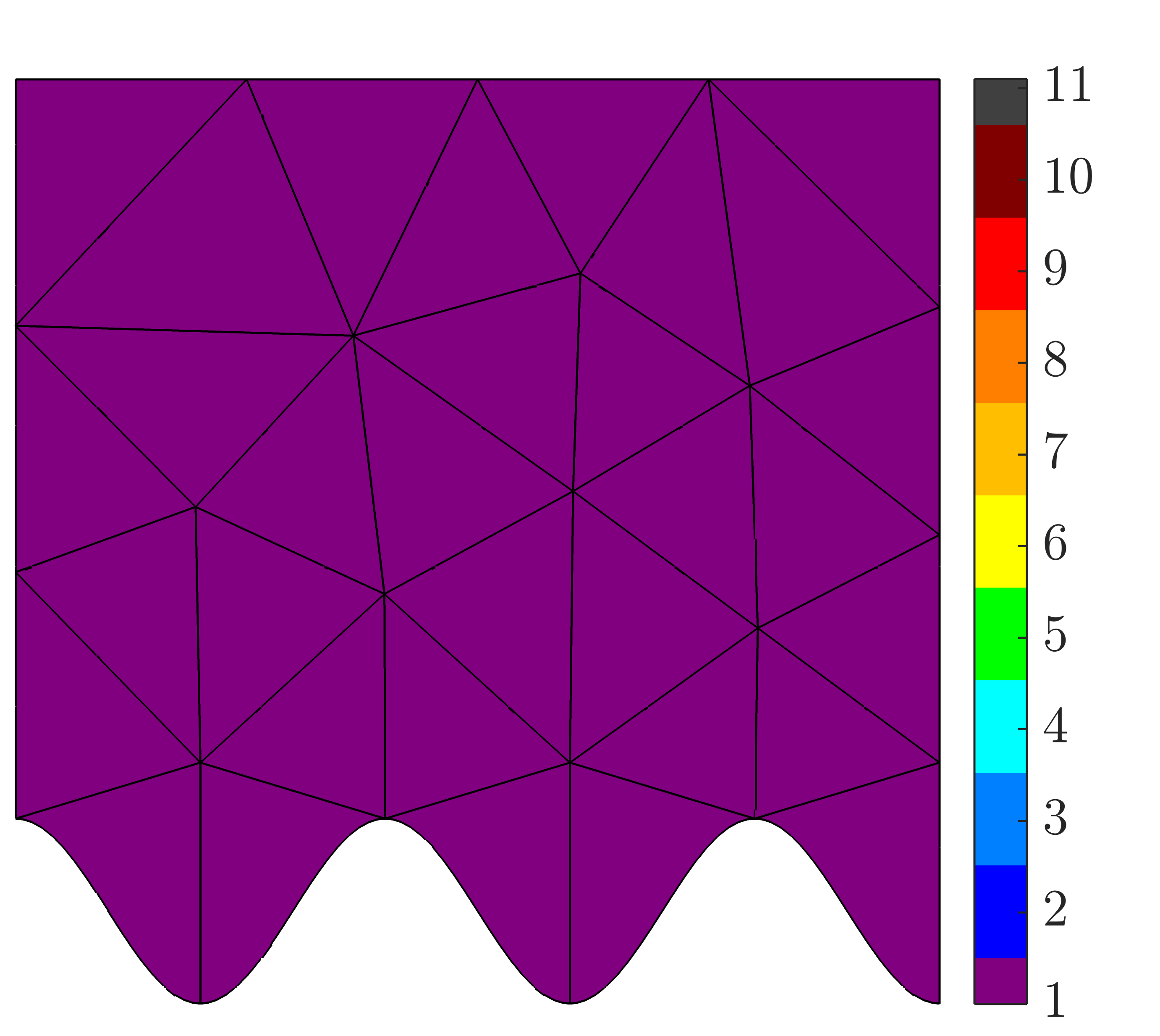}}
	\subfigure[Estimated error]	{\includegraphics[width=0.32\textwidth]{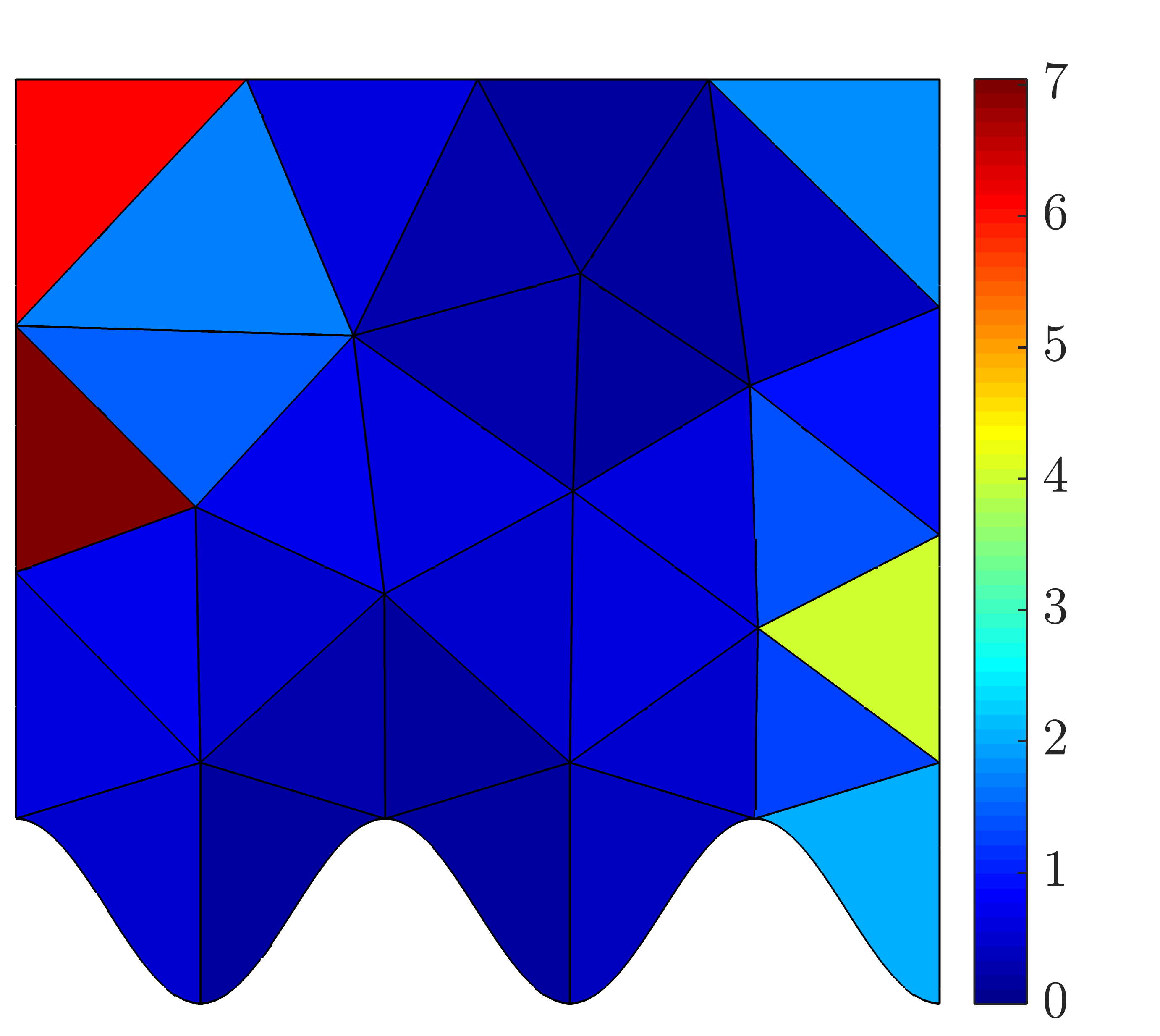}}
	\subfigure[Exact Error]		{\includegraphics[width=0.32\textwidth]{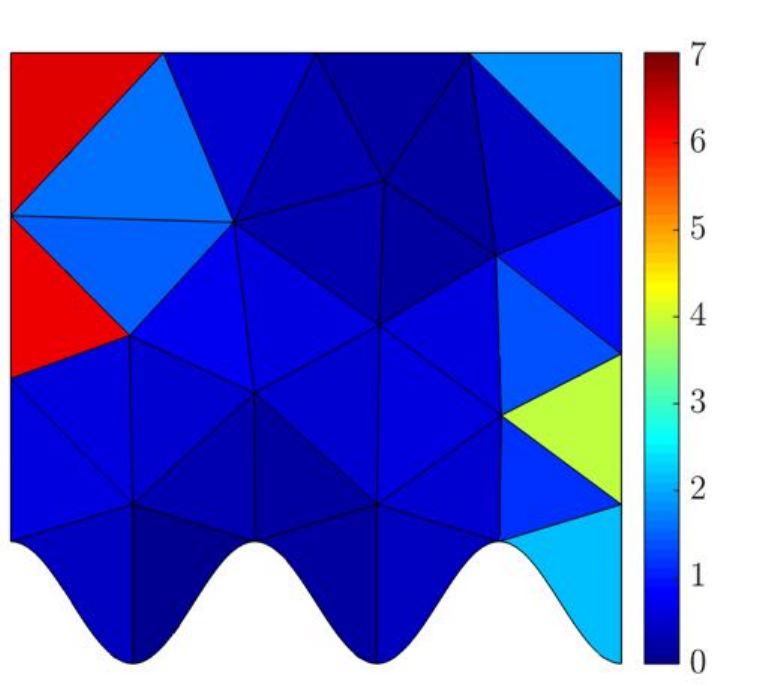}}
	\caption{First iteration of the degree adaptivity procedure with geometric update using HDG NEFEM elements.}
	\label{fig:comparisonNEFEM_Iteration1}
\end{figure}

Figure~\ref{fig:comparisonNEFEM_Iteration1} shows the original mesh, the estimated error and the exact error, computed by using the known analytical solution. 
It is worth emphasising that, even when the degree of the functional approximation used is linear ($k=1$) the exact boundary representation is considered, as shown in Figure~\ref{fig:comparisonNEFEM_Iteration1} (a). The results show a very similar distribution for the estimated and exact errors. 

\begin{figure}[!bt]	
	\centering
	\subfigure[Degree]			{\includegraphics[width=0.32\textwidth]{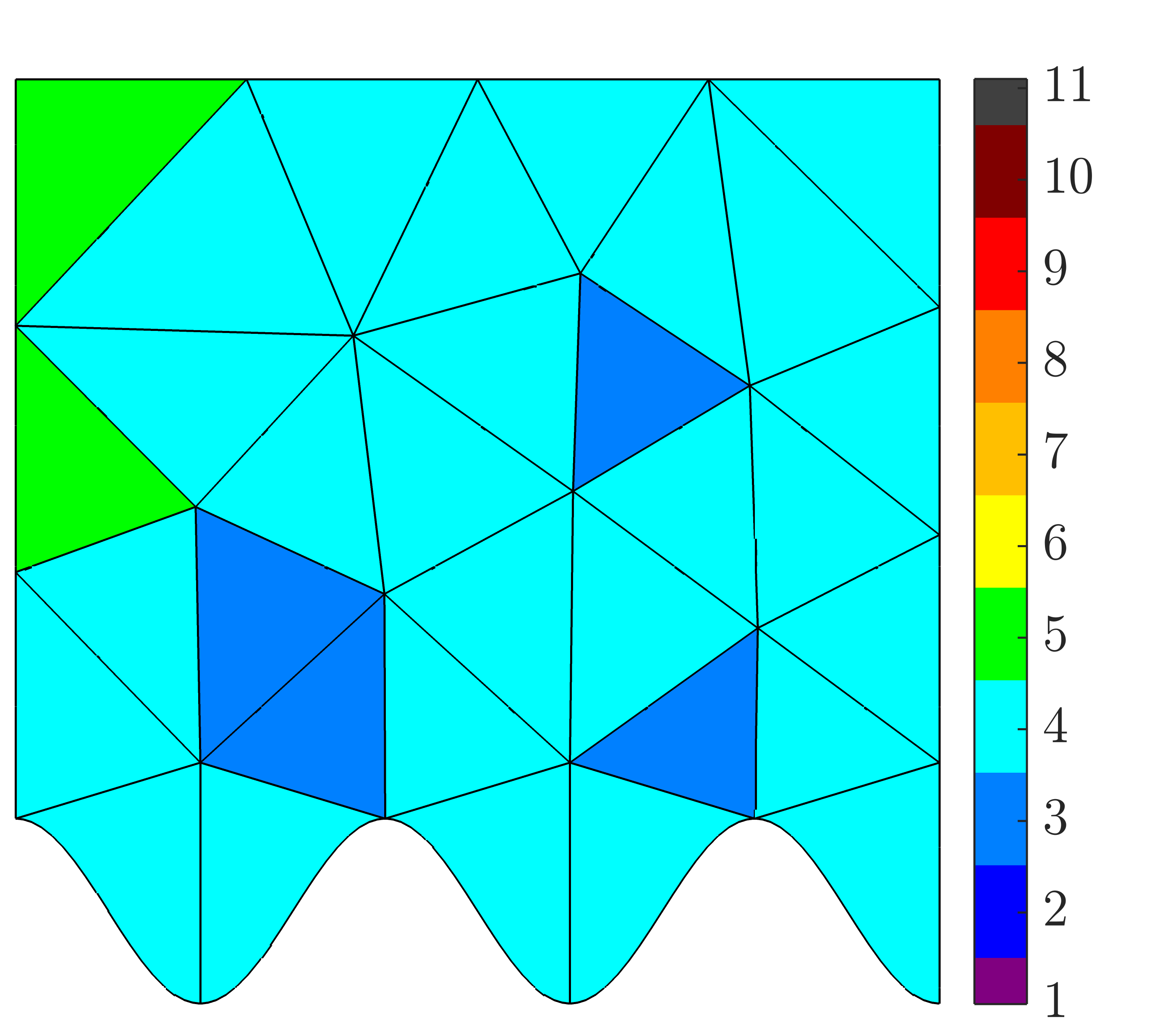}}
	\subfigure[Estimated error]	{\includegraphics[width=0.32\textwidth]{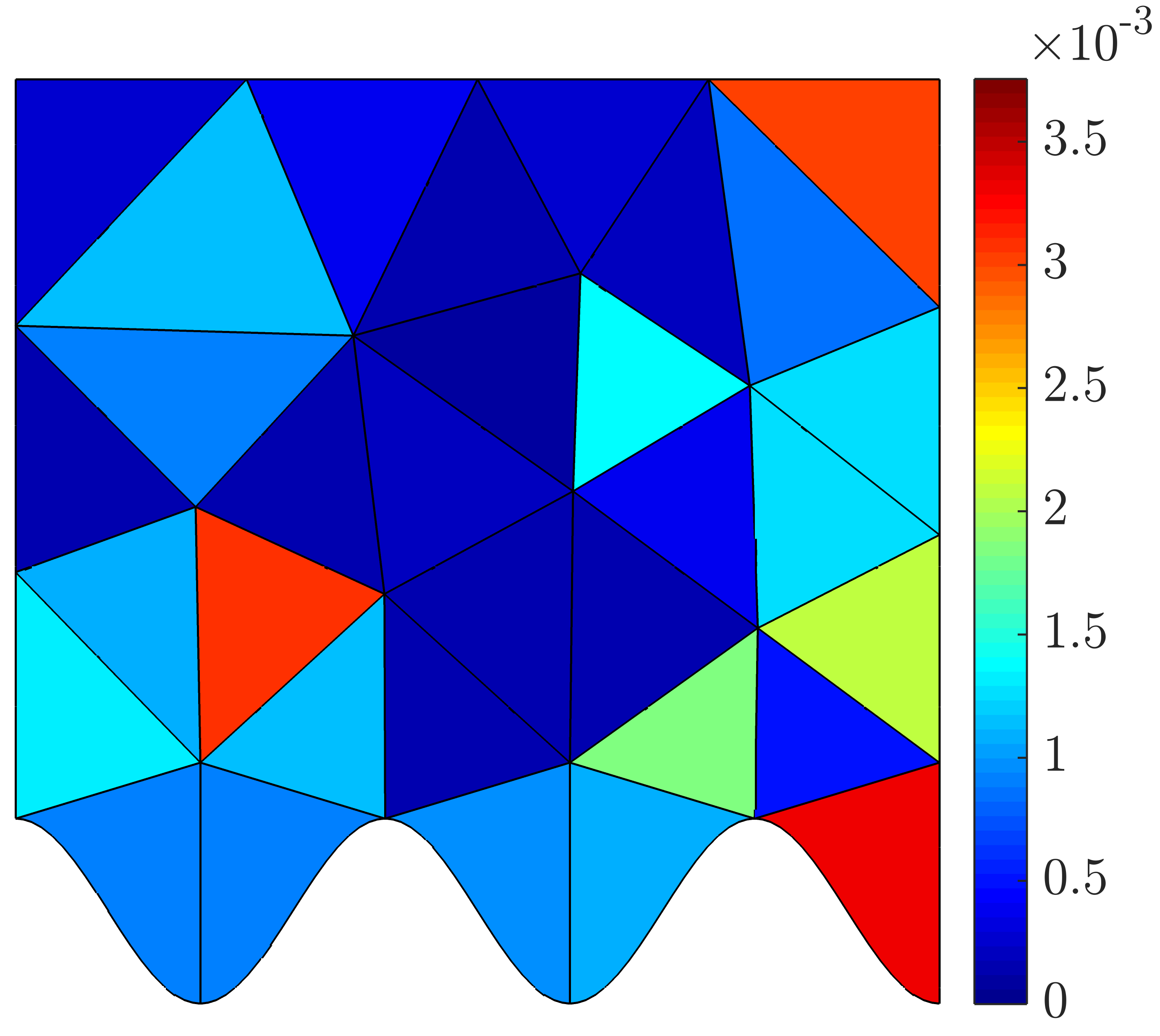}}
	\subfigure[Exact Error]		{\includegraphics[width=0.32\textwidth]{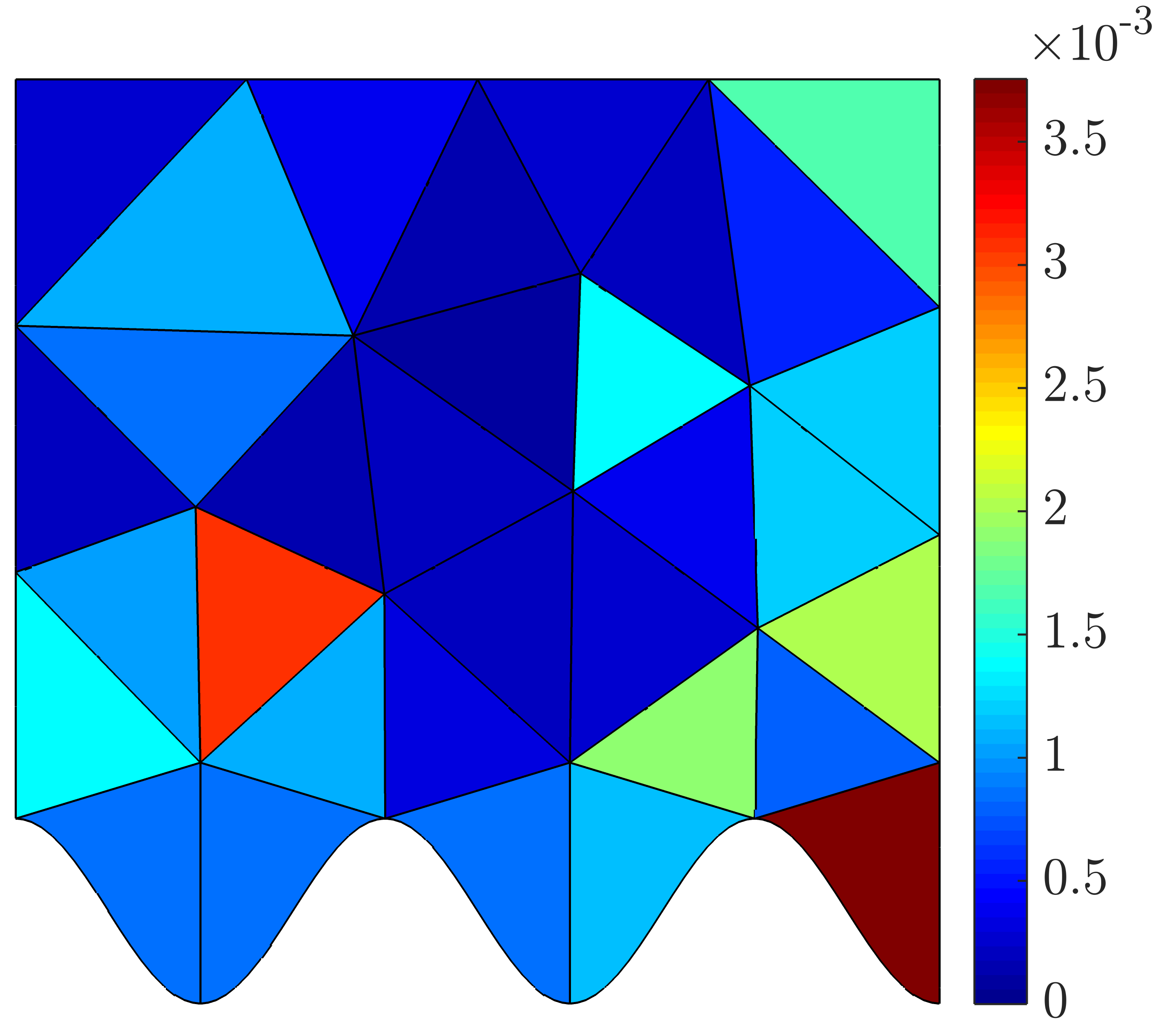}}
	\caption{Third iteration of the degree adaptivity procedure with geometric update using HDG NEFEM elements.}
	\label{fig:comparisonNEFEM_Iteration3}
\end{figure}

In this case, after only three iterations of the degree adaptive process convergence is achieved. The degree of approximation used in each element, the estimated and the exact errors in each element are represented in Figure~\ref{fig:comparisonNEFEM_Iteration3}.

\begin{remark}
	As discussed in Section~\ref{sc:geometricUpdate}, an alternative, not employed in practice due to the high cost induced by the re-generation of the mesh at each iteration of the adaptive process, consists of changing both the degree of approximation for the solution and for the geometry during the adaptivity process, as illustrated in Figure~\ref{fig:pAdaptivityFEM1eMeshes}. 
\end{remark}

\begin{figure}[!tb]	
	\centering
	\subfigure[Desired error $0.5 \times 10^{-2}$]{\includegraphics[width=0.49\textwidth]{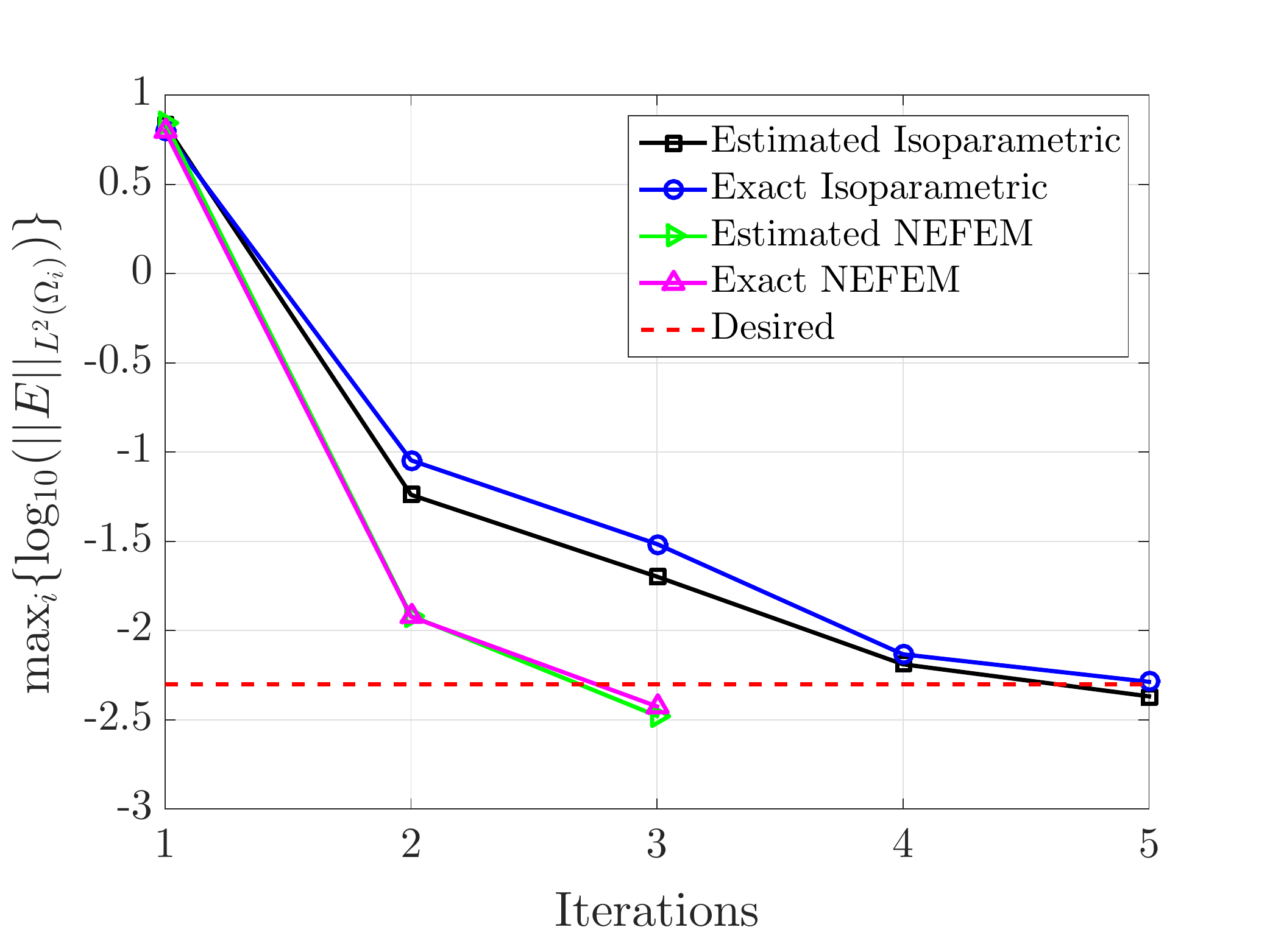}}
	\subfigure[Desired error $0.5 \times 10^{-3}$]{\includegraphics[width=0.49\textwidth]{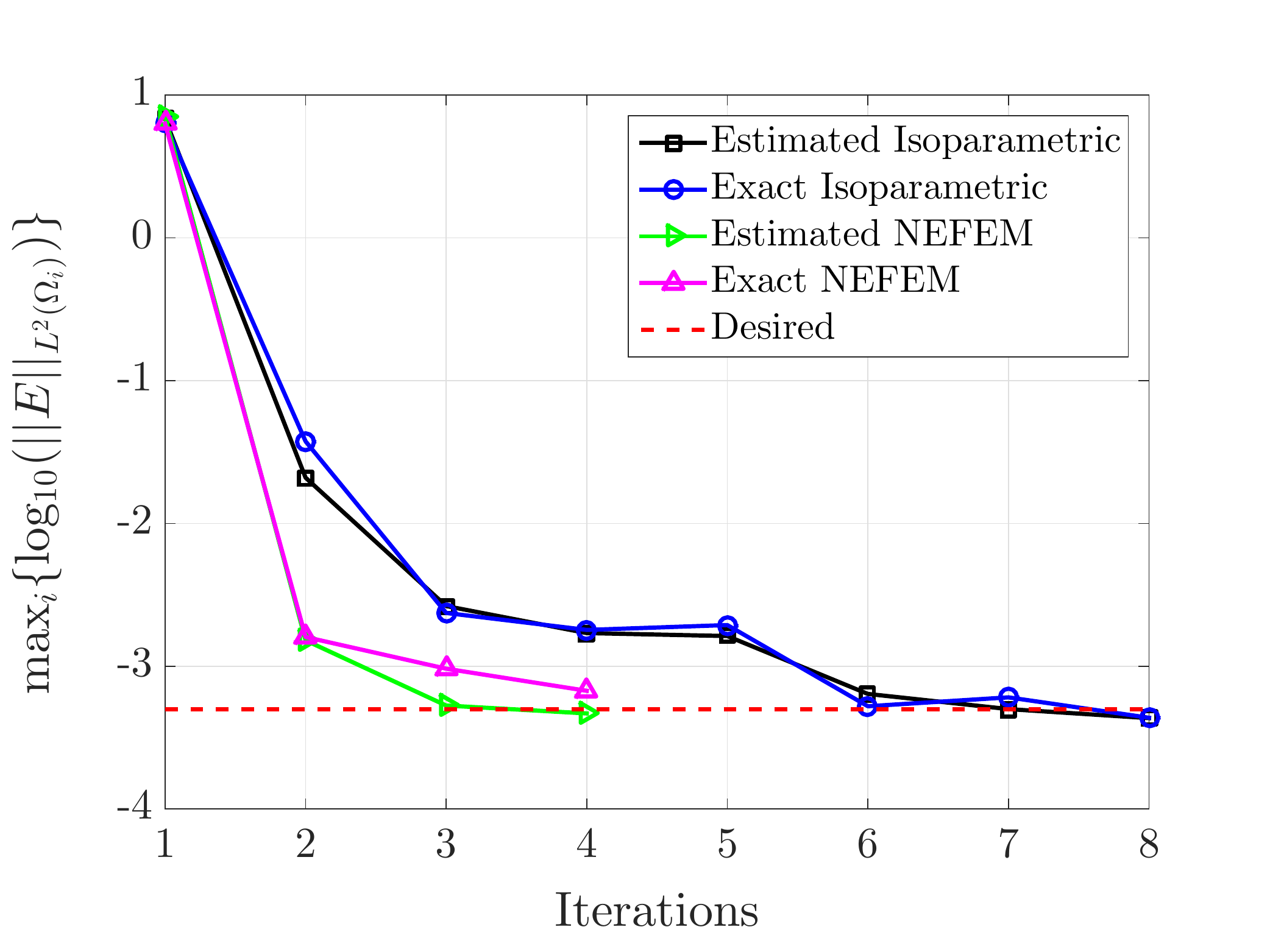}}
	\caption{Estimated and exact errors for isoparametric and NEFEM.}
	\label{fig:adaptivityErrorFEMNEFEM}
\end{figure}

To illustrate the superiority of NEFEM, Figure~\ref{fig:adaptivityErrorFEMNEFEM}  shows the evolution of the maximum estimated error in each element and the maximum exact error in each element for isoparametric and NEFEM approaches and for two magnitudes of the desired error. 
Figure~\ref{fig:adaptivityErrorFEM0} (a) corresponds to the case previously illustrated, where the desired error is $0.5 \times 10^{-2}$, whereas Figure~\ref{fig:adaptivityErrorFEM0} (b) shows the same study but with a desired error of $0.5 \times 10^{-3}$. In both cases the superiority of NEFEM is clear as the desired error is achieved by substantially reducing the computational cost (i.e. the number of iterations of the degree adaptive process and the number of degrees of freedom required to achieve the required accuracy).

In addition, it is worth emphasising that the isoparametric approach requires communication with the CAD model in each iteration to re-generated the high-order nodal distribution. These nodal distributions in curved elements must be specifically designed to ensure optimal convergence of the isoparametric approach~\cite{Bernardi:89,Babuska-CB:95}, while for NEFEM the Cartesian approximation of the solution ensures that the accuracy of the approximation is much less sensitive to the quality of the nodal distribution.

%==========================================================================
\section{Numerical Examples} \label{sc:examples}
%==========================================================================

This section presents four numerical examples to illustrate the potential of NEFEM when combined with HDG to perform a degree adaptive process. The examples involve geometries with curved boundaries and where coarse meshes are considered to show the robustness of the proposed methodology. In all the examples the high-order isoparametric and NEFEM meshes are generated using the techniques described in~\cite{xie2013generation,poya2016unified} and~\cite{sevilla2016generation} respectively.

%==========================================================================
\subsection{Flow in a channel with randomly distributed ellipses} \label{sc:exampleEllipses}
%==========================================================================

\begin{figure}[!b]
	\centering
	\includegraphics[width=0.8\textwidth]{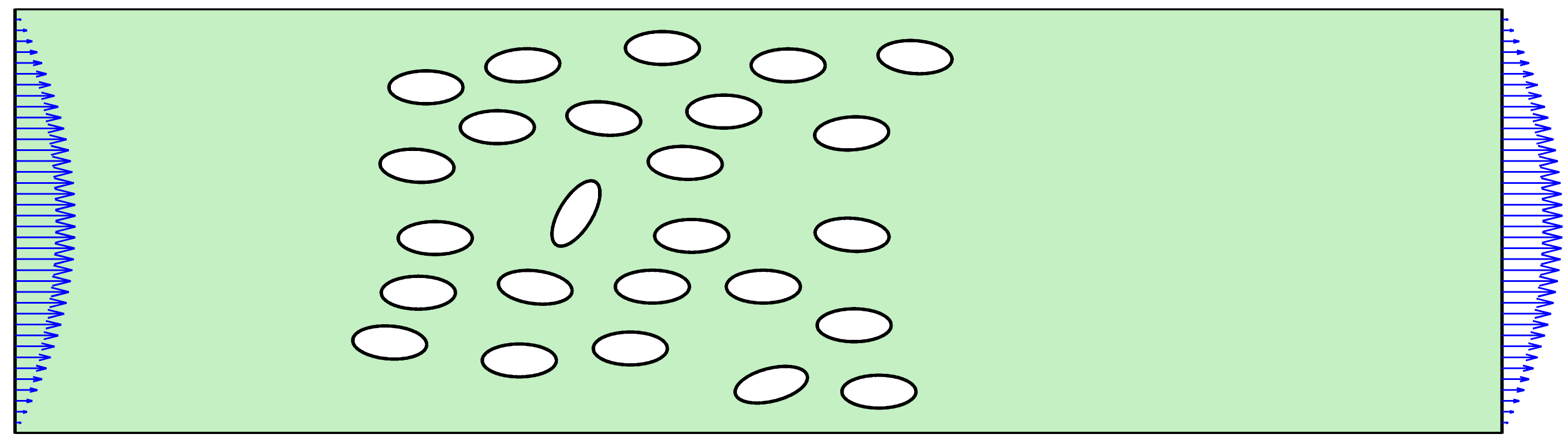}
	\caption{Computational domain and boundary conditions for the problem involving a flow in a channel with randomly distributed ellipses.}
	\label{fig:domainEllipses}
\end{figure}

The first example, similar to a test case presented in~\cite{liu2008new}, considers the flow around a set of randomly distributed set of 25 ellipses in a channel. Dirichlet boundary conditions are considered in the whole domain corresponding to a parabolic velocity profile on the left (inflow) and right (outflow) boundaries and zero (no-slip) Dirichlet boundary condition on the top and bottom walls and on the boundary of the ellipses, as illustrated in Figure~\ref{fig:domainEllipses}.

A coarse mesh with 2,443 triangular elements is first considered. As no analytical solution is available, a reference solution is computed in a much finer mesh with 28,150 elements and by employing a degree of approximation $k=4$. This reference solution is used to measure the accuracy of the adaptive computations performed in much coarse meshes.

\begin{figure}[!tb]
	\centering
	\subfigure[Degree]{\includegraphics[width=0.8\textwidth]{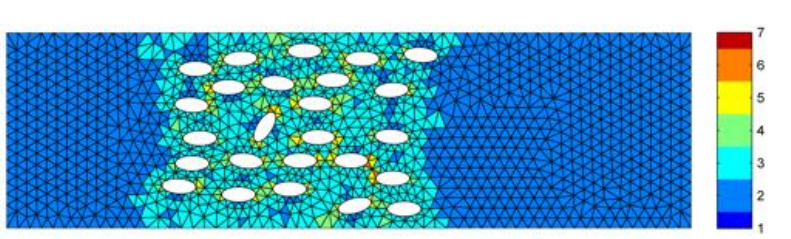}}
	\subfigure[Estimated error]{\includegraphics[width=0.8\textwidth]{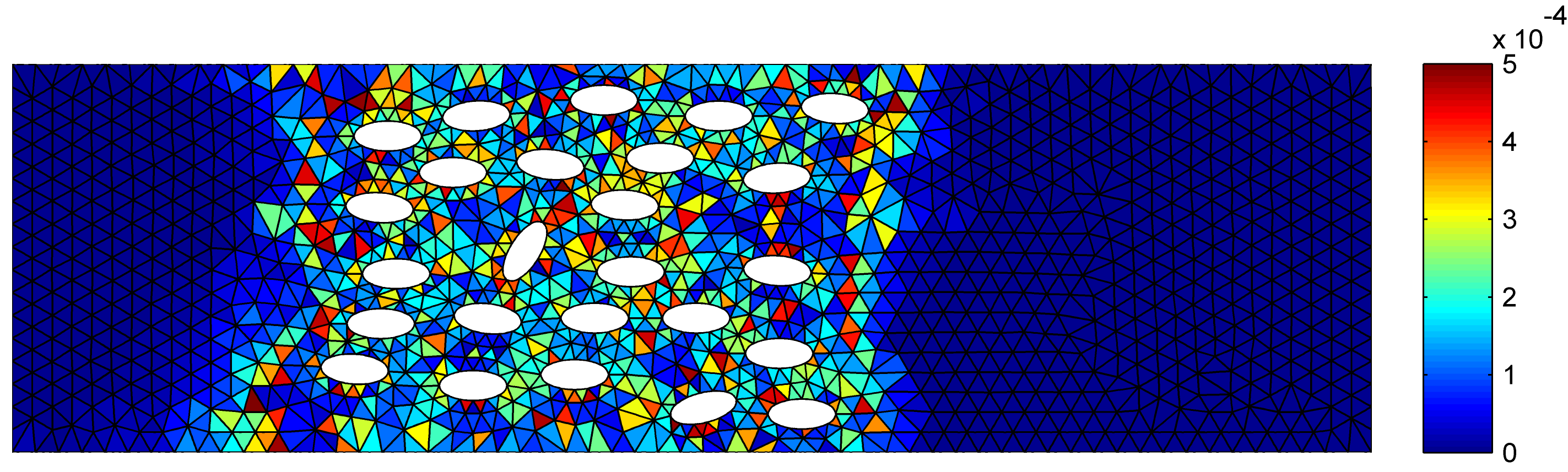}}
	\subfigure[Exact error]{\includegraphics[width=0.8\textwidth]{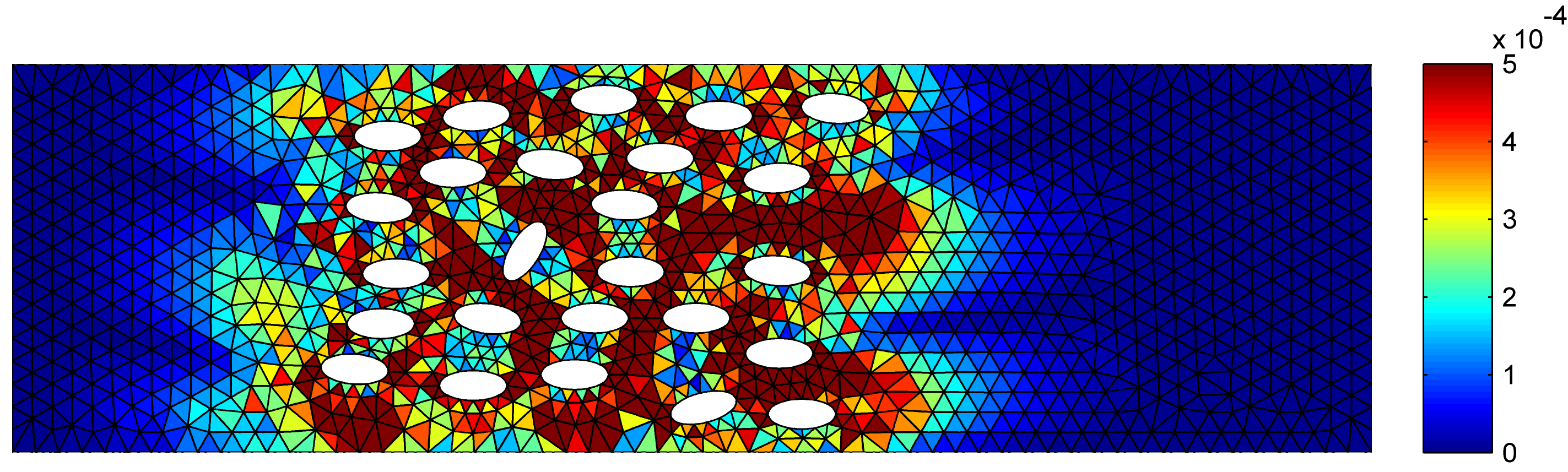}}
	\caption{Final degree distribution, estimated and reference errors for an adaptive computation with isoparametric HDG and quadratic approximation of the curved boundaries in a coarse mesh with 2,443 elements.}
	\label{fig:ellipsesFEMH1p2}
\end{figure}

An adaptive process is performed using a quadratic approximation of the curved boundaries and standard isoparametric elements with a desired error of $0.5 \times 10^{-3}$. Figure~\ref{fig:ellipsesFEMH1p2} (a) shows the computational mesh and degree of approximation after five iterations.
In the vicinity of the ellipses the majority of elements have a cubic degree of approximation where the elements in contact with the ellipses need a higher order of approximation to capture all the flow features. The highest order of approximation is $k=6$, used, as expected, with the regions with higher curvature of the boundary. Figure~\ref{fig:ellipsesFEMH1p2} (b) and (c) show the estimated and the reference errors after the adaptive process converged. The discrepancy between the estimated and the reference errors is clearly observed. Despite the adaptive process converges, meaning that all elements have an estimated error below the desired error, a total of 408 elements have a reference error above the desired tolerance. 

When the same computation is performed by considering a cubic approximation of the geometry (not reported for brevity), the adaptive process converges again in five iterations. he highest order of approximation used in a few elements is now $k=7$, indicating that a different geometric representation leads to a different degree of approximation required to achieve convergence. In addition, the error estimator is again not reliable as there are 15 elements where the reference error is above the desired error. 

\begin{figure}[!tb]
	\centering
	\subfigure[Degree]{\includegraphics[width=0.8\textwidth]{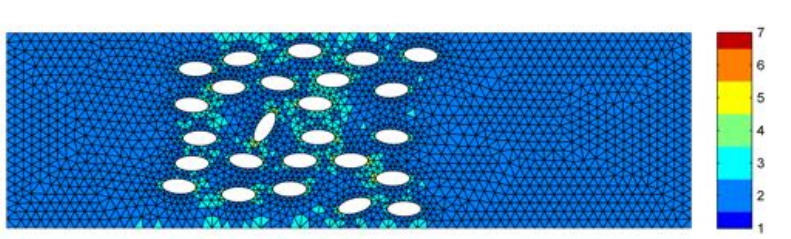}}
	\subfigure[Estimated error]{\includegraphics[width=0.8\textwidth]{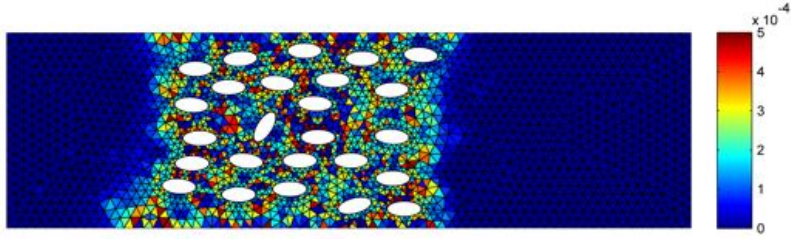}}
	\subfigure[Exact error]{\includegraphics[width=0.8\textwidth]{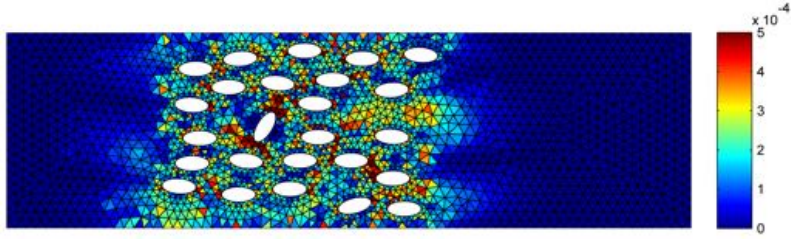}}
	\caption{Final degree distribution, estimated and reference errors for an adaptive computation with isoparametric HDG and quadratic approximation of the curved boundaries in a finer mesh with 4,048 elements.}
	\label{fig:ellipsesFEMH2p2}
\end{figure}

It is apparent that an adaptive computation with isoparametric elements requires an initial pre-adaptation of the mesh and the degree of approximation used to approximate the geometry in order to obtain a reliable error estimator. Next, the finer mesh with 4,048 triangular elements depicted in Figure~\ref{fig:ellipsesFEMH2p2} (a) is considered. When a quadratic approximation of the geometry is considered, the adaptive process converges in four iterations but the results are still not satisfactory as there are 155 elements where the estimated error, shown in Figure~\ref{fig:ellipsesFEMH2p2} (b), is above the desired error, represented in Figure~\ref{fig:ellipsesFEMH2p2} (c). 

Finally, if a cubic approximation of the geometry is employed the adaptive process converges in only three iterations with one element still showing an error above the desired tolerance. The results clearly indicate that in the presence of curved boundaries the level of pre-adaptation required negates all the advantages of an an automatic adaptive process as the initial mesh has to be designed in such a way so that the error in the first computation is near the desired error. 

\begin{figure}[!tb]
	\centering
	\subfigure[Degree]{\includegraphics[width=0.8\textwidth]{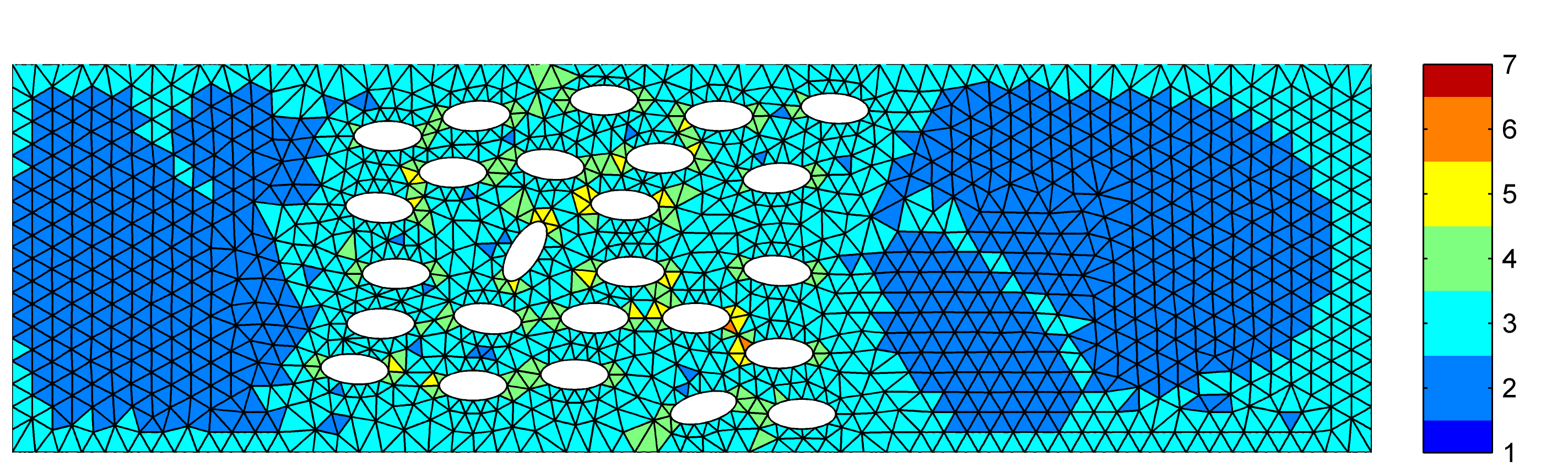}}
	\subfigure[Estimated error]{\includegraphics[width=0.8\textwidth]{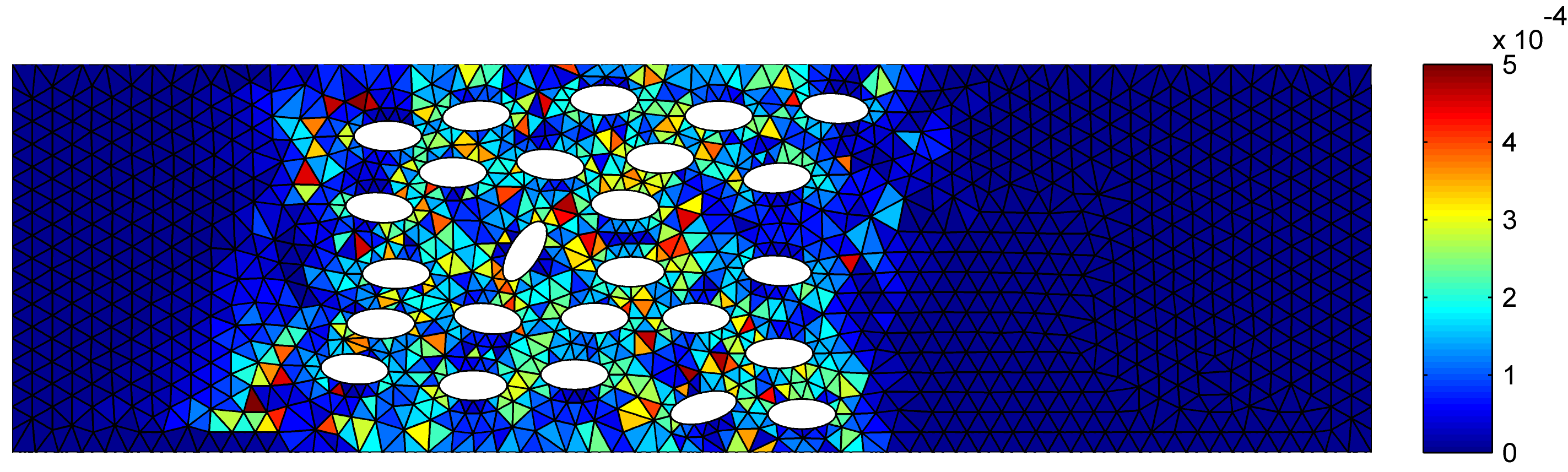}}
	\subfigure[Exact error]{\includegraphics[width=0.8\textwidth]{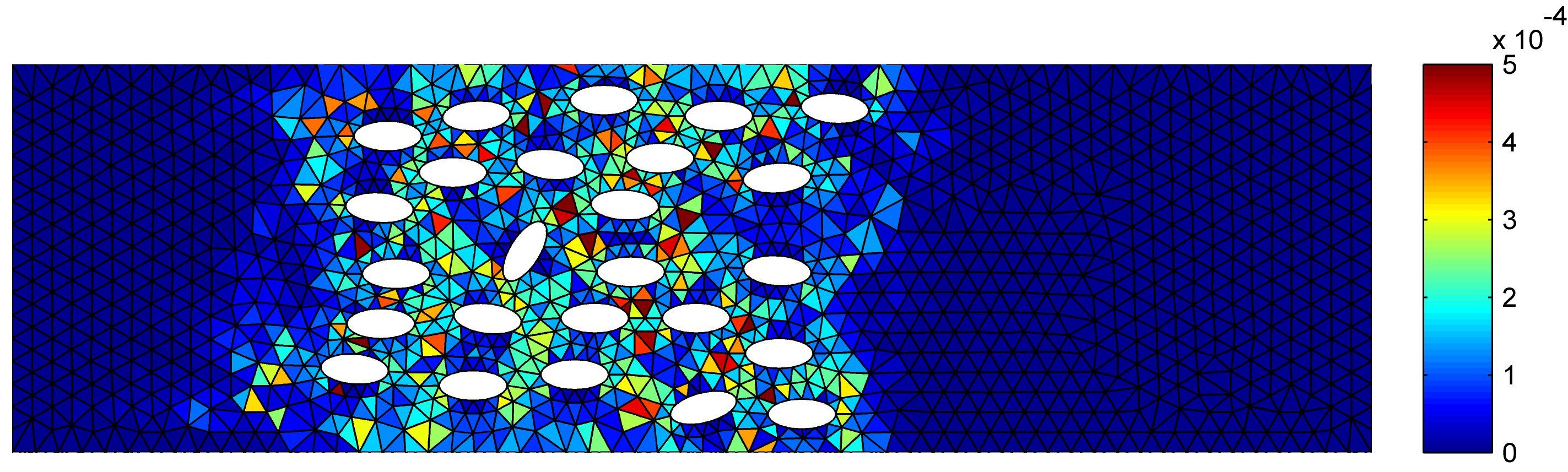}}
	\caption{Final degree distribution, estimated and reference errors for an adaptive computation with HDG-NEFEM in a coarse mesh with 2,443 elements.}
	\label{fig:ellipsesNEFEMH1p1}
\end{figure}

To show the potential of NEFEM in this scenario, an adaptive process is performed employing the coarse mesh with 2,443 triangular elements and starting with a degree of approximation $k=1$. The adaptive process converges in four iterations. The final degree of approximation used in each element is shown in Figure~\ref{fig:ellipsesNEFEMH1p1} (a), with two elements having the maximum degree of approximation, $k=6$, required to achieve the desired error. 
The estimated and reference errors, depicted in Figures~\ref{fig:ellipsesNEFEMH1p1} (b) and (c), respectively, shows a consistent behaviour that illustrates the reliability of the proposed strategy to estimate the error due to the use of the exact boundary representation. It is worth noting that in the majority of the elements surrounding the ellipses a degree of approximation $k=3$ is enough to obtain the required error, illustrating why cubic isoparametric elements outperformed the use of quadratic elements in the previous computations.

\begin{figure}[!tb]
	\centering
	\includegraphics[width=0.8\textwidth]{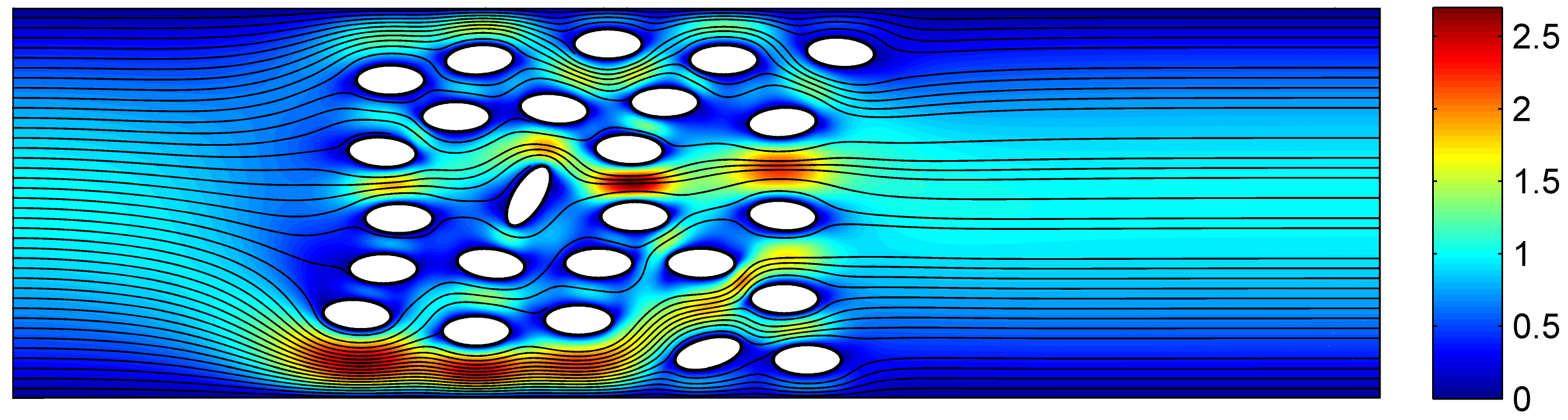}
	\caption{Magnitude of velocity and isolines computed with HDG-NEFEM on the mesh shown in Figure~\ref{fig:ellipsesNEFEMH1p1} (a) after four iterations of the degree adaptive process.}
	\label{fig:velocityNEFEMH1P1}
\end{figure}

The velocity field computed with NEFEM on the mesh shown in Figure~\ref{fig:ellipsesNEFEMH1p1} (a) is depicted in Figure~\ref{fig:velocityNEFEMH1P1}.

\begin{figure}[!htb]	
	\centering
	\subfigure[$q$=2]{\includegraphics[width=0.49\textwidth]{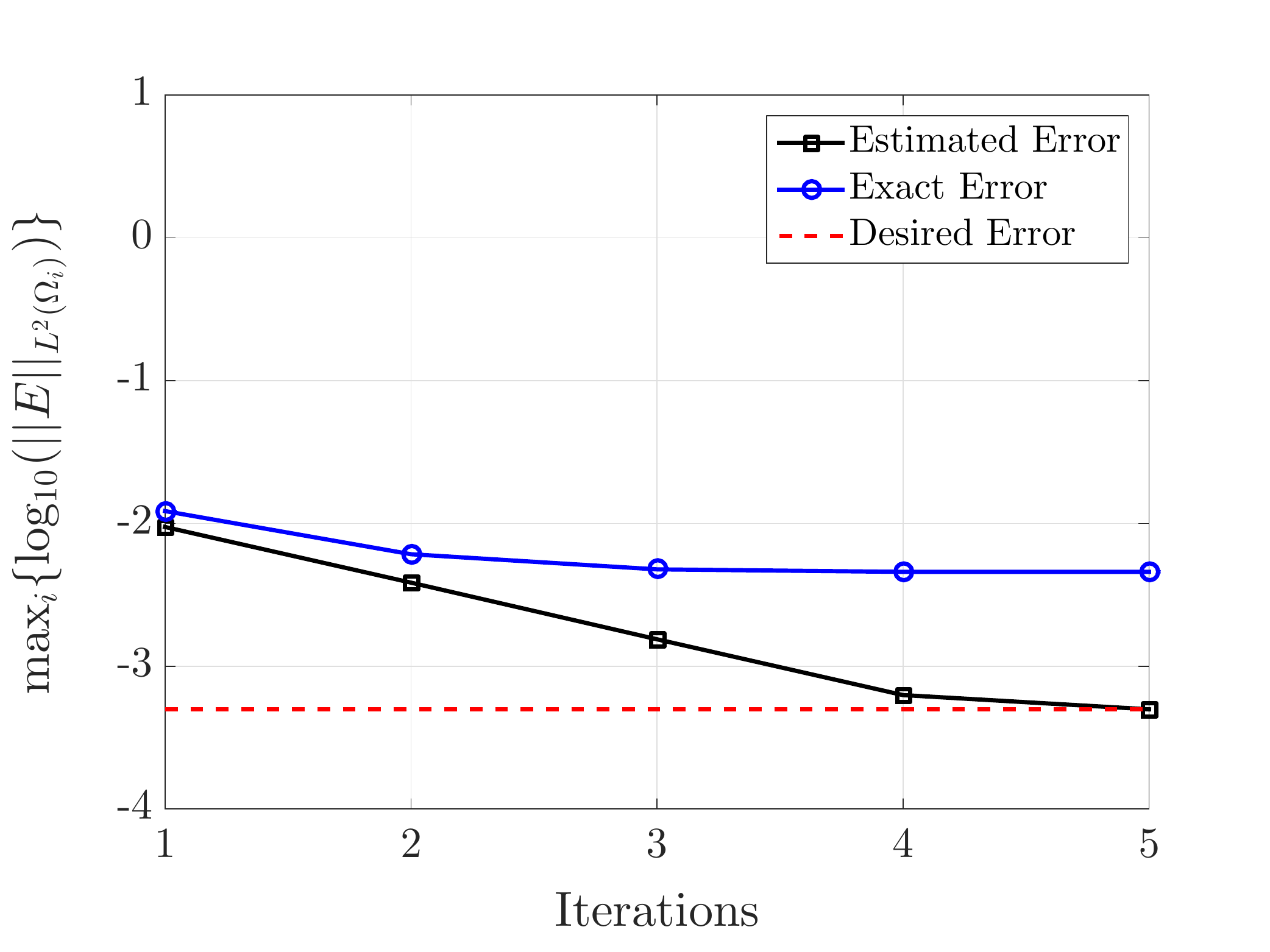}}
	\subfigure[$q$=3]{\includegraphics[width=0.49\textwidth]{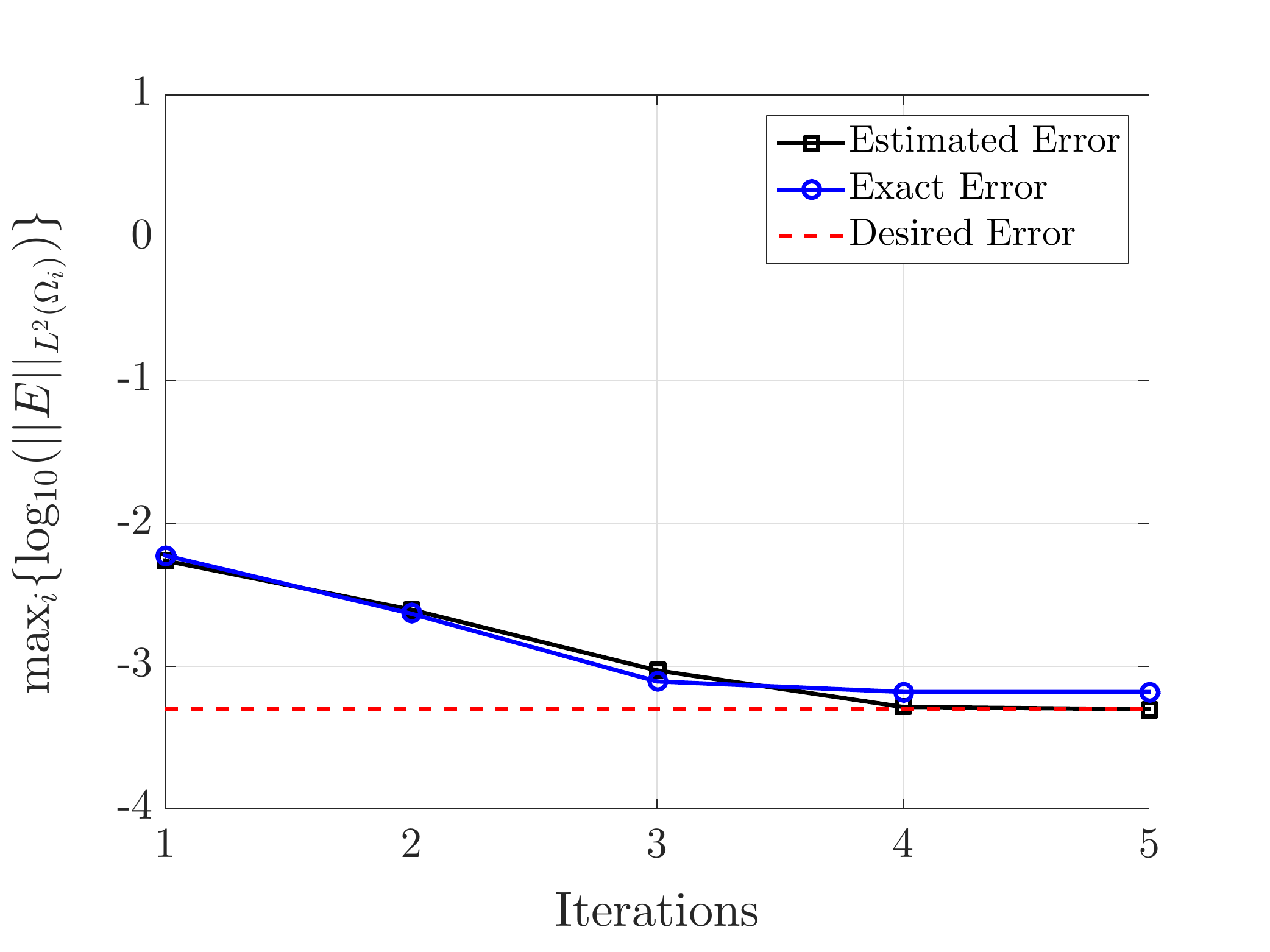}}
	\caption{Evolution of the estimated and exact errors during the degree adaptivity process for $q=2$ and $q=3$ with isoparametric HDG in the coarse mesh with 2,443 elements and a desired error of $0.5 \times 10^{-3}$.}
	\label{fig:adaptivityErrorEllipsesFEMH1Em3}
\end{figure}

To better analyse the effect of an accurate boundary representation in a degree adaptive procedure, Figure~\ref{fig:adaptivityErrorEllipsesFEMH1Em3} shows the evolution of the estimated and exact errors as a function of the number of iterations of the adaptive process by using the coarse mesh with 2,443 elements.
For a desired error of $0.5 \times 10^{-3}$ a quadratic approximation of the geometry prevents convergence of the exact error whereas better results are obtained with a cubic representation of the geometry. It is worth noting that in both cases the exact error stagnates, indicating that the geometric error dominates over the interpolation error. Even when a cubic representation of the geometry is considered, the exact error is not decreasing after the third iteration.

\begin{figure}[!htb]	
	\centering
	\subfigure[$q$=2]{\includegraphics[width=0.49\textwidth]{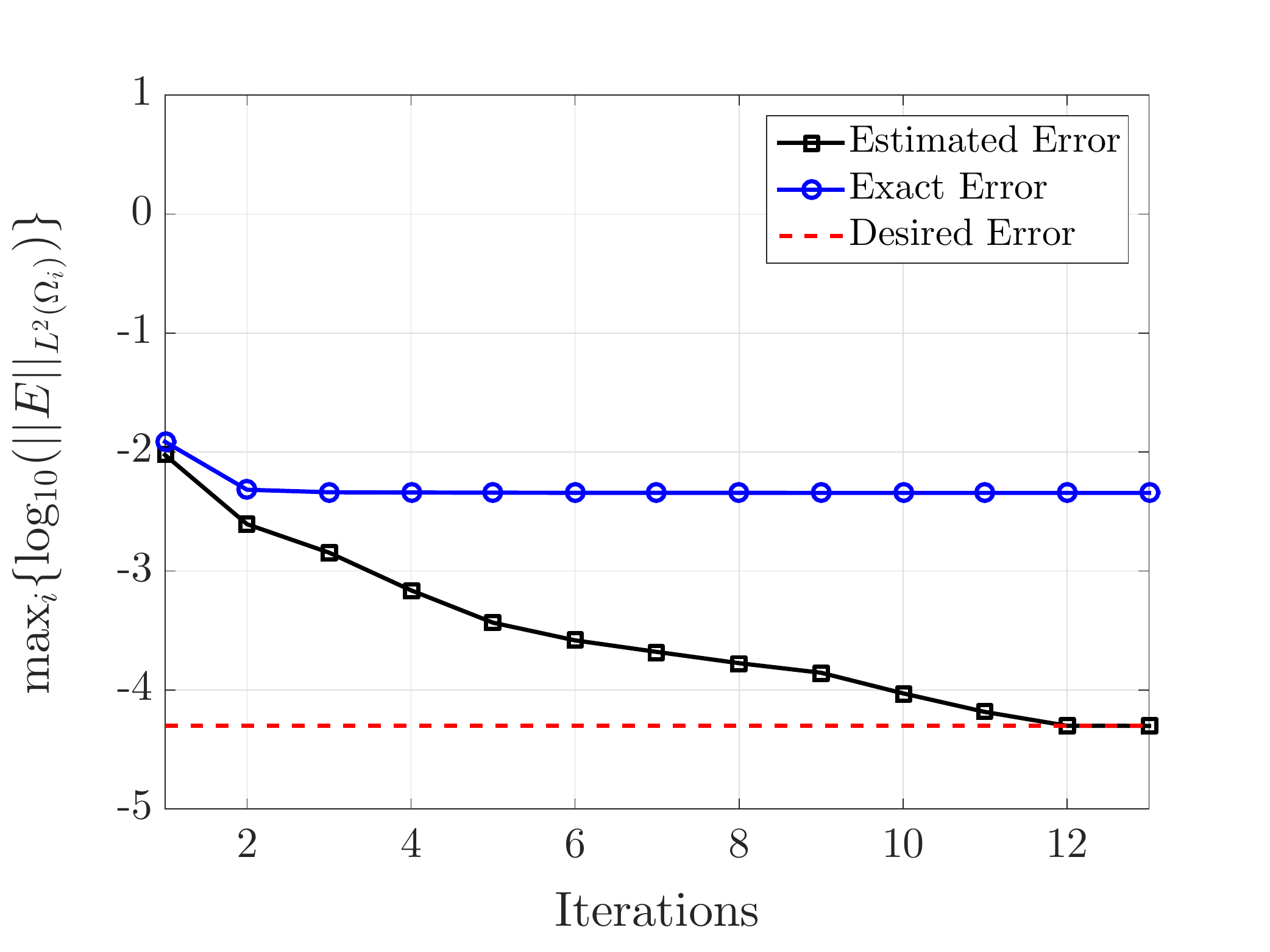}}
	\subfigure[$q$=3]{\includegraphics[width=0.49\textwidth]{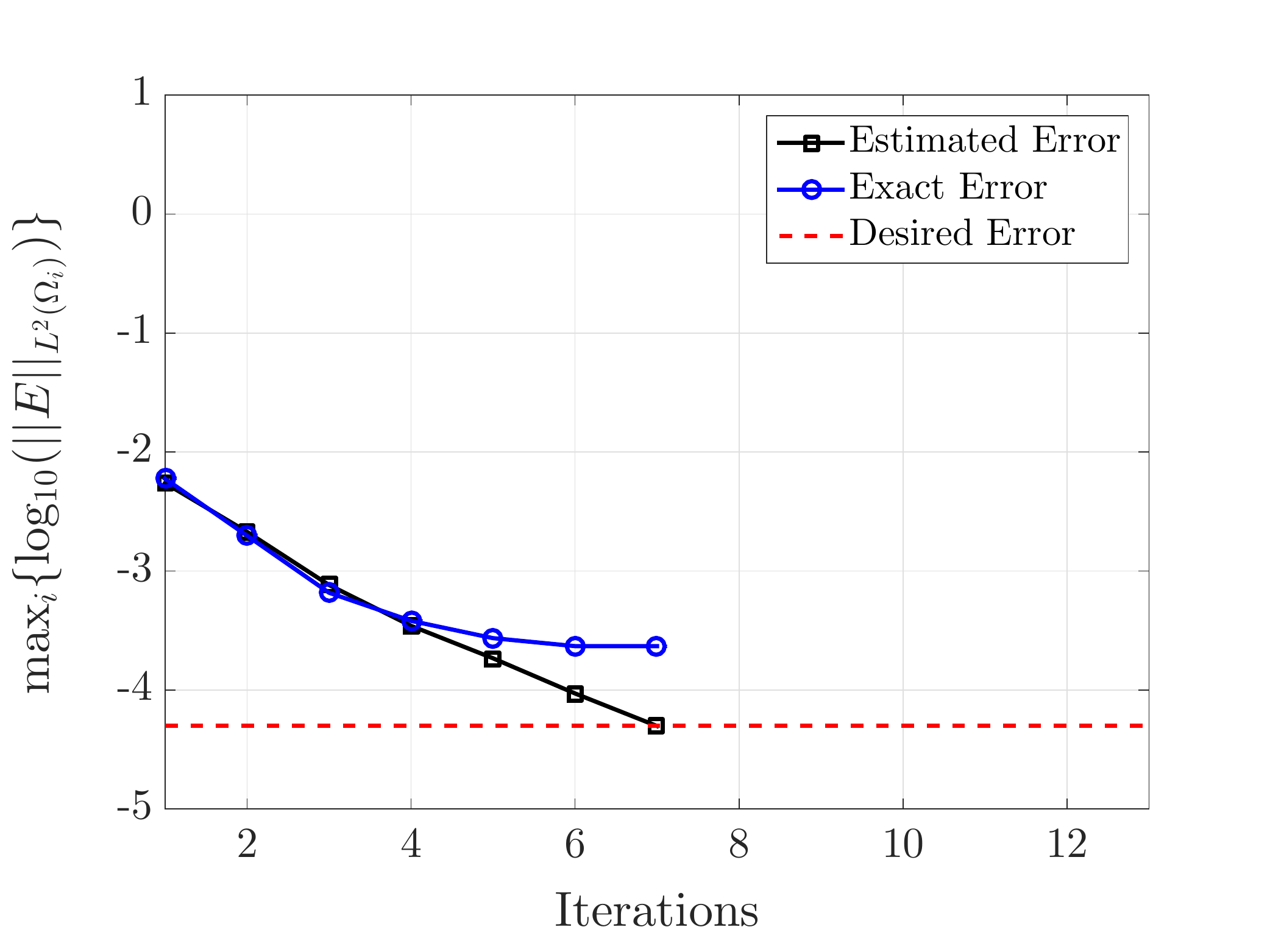}}
	\caption{Evolution of the estimated and exact errors during the degree adaptivity process for $q=2$ and $q=3$ with isoparametric HDG in the coarse mesh with 2,443 elements and a desired error of $0.5 \times 10^{-4}$.}
	\label{fig:adaptivityErrorEllipsesFEMH1Em4}
\end{figure}

To further illustrate the limitations of an isoparametric formulation during a degree adaptivity procedure, the same analysis is repeated with a lower desired error, namely $0.5 \times 10^{-4}$. Figure~\ref{fig:adaptivityErrorEllipsesFEMH1Em4} shows the evolution of the estimated and exact errors as a function of the number of iterations of the adaptive process by using the coarse mesh with 2,443 elements. 
The results show that a cubic approximation of the geometry is not enough because the difference between the estimated and exact error in the final computation with cubic elements is almost an order of magnitude.  This suggests, once more that the initial mesh has to be pre-adapted to achieve a reliable degree adaptive process. 

\begin{figure}[!htb]	
	\centering
	\subfigure[Desired error $0.5 \times 10^{-3}$]{\includegraphics[width=0.49\textwidth]{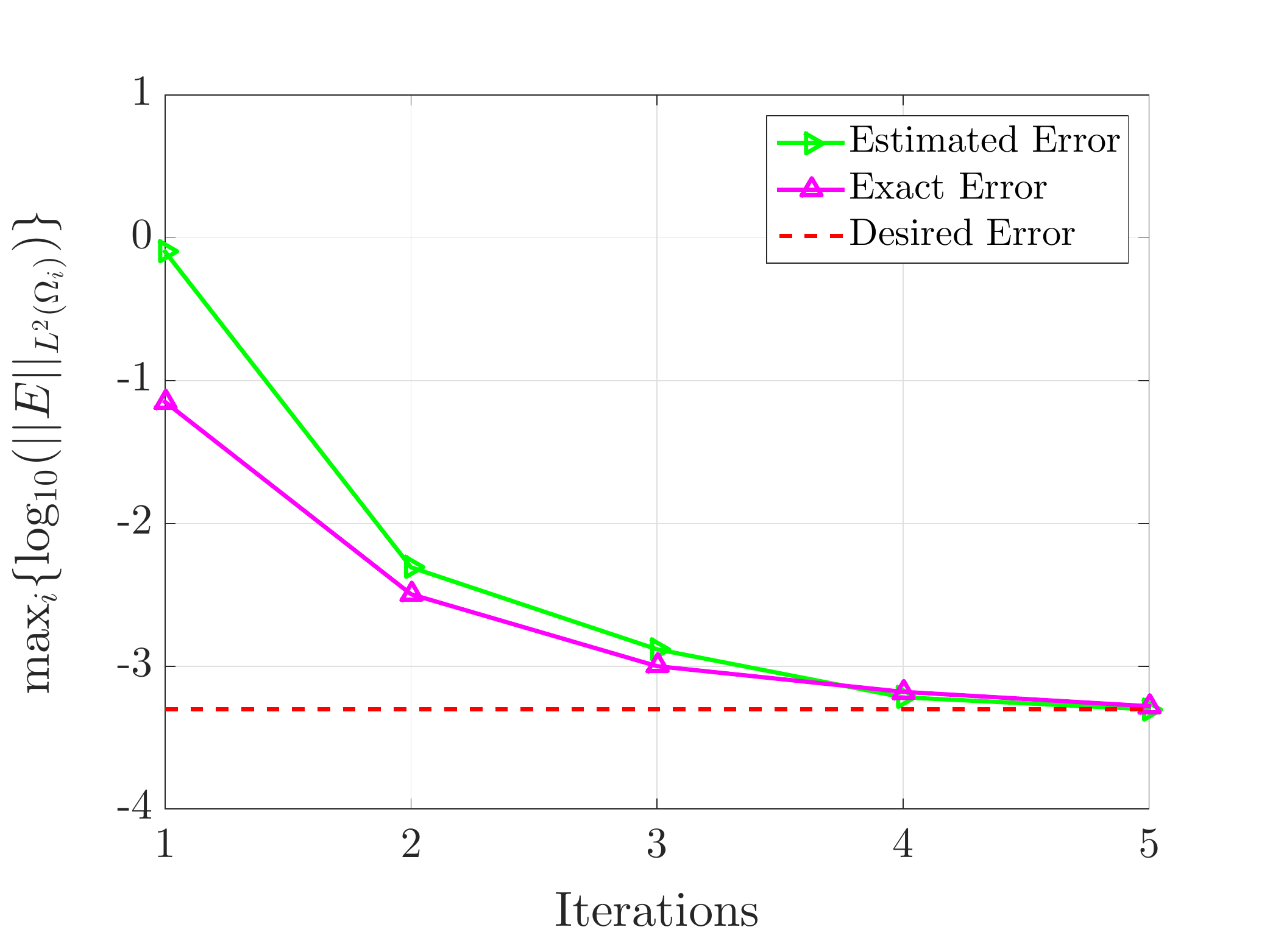}}
	\subfigure[Desired error $0.5 \times 10^{-4}$]{\includegraphics[width=0.49\textwidth]{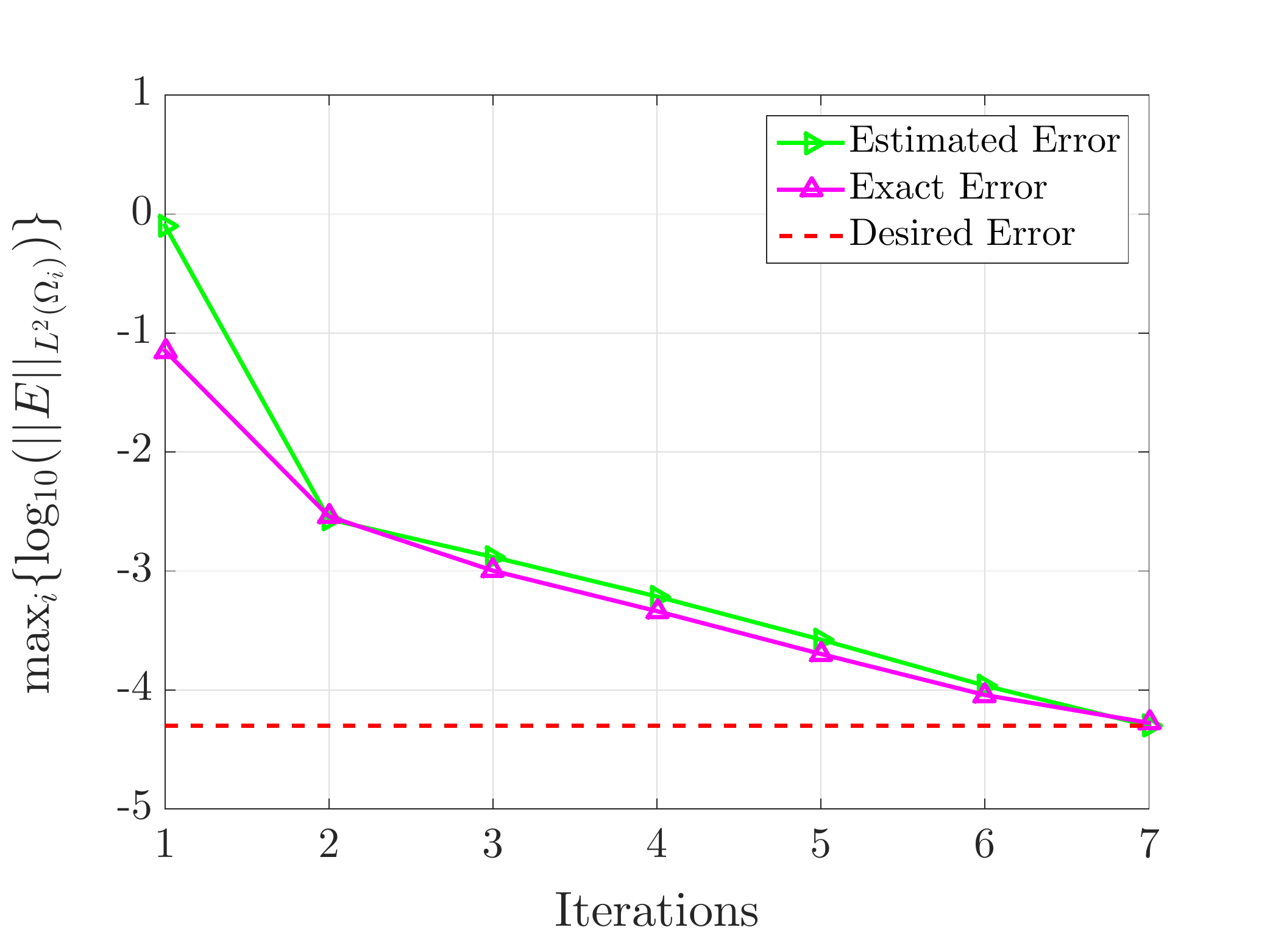}}
	\caption{Evolution of the estimated and exact errors during the degree adaptivity process for NEFEM HDG in the coarse mesh with 2,443 elements.}
	\label{fig:adaptivityErrorEllipsesNEFEM}
\end{figure}

Finally, Figure~\ref{fig:adaptivityErrorEllipsesNEFEM} shows the results obtained with NEFEM using the coarse mesh with 2,443 elements and starting with a linear approximation of the solution $k=1$.
The robustness of the proposed approach is clearly illustrated as convergence of both the estimated and exact errors is achieved in the coarse mesh even for a desired error of $0.5 \times 10^{-4}$. It is worth emphasising that with NEFEM the adaptive process provides a reliable error estimator even when the desired error is several orders of magnitude lower than the error of the computation in the first mesh. The results clearly indicate that no pre-adaptation of the mesh is required with NEFEM as the geometry is exactly represented irrespectively of the spatial discretisation. Therefore, the adaptive process is purely driven by the functional approximation and not by the geometric error as it happens with an isoparametric formulation.  

Further numerical experiments, not reported here, indicate that NEFEM is also superior to an isoparametric approach where the mesh is  re-generated at each iteration of the adaptive procedure by using the same degree of the approximation for both the geometry and the solution. In all cases, not only the time required by NEFEM is lower (due to the extra time required to communicate with the CAD model and the mesh generator) but also due to the fact that more iterations of the adaptive process are required with an isoparametric formulation.

%--------------------------------------------------------------------------
\subsection{Flow in a channel with wavy boundaries}  \label{sc:exampleSection}
%--------------------------------------------------------------------------

\begin{figure}[!htb]	
	\centering
	\subfigure[On phase oscillations]{\includegraphics[width=0.8\textwidth]{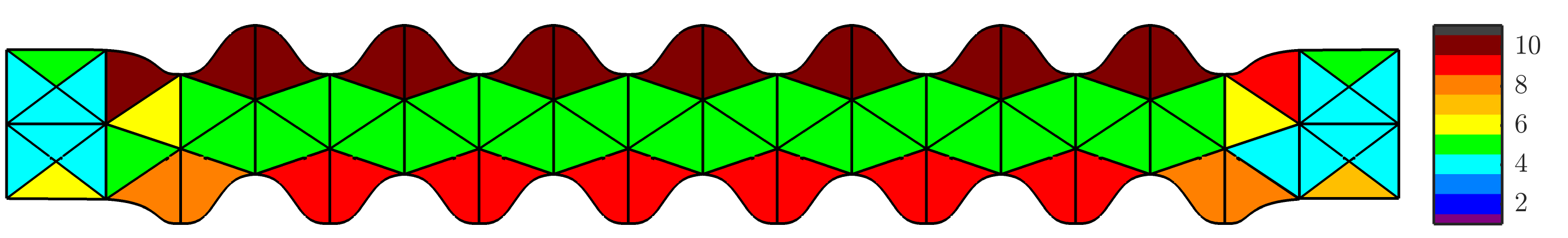}}
	\subfigure[Out of phase oscillations]{\includegraphics[width=0.8\textwidth]{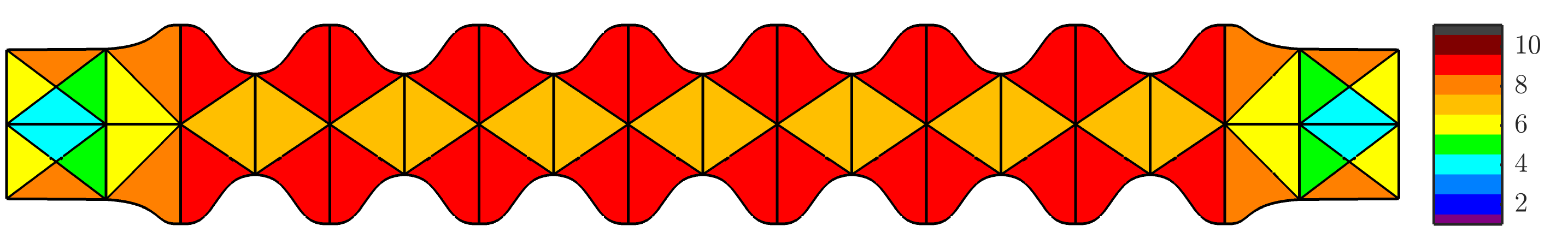}}
	\caption{Mesh and degree of approximation of the converged degree adaptive procedure with HDG-NEFEM for the computation of the flow in a channel with oscillations of the top and bottom boundaries.}
	\label{fig:waveChannel1P}
\end{figure}

The next example considers the flow in a channel with oscillatory boundaries. This problem, of interest to the micro and nano-fluidics community, is often considered to study the flow structure induced by the different phase of the oscillations of the top and bottom boundaries~\cite{khuri2006stokes,wang2011stokes}. The two extreme cases are considered here, where the oscillations are exactly on phase and completely out of phase. Figure~\ref{fig:waveChannel1P} (a) shows the coarse computational mesh employed with HDG-NEFEM and the final degree of approximation obtained after ten iterations of the degree adaptive process for the case where the boundary oscillations are on phase. Similarly, Figure~\ref{fig:waveChannel1P} (b) shows the mesh and degree of approximation obtained after eight iterations of the degree adaptive process for the case where the boundary oscillations are completely out of phase.

\begin{figure}[!htb]	
	\centering
	\subfigure[On phase oscillations]{\includegraphics[width=0.8\textwidth]{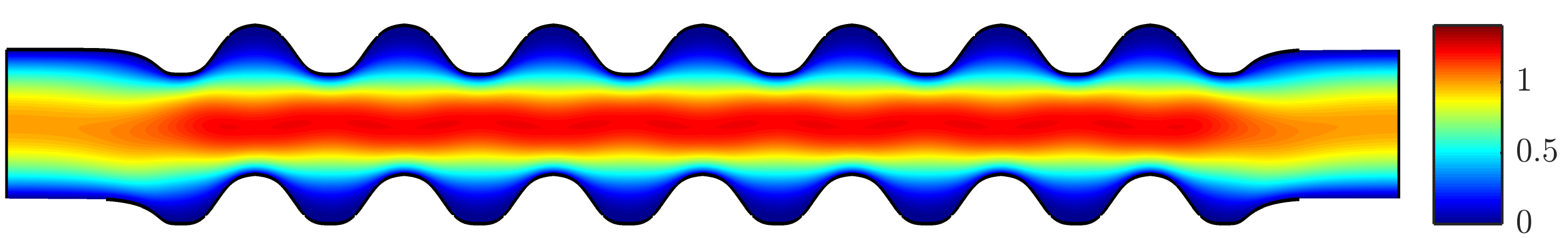}}
	\subfigure[Out of phase oscillations]{\includegraphics[width=0.8\textwidth]{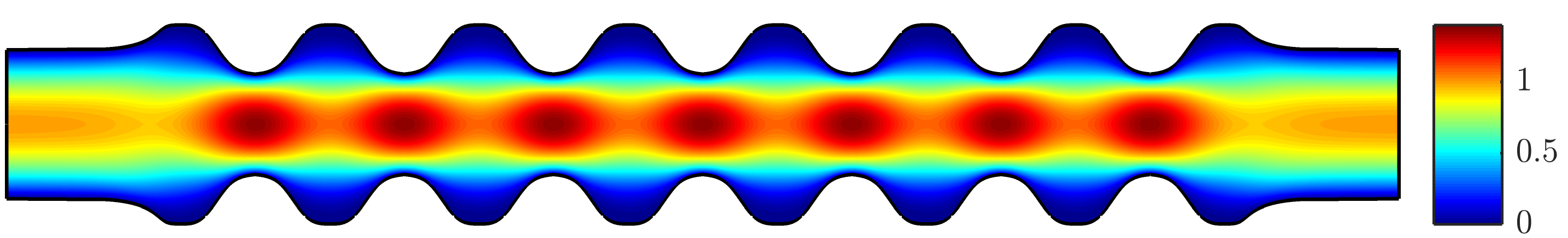}}
	\caption{Magnitude of velocity computed with HDG-NEFEM after convergence of the degree adaptive procedure on the meshes shown in Figure~\ref{fig:waveChannel1P}.}
	\label{fig:waveChannelVelo}
\end{figure}

The velocity fields obtained for both cases are represented in Figure~\ref{fig:waveChannelVelo} showing the ability of the proposed approach to capture the different flow structure induced by the oscillatory boundaries.

%--------------------------------------------------------------------------
\subsection{Flow in a porous media}  \label{sc:examplePorous}
%--------------------------------------------------------------------------

\begin{figure}[!htb]	
	\centering
	\subfigure[Degree]{\includegraphics[width=0.49\textwidth]{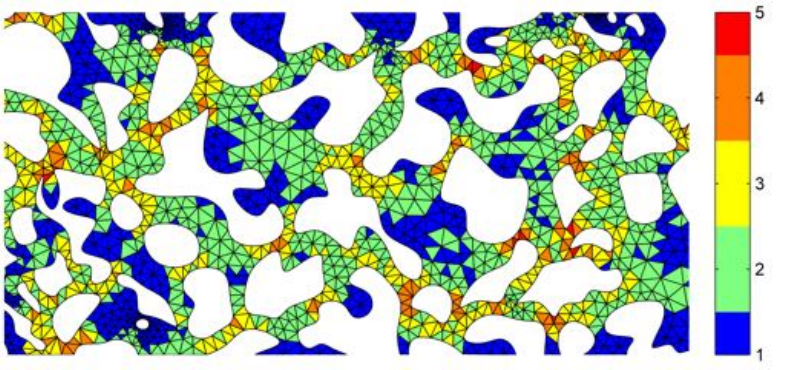}}
	\subfigure[Velocity]{\includegraphics[width=0.49\textwidth]{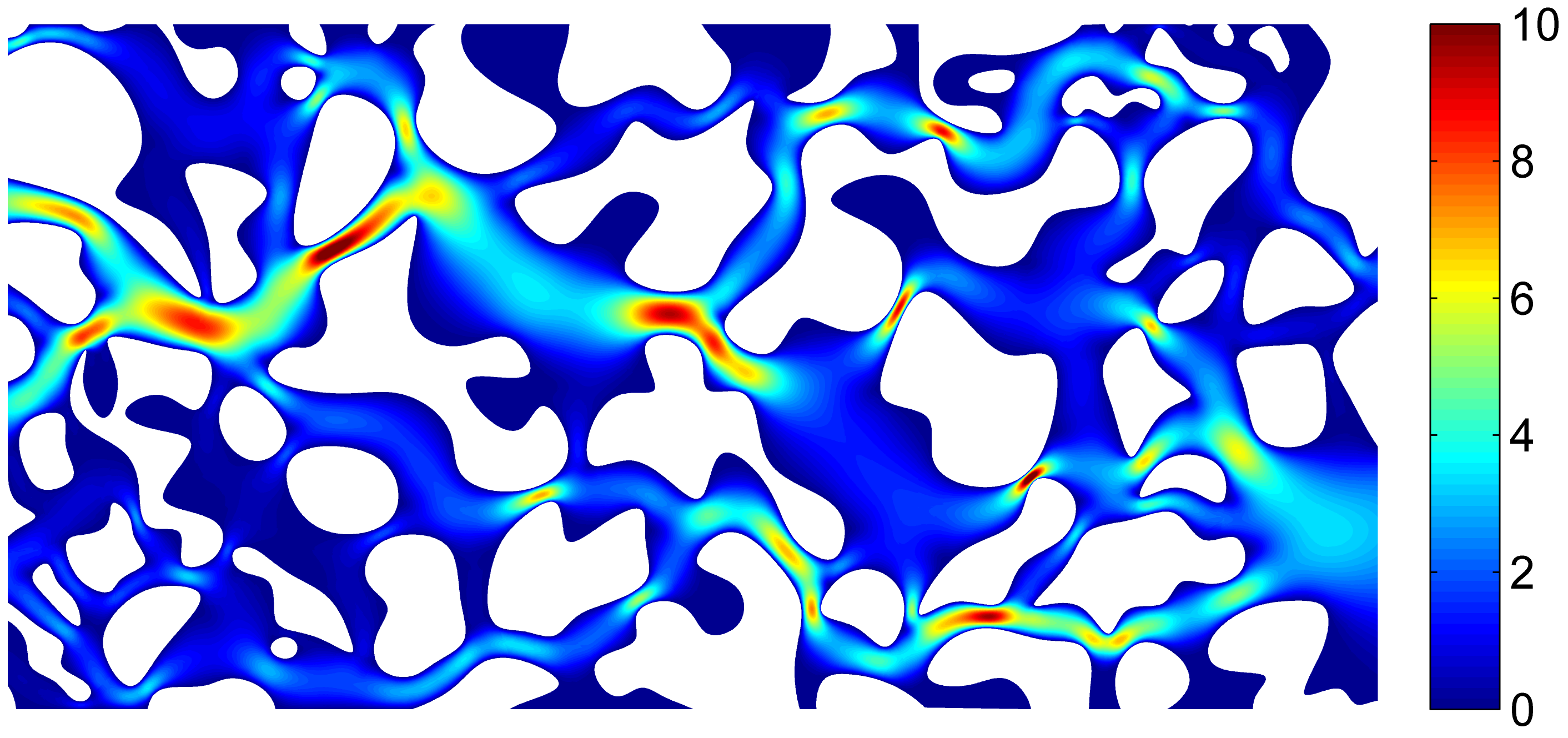}}
	\caption{Mesh and degree of approximation of the converged degree adaptive procedure with HDG-NEFEM and velocity field.}
	\label{fig:pore}
\end{figure}

The next example, taken from~\cite{sirivithayapakorn2003transport}, considers the flow in the interstices of a porous media. The geometry consists of the surroundings of a large number of particles in the porous media. Figure~\ref{fig:pore} shows the mesh and degree of approximation after eight iterations of the degree adaptivity procedure.

It is worth noting that a linear degree of approximation is used in many elements in contact with curved boundaries. This shows that the proposed adaptivity strategy is completely driven by the complexity of the solution and not by the complexity of the geometry.

%--------------------------------------------------------------------------
\subsection{Flow in a channel with thin obstacles}  \label{sc:exampleObstacles}
%--------------------------------------------------------------------------

The last example shows a further benefit of using NEFEM by demonstrating its unique ability to obtain accurate solutions with ultra-coarse meshes even when geometric features, smaller than the element size, are present in the boundary representation of the computational domain.

\begin{figure}[!htb]	
	\centering
	\subfigure[Whole domain]{\includegraphics[width=\textwidth]{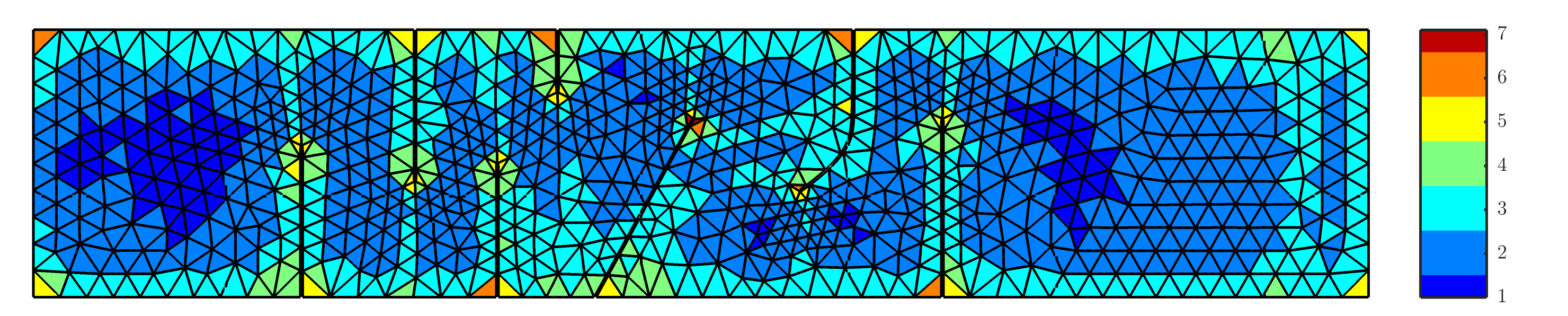}}
	\subfigure[Zoom]{\includegraphics[width=\textwidth]{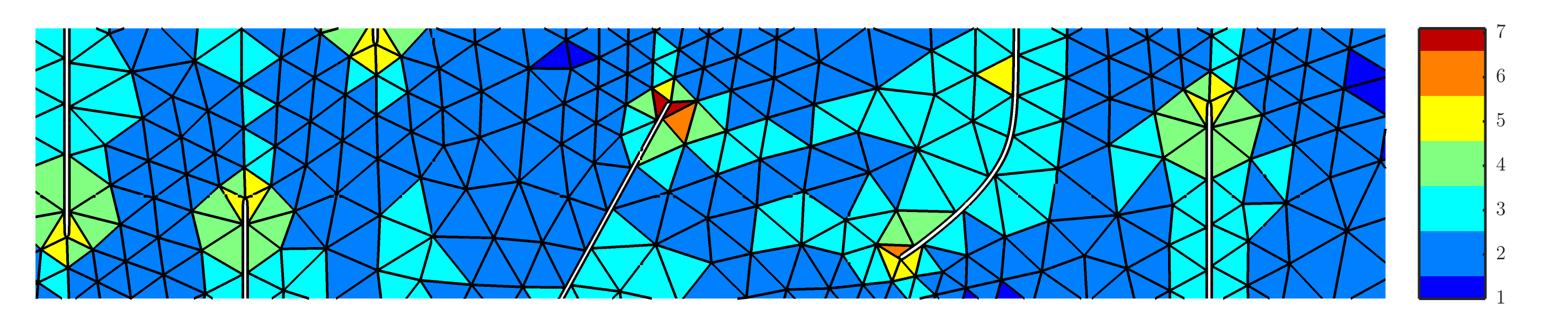}}
	\caption{Mesh and degree of approximation of the converged degree adaptive procedure with HDG-NEFEM for the computation of the flow in a channel with high curvature.}
	\label{fig:obstacles_P}
\end{figure}

The flow in a channel with a number of thin obstacles is considered. The thickness of the obstacles is approximately 0.08 whereas the minimum element size of the mesh that has been generated, using the technique proposed in~\cite{sevilla2016generation}, is 0.32. The degree adaptive process is started, as in previous examples, with a linear approximation of the solution and convergence is achieved in four iterations. The final degree of approximation in each element is represented in Figure~\ref{fig:obstacles_P} (a).
Figure~\ref{fig:obstacles_P} shows a detailed view of a region in the channel showing the elements near the end of some of the obstacles. This plot shows not only that the element size is independent on the geometric complexity but it also demonstrates the robustness of the proposed degree adaptive technique. The adaptivity process is clearly driven only by the complexity of the solution as a different degree of approximation is employed in elements with almost identical geometric complexity due to the different complexity of the solution.

\begin{figure}[!htb]	
	\centering
	\includegraphics[width=0.8\textwidth]{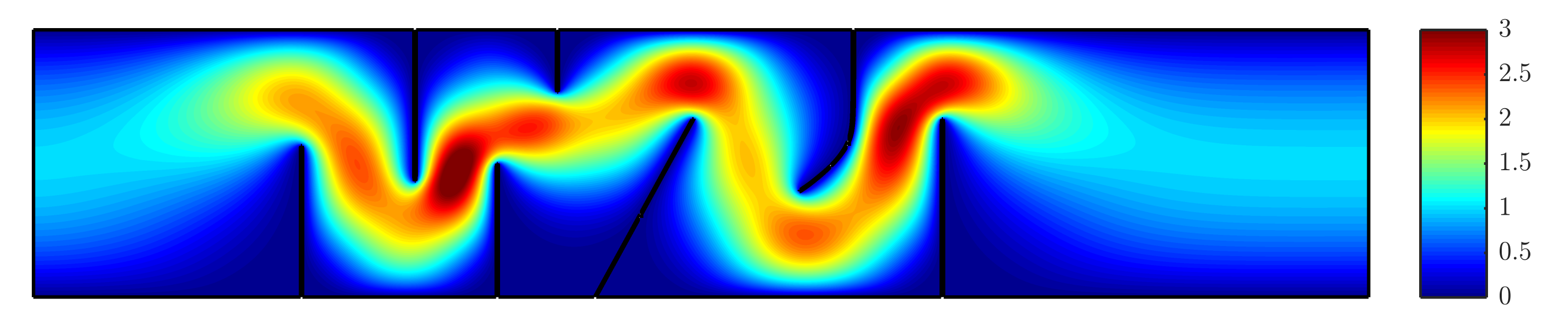}
	\caption{Velocity computed with HDG-NEFEM after convergence of the degree adaptive procedure in the mesh of Figure~\ref{fig:obstacles_P}.}
	\label{fig:obstaclesVelo}
\end{figure}

The velocity field, obtained on the mesh shown in Figure~\ref{fig:obstacles_P}, is represented in Figure~\ref{fig:obstaclesVelo}.

%==========================================================================
\section{Concluding remarks} \label{sc:conclusions}
%==========================================================================

A new degree adaptive methodology that combines the advantages of the HDG formulation and NEFEM has been presented. The proposed method results in a cheap and reliable error estimator due to the cheap computation of a post-processed solution provided by HDG and the ability to exactly represent curved boundaries irrespectively of the polynomial degree used for the functional approximation that is characteristic of NEFEM.

The proposed approach is compared against two alternative options to perform degree adaptivity. The first approach, broadly used in practice, consists of keeping the shape of the curved elements during the degree adaptivity process. It is found that this approach leads to an unreliable error estimator. The numerical examples show that even when the estimated error is below the required tolerance, the exact error can be orders of magnitude higher. The second approach, not used in practice, consists of changing the shape of the curved elements during the adaptive process. The main drawback is its high cost due to the need to constantly communicate with the CAD model and re-generate the nodal distributions for curved elements. 

The proposed approach considers, for the first time, the implementation of the NEFEM rationale in an HDG framework.  A number of numerical examples have been presented to compare the performance of the proposed methodology and to shows its superiority on a number of problems involving domains with curved boundaries.

%==============================================================
\section*{Acknowledgements} 
%==============================================================

The authors gratefully acknowledge the financial support of the Ministerio de Econom\'ia y Competitividad (grant number DPI2014-51844-C2-2-R). The first author also gratefully acknowledges the financial support provided by the S\^{e}r Cymru National Research Network for Advanced Engineering and Materials (grant number NRN045). The second author is also grateful for the financial support provided by the Generalitat de Catalunya (grant number 2014-SGR-1471).

%________________________________________________________________________
\bibliographystyle{abbrv}
\bibliography{references}
\end{document}